\documentstyle[11pt]{article}

\newread\epsffilein    
\newif\ifepsffileok    
\newif\ifepsfbbfound   
\newif\ifepsfverbose   
\newif\ifepsfdraft     
\newdimen\epsfxsize    
\newdimen\epsfysize    
\newdimen\epsftsize    
\newdimen\epsfrsize    
\newdimen\epsftmp      
\newdimen\pspoints     
\def\epsfbox#1{\global\def\epsfllx{72}\global\def\epsflly{72}%
   \global\def\epsfurx{540}\global\def\epsfury{720}%
   \def\lbracket{[}\def\testit{#1}\ifx\testit\lbracket
   \let\next=\epsfgetlitbb\else\let\next=\epsfnormal\fi\next{#1}}%
\def\epsfgetlitbb#1#2 #3 #4 #5]#6{\epsfgrab #2 #3 #4 #5 .\\%
   \epsfsetgraph{#6}}%
\def\epsfnormal#1{\epsfgetbb{#1}\epsfsetgraph{#1}}%
\def\epsfgetbb#1{%
\openin\epsffilein=#1
\ifeof\epsffilein\errmessage{I couldn't open #1, will ignore it}\else
   {\epsffileoktrue \chardef\other=12
    \def\do##1{\catcode`##1=\other}\dospecials \catcode`\ =10
    \loop
       \read\epsffilein to \epsffileline
       \ifeof\epsffilein\epsffileokfalse\else
          \expandafter\epsfaux\epsffileline:. \\%
       \fi
   \ifepsffileok\repeat
   \ifepsfbbfound\else
    \ifepsfverbose\message{No bounding box comment in #1; using defaults}\fi\fi
   }\closein\epsffilein\fi}%
\def\epsfclipoff{\def\epsfclipstring{\ifepsfdraft\space clip\fi}}%
\epsfclipoff
\def\epsfsetgraph#1{%
   \epsfrsize=\epsfury\pspoints
   \advance\epsfrsize by-\epsflly\pspoints
   \epsftsize=\epsfurx\pspoints
   \advance\epsftsize by-\epsfllx\pspoints
   \epsfxsize\epsfsize\epsftsize\epsfrsize
   \ifnum\epsfxsize=0 \ifnum\epsfysize=0
      \epsfxsize=\epsftsize \epsfysize=\epsfrsize
      \epsfrsize=0pt
     \else\epsftmp=\epsftsize \divide\epsftmp\epsfrsize
       \epsfxsize=\epsfysize \multiply\epsfxsize\epsftmp
       \multiply\epsftmp\epsfrsize \advance\epsftsize-\epsftmp
       \epsftmp=\epsfysize
       \loop \advance\epsftsize\epsftsize \divide\epsftmp 2
       \ifnum\epsftmp>0
          \ifnum\epsftsize<\epsfrsize\else
             \advance\epsftsize-\epsfrsize \advance\epsfxsize\epsftmp \fi
       \repeat
       \epsfrsize=0pt
     \fi
   \else \ifnum\epsfysize=0
     \epsftmp=\epsfrsize \divide\epsftmp\epsftsize
     \epsfysize=\epsfxsize \multiply\epsfysize\epsftmp   
     \multiply\epsftmp\epsftsize \advance\epsfrsize-\epsftmp
     \epsftmp=\epsfxsize
     \loop \advance\epsfrsize\epsfrsize \divide\epsftmp 2
     \ifnum\epsftmp>0
        \ifnum\epsfrsize<\epsftsize\else
           \advance\epsfrsize-\epsftsize \advance\epsfysize\epsftmp \fi
     \repeat
     \epsfrsize=0pt
    \else
     \epsfrsize=\epsfysize
    \fi
   \fi
   \ifepsfverbose\message{#1: width=\the\epsfxsize, height=\the\epsfysize}\fi
   \epsftmp=10\epsfxsize \divide\epsftmp\pspoints
   \vbox to\epsfysize{\vfil\hbox to\epsfxsize{%
      \ifnum\epsfrsize=0\relax
        \includegraphics{\ifepsfdraft}%
      \else
        \epsfrsize=10\epsfysize \divide\epsfrsize\pspoints
        \includegraphics{\ifepsfdraft}%
      \fi
      \hfil}}%
\global\epsfxsize=0pt\global\epsfysize=0pt}%
{\catcode`\%=12 \global\let\epsfpercent=
\long\def\epsfaux#1#2:#3\\{\ifx#1\epsfpercent
   \def\testit{#2}\ifx\testit\epsfbblit
      \epsfgrab #3 . . . \\%
      \epsffileokfalse
      \global\epsfbbfoundtrue
   \fi\else\ifx#1\par\else\epsffileokfalse\fi\fi}%
\def\epsfempty{}%
\def\epsfgrab #1 #2 #3 #4 #5\\{%
\global\def\epsfllx{#1}\ifx\epsfllx\epsfempty
      \epsfgrab #2 #3 #4 #5 .\\\else
   \global\def\epsflly{#2}%
   \global\def\epsfurx{#3}\global\def\epsfury{#4}\fi}%
\def\epsfsize#1#2{\epsfxsize}
\newtheorem{sect}{}[section]

\begin{document}

\title{Computations of Quandle Cocycle Invariants of Knotted Curves and Surfaces}

\author{
J. Scott Carter
\\University of South Alabama \\
Mobile, AL 36688 \\ carter@mathstat.usouthal.edu
\and
Daniel Jelsovsky 
\\ University of South Florida
\\ Tampa, FL 33620  \\ jelsovsk@math.usf.edu
\and
Seiichi Kamada
\\  Osaka City University \\ 
Osaka 558-8585, JAPAN\\
kamada@sci.osaka-cu.ac.jp   
\and 
Masahico Saito
\\ University of South Florida
\\ Tampa, FL 33620  \\ saito@math.usf.edu
}

\maketitle

\begin{abstract}
State-sum invariants for knotted curves and surfaces using
quandle cohomology were introduced by Laurel Langford and the authors
in \cite{CJKLS}. In this paper we present
methods to compute the invariants and sample computations. 
Computer calculations of cohomological dimensions for some quandles
are presented.  For classical knots, Burau representations 
together with {\sc Maple} programs are used to evaluate the invariants for knot table. 
For knotted surfaces in $4$-space, movie methods and surface braid 
theory are used. Relations between the invariants and symmetries of
knots are discussed.

\vspace{1cm}

{\bf MRCN} : Primary 57M25, 57Q45; Secondary 55N99, 18G99.


\end{abstract}

\newpage

\section{Introduction}

In \cite{CJKLS}, we (with Laurel Langford) defined a state-sum invariant 
 of classical links and of knotted orientable surfaces.
 The invariant uses the cohomology theory of racks and quandles 
developed in \cite{FRS1,FRS2,Flower,Greene} 
 as its input.
 We modified the cohomology theory slightly to allow for type I
 Reidemeister moves and their higher dimensional analogues. 
Relations to linking numbers were given for some cocycles,
and it was shown that an invariant can detect non-invertibility
of the $2$-twist spun trefoil 
\cite{Rolf,Zeeman}. 
The nature of these invariants, however, is still a mystery.

The purpose of this paper is to present 
computational methods
in a variety of contexts. 
The computational results have topological implications, such as 
non-invertibility for some knotted surfaces.
For classical knots
we use Burau representations 
of the braid group and finite quotients of the 
Alexander quandles to give computations. 
For small quandles these are well suited to desktop computer calculations.
In the case of knotted surfaces in ${\bf R}^4$ we develop 
methods of computations using the theories 
of surface braids and  movies.
Our results are a combination of 
the above mentioned theories and 
computer calculations. 
The latter are supported by {\sc Maple} and {\sc Mathematica}. 
Here we have concentrated on several important families of knots 
and knotted surfaces.
In the classical case, we have computed several invariants
 in the knot  table 
up to 9 crossings and torus knots.
 In the knotted surface case we have calculated for 
twist-spun torus knots 
(where the movies and surface braids follow some patterns),
 and for the movie of  
a  deform-spun figure-8 knot. 
There 
are advantages 
to both the movie and the surface braid methods.

\begin{sect}{\bf Organization.\ }
{\rm 
In Section~\ref{revsec} the basic definitions are reviewed from \cite{CJKLS}.
Section~\ref{cocysec} presents 
 cohomological dimensions  
for some Alexander quandles, after reviewing
other calculations.
Using some of these cocycles, invariants for knots
in the table are computed in Section~\ref{tablesec}. 
Another good family of classical knots
is torus knots. We prove some periodicity and computations of invariants
for torus knots in Section~\ref{torussec}.
In Section~\ref{defate}, we give an explicit computation for a 
deform-spun figure-8 knot. 
This knotted sphere has 6 critical points, 16 triple points and no branch points. 
For twist spun torus knots in dimension $4$, 
we use movie methods (Section~\ref{moviesec}) 
and surface braid theory (Section~\ref{sfcebrsec}). 
Relations between the invariant and symmetries of knots 
are 
discussed in Section~\ref{symsect}.
}\end{sect}

\begin{sect}{\bf Acknowledgements.\/} {\rm 
Jos\'{e} Barrionuevo, 
Edwin Clark, and Cornelius Pillen 
had helpful programming hints
 for the computation of quandle cocycles.
Seiichi Kamada is being supported by a Fellowship from the 
Japan Society for the Promotion of Science.
}\end{sect}

\section{Definitions of Quandle Cocycle Invariants}
\label{revsec}

\begin{sect}{\bf Definition.\ }
{\rm 
A {\it quandle}, $X$, is a set with a binary operation $(a, b) \mapsto a * b$
such that

(I) For any $a \in X$,
$a* a =a$.

(II) For any $a,b \in X$, there is a unique $c \in X$ such that 
$a= c*b$.

(III) 
For any $a,b,c \in X$, we have
$ (a*b)*c=(a*c)*(b*c). $
 
\noindent
A {\it rack} is a set with a binary operation that satisfies 
(II) and (III).
A typical example of a quandle is a group $X=G$ with 
$n$-fold conjugation 
as the quandle operation: $a*b=b^{-n} a b^n$. 
Racks and quandles have been studied in, for example, 
\cite{Brieskorn},\cite{FR},\cite{Joyce},\cite{K&P}, and \cite{Matveev}.
The axioms for a quandle correspond respectively to the 
Reidemeister moves of type I, II, and III.
(see  also \cite{FR},\cite{K&P}). 
Indeed, knot diagrams were one of the motivations 
 to define such an algebraic structure. 
}
\end{sect}

\begin{figure}
\begin{center}
\mbox{
\epsfxsize=3in 
\epsfbox{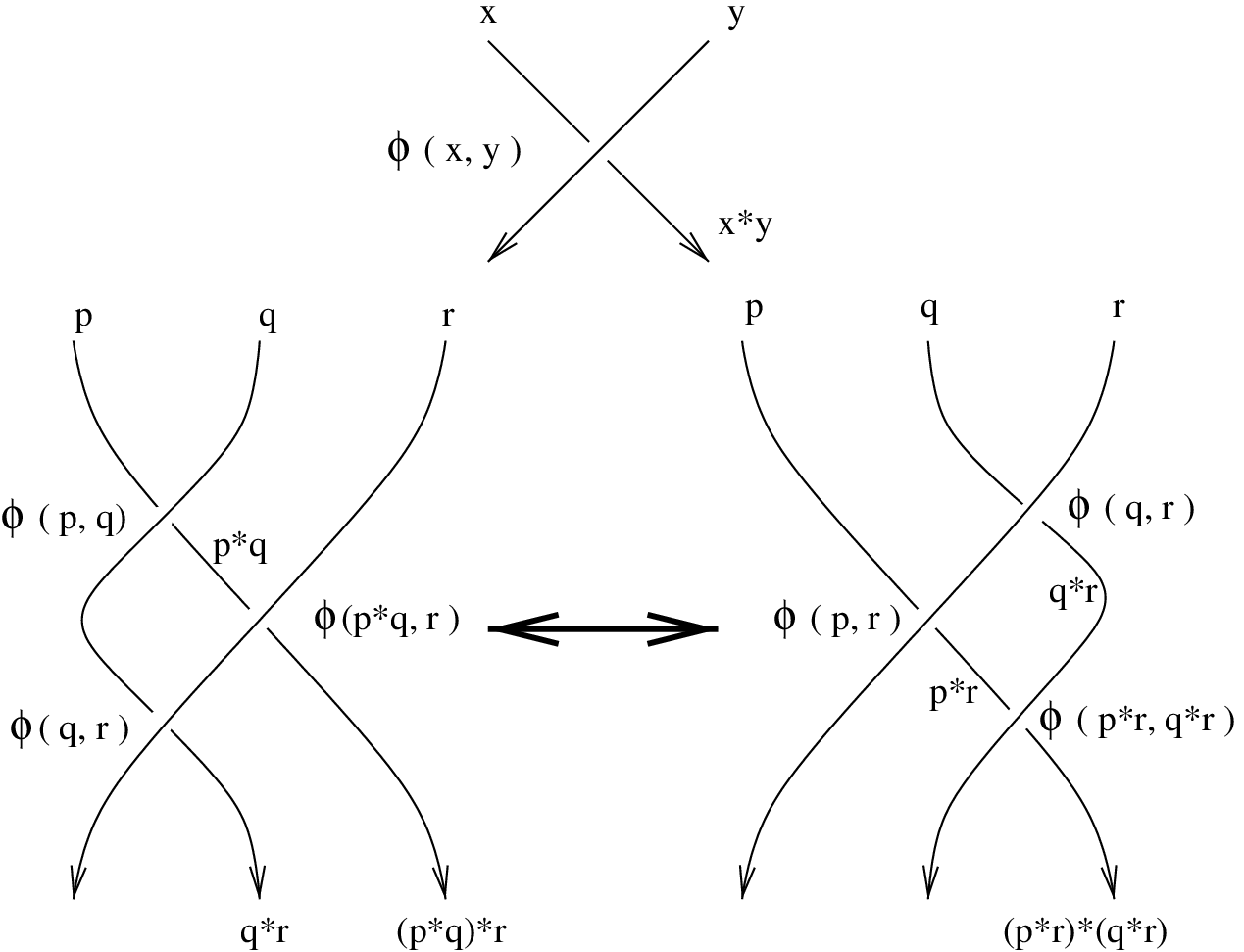}
}
\end{center}
\caption{The 2-cocycle condition and the Reidemeister type III move}
\label{2cocy}
\end{figure}

\begin{sect}{\bf Definition.\/} {\rm
Let $X$ be a rack, 
and let 
 $A$ be an abelian group,  
written additively. 
The cochain group $C^n=C^n(X;A)$ is the 
abelian group 
of functions $f: {\mbox{\rm FA}}(X^n) \rightarrow A$
from the free abelian group generated by $n$-tuples of elements 
 from $X$ to the 
abelian group $A$. 
The {\it coboundary homomorphism} 
$ \delta : C^{n} \rightarrow C^{n+1}$
is defined by
$$ (\delta f)( x_0, \cdots, x_n )
= \sum_{i=1}^{n}  
(-1)^{i-1} 
f( x_0, \cdots, \hat{x}_i , \cdots, x_n) $$
$$+
\sum_{j=1}^n (-1)^j f( x_0 * x_j , \cdots, x_{j-1} * x_j , x_{j+1}, 
\cdots, x_n ). $$
} \end{sect}
(Note: Neither sum includes a $0\/$th term as these terms cancel.)

The rest of the section is a review from \cite{CJKLS}.

\begin{sect} {\bf Lemma.\/} 
The cochain group and the boundary homomorphism 
form  
a cochain complex.
\end{sect}
{\it Proof.\/}
It is a routine calculation (that depends on axiom III of the rack)
 to check that 
$\delta \circ \delta =0$.   $\Box$

\begin{sect}{\bf Definition.\/} {\rm
The cohomology groups of the above complex 
are 
called {\it the rack cohomology groups} and 
are 
denoted by $H_{\mbox{\rm rack}}^n(X, A)$.
Also, the groups of cocycles and coboundaries are denoted by 
$Z_{\mbox{\rm rack}}^n(X,A)$ 
and $B_{\mbox{\rm rack}}^n(X,A)$ respectively. 
Their elements are called {\it $n$-cocycles} and {\it $n$-coboundaries}, 
respectively. 
This definition coincides with the 
cohomology theory 
defined 
in \cite{FRS1} and \cite{FRS2}.

\begin{figure}
\begin{center}
\mbox{
\epsfxsize=3in 
\epsfbox{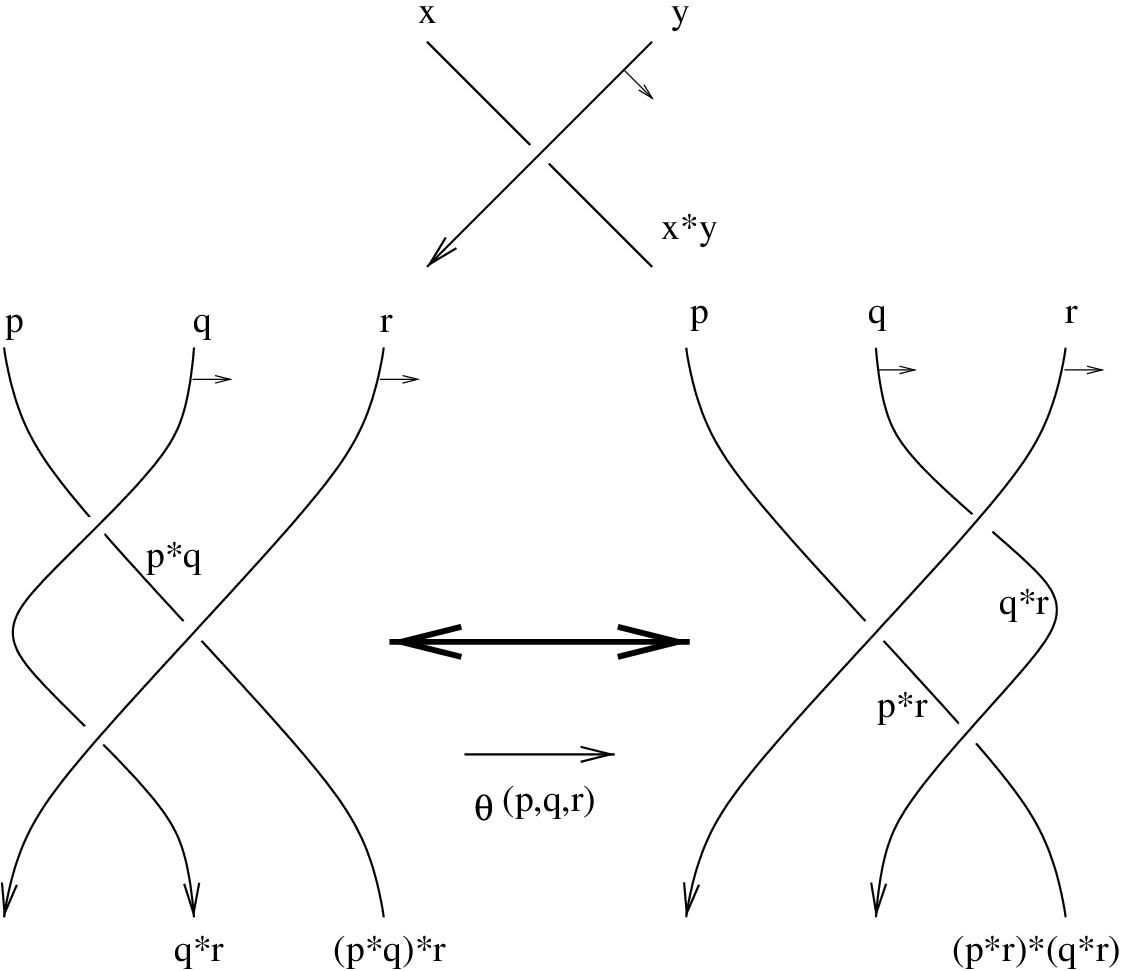}
}
\end{center}
\caption{Assigning a 3-cocycle to a type III move (triple point)}
\label{typeIII}
\end{figure}

For applications, 
we are interested in the case 
when $X$ is a quandle, so we will 
intersect the cocycles and coboundaries 
with a subset that  captures axiom (I) and
 its consequences in higher dimensions.
Let 
 $P^n= \{ f \in C^n: f(\vec{x})=0 $ for all $\vec{x}$ such that 
$x_j=x_{j+1}$ for some $j$ $\}$.
Let $Z^n = Z^n_{\mbox{\rm rack}} \cap P^n$,  and
$B^n = B^n_{\mbox{\rm rack}} \cap P^n$. A straightforward calculation gives:
if $f \in P^n$, then 
$\delta f \in P^{n+1}$
if $X$ is a quandle.
Define
 $$H^n_Q(X,A) = H^n(X,A)=
 (P^n \cap Z_{\mbox{\rm rack}}^n)/ (P^n \cap B_{\mbox{\rm rack}}^n).$$
This group is called the {\it quandle cohomology group}. 
The elements $f\in Z^n(X,A)$ are called {\it quandle $n$-cocycles} or simply 
{\it $n$-cocycles}.
}\end{sect}

\begin{figure}
\begin{center}
\mbox{
\epsfxsize=5in
\epsfbox{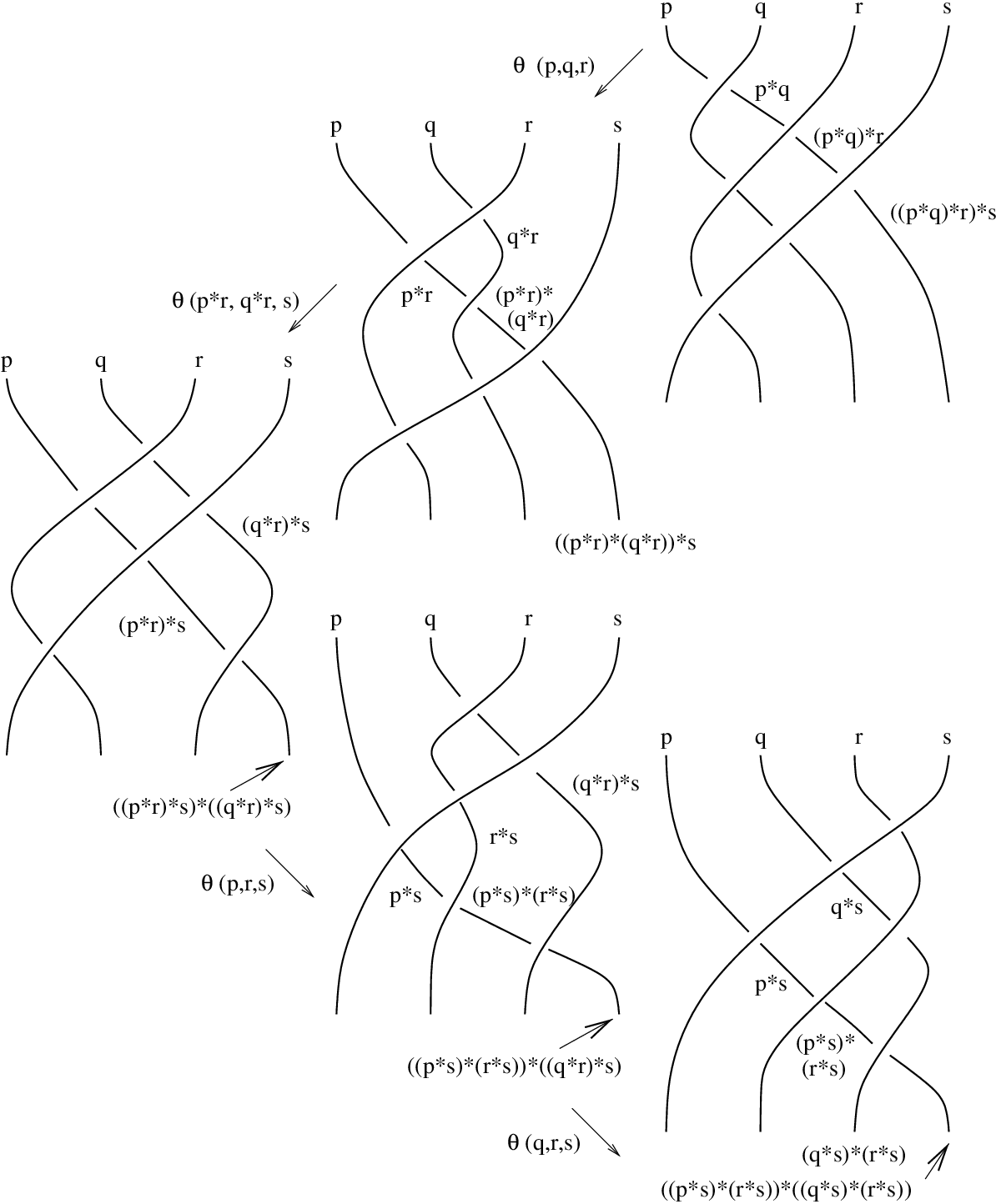}
}
\end{center}
\caption{The tetrahedral move and a cocycle relation, LHS  }
\label{tetraL}
\end{figure}

\begin{figure}
\begin{center}
\mbox{
\epsfxsize=5in
\epsfbox{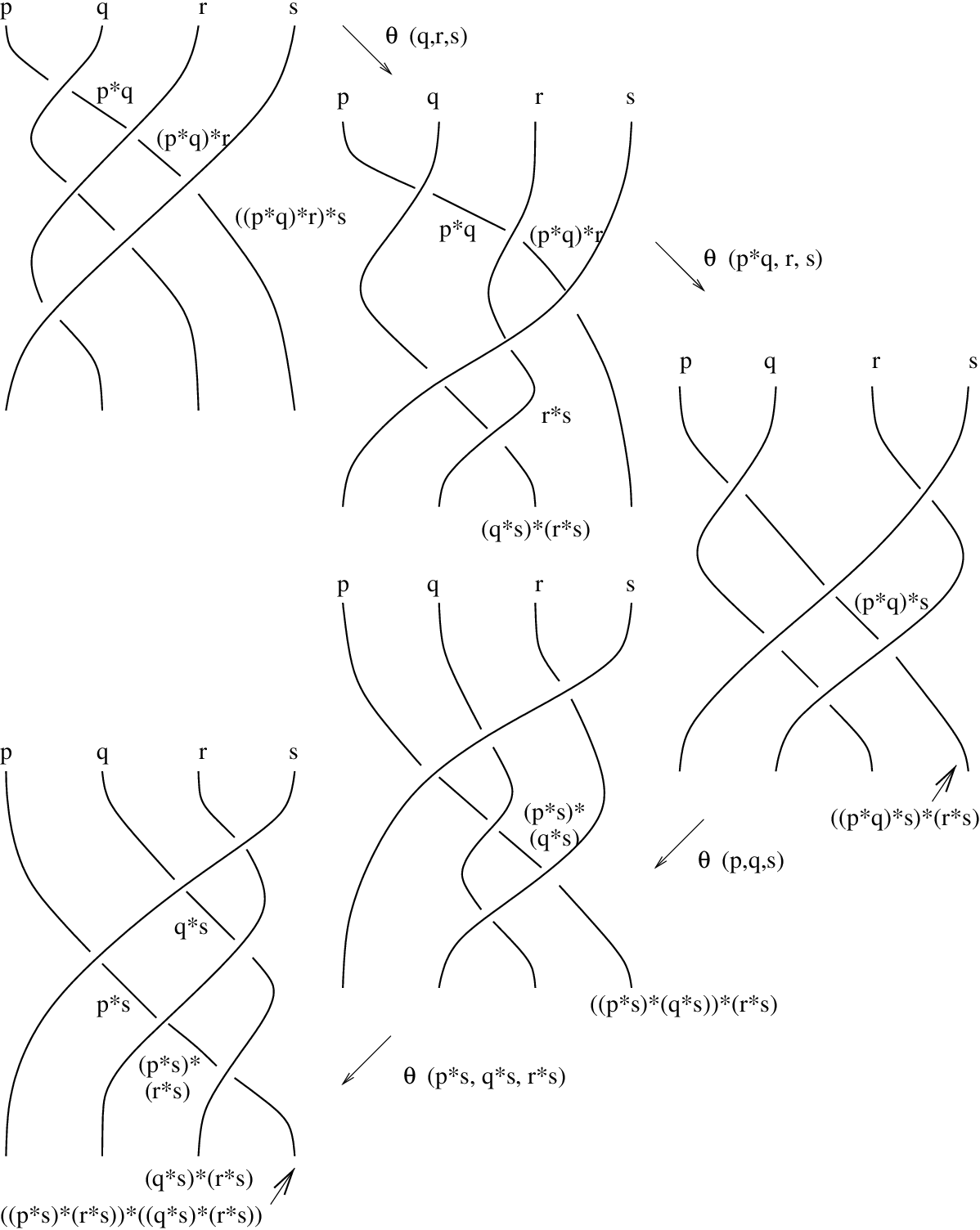}
}
\end{center}
\caption{The tetrahedral move and a cocycle relation, RHS   }
\label{tetraR}
\end{figure}
 
Figures~\ref{2cocy}, \ref{typeIII},  
 \ref{tetraL}, and \ref{tetraR} indicate the relation
 between the cocycle conditions and the moves to classical knots and 
 knotted surfaces.
In Fig.~\ref{2cocy}, $2$-cocycles are assigned to crossings.
Then the type III Reidemeister move corresponds to the $2$-cocycle condition.
For knotted surfaces, $3$-cocycles are assigned to triple points on
projections, which are the type III Reidemeister moves in a movie description.
Figure~\ref{typeIII} shows the assignment of a $3$-cocycle to a type III move.
 Figures~\ref{tetraL} and \ref{tetraR} are movie descriptions of one of the 
generalized Reidemeister moves, called Roseman moves,
which corresponds to the $3$-cocycle condition.

\begin{sect} {\bf Definition.\/} {\rm 
A {\it color} (or {\it coloring}) 
on an oriented  classical knot diagram is a
function ${\cal C} : R \rightarrow X$, where $X$ is a fixed 
quandle
and $R$ is the set of over-arcs in the diagram,
satisfying the  condition
depicted 
in the top
of Fig.~\ref{typeIII}. 
In the figure, a 
crossing with
over-arc, $r$, has color ${\cal C}(r)= y \in X$. 
The under-arcs are called $r_1$ and $r_2$ from top to bottom; they are colored 
${\cal C}(r_1)= x$ and ${\cal C}(r_2)=x*y$.
If the pair of the co-orientation of the over-arc and  that of the under-arc
matches the (right-hand) orientation of the plane, then the 
crossing is called {\it positive}; otherwise it is {\it negative}.
Note that locally the colors do not depend on the 
orientation of the under-arc.

Throughout this paper we consider only finite quandles.
} \end{sect}

\begin{figure}
\begin{center}
\mbox{
\epsfxsize=4in
\epsfbox{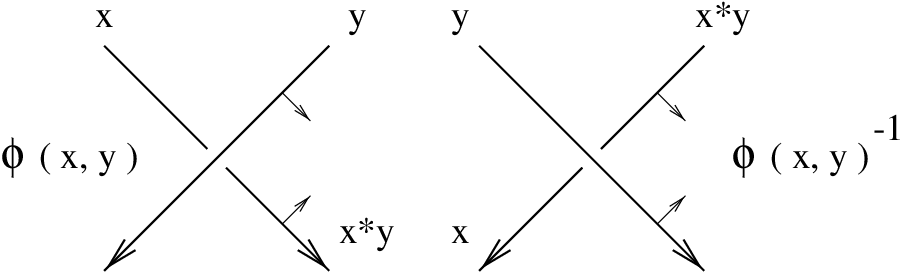}
}
\end{center}
\caption{Weights for positive and negative crossings  }
\label{twocrossings}
\end{figure}

\begin{sect} {\bf Definition.\/} {\rm
Assume that a finite quandle $X$ is given.
Pick a 
quandle 
2-cocycle 
$\phi \in  Z^2(X, A),$ 
and write the coefficient 
group, $A$, multiplicatively. 
Consider a crossing in the diagram.
For each coloring of the diagram, evaluate 
the 2-cocycle on two of the three quandle colors that 
appear near the crossing. One such color is the color on the upper arc and 
is the second argument of the 2-cocycle. 
The other color 
is  the color on 
the under-arc away from which the normal 
arrow points; this is the first argument of the cocycle.

In Fig.~\ref{twocrossings}, 
the two possible oriented and co-oriented 
crossings are depicted. The left is a positive crossing,
and the right is negative.
Let $\tau$ denote a crossing, and ${\cal  C}$ denote a coloring.
When the colors of 
the segments are as indicated, the 
{\it (Boltzmann) weights 
of the crossing},
$B(\tau, {\cal C}) = \phi(x,y)^{\epsilon (\tau)}$,
 are as shown. These weights are assignments of cocycle
values  
to the colored crossings where the arguments are as 
defined in the previous paragraph.

The {\it partition function}, or a {\it state-sum}, 
is the expression 
$$
\sum_{{\cal C}}  \prod_{\tau}  B( \tau, {\cal C}).
$$
The product is taken over all crossings of the given diagram,
and the sum is taken over all possible colorings.
The values of the partition function 
are  taken to be in  the group ring ${\bf Z}[A]$ where $A$ is the coefficient 
group.
} \end{sect}

\begin{sect} {\bf Theorem.\/} 
The partition function is invariant under Reidemeister moves, 
so that it defines an invariant of knots and links.
Thus it will be  denoted by $\Phi (K)$
(or $\Phi_{\phi}(K)$ to specify the $2$-cocycle $\phi$ used).
\end{sect}

With regard to the partition function, we easily obtain the following result.

\begin{sect} {\bf Proposition.\/} 
If $\Phi_{\phi}$ and $\Phi_{\phi '} $ denote the state-sum invariants 
defined from cohomologous cocycles  $\phi$ and $\phi'$
then $\Phi_{\phi} =\Phi_{\phi '} $ (so that $\Phi_{\phi} (K)=\Phi_{\phi '}(K)$
 for any link $K$). 
In particular, 
the state-sum is equal to the number of
colorings of
a given 
knot diagram    
if the $2$-cocycle used for the Boltzmann weight is a coboundary.

\end{sect}

Before we define a similar invariant for knotted surfaces,
 we recall the 
notion of 
knotted surface diagrams. 
See \cite{CS:book} for details and examples. 
Let $f:F \rightarrow {\bf R}^4$ denote a smooth embedding of a closed
surface $F$ into 4-dimensional space.
By deforming the map $f$ slightly by an ambient isotopy 
in 
${\bf R}^4$
if necessary,         
we may assume that
$p \circ f$ is a general position map, 
where  $p: {\bf R}^4 \rightarrow {\bf R}^3$
denotes the 
orthogonal
projection onto an affine subspace.

Along the double curves, one of the sheets (called 
the {\it over-sheet}) lies farther than the other ({\it under-sheet})
with respect to the projection direction.
The {\it under-sheets}
 are coherently broken in the projection,
and such broken surfaces are called {\it knotted surface diagrams}.

When the surface is oriented, we take normal vectors $\vec{n}$
to the projection of the surface such that the triple
$(\vec{v}_1, \vec{v}_2, \vec{n})$ matches the 
right-handed orientation
of 3-space,
where $(\vec{v}_1, \vec{v}_2)$ defines the orientation of the surface.
Such normal vectors are defined on the projection at all points other than
the isolated branch points.

\begin{sect} {\bf Definition.\/} \label{4dcolor} {\rm 
A {\it color} on an oriented  (broken) knotted surface diagram is a
function ${\cal C} : R \rightarrow X$, where $X$ is a fixed 
finite 
quandle
and  where $R$ is the set of regions in the broken surface diagram,
satisfying the following condition at the double point set.  

At a double point curve, two coordinate planes intersect locally.
One is the 
over-sheet $r$, the other is the under-sheet, and the under-sheet is 
broken into two components, say $r_1$ and $r_2$.
 A normal of the over-sheet $r$ points to
one of the components, say $r_2$. 
If ${\cal C}(r_1) = x \in X$, ${\cal C} (r) = y$, then we require that 
${\cal C} (r_2) = x*y$. 
} \end{sect}

\begin{sect} {\bf Lemma.\/} 
The above condition is compatible at each triple point.
\end{sect}
{\it Proof.\/}
The meaning of this lemma is as follows. 
There are 6 double curves near a triple point, giving 
6 conditions on colors assigned.
 It can be checked in a straightforward manner 
that these conditions do not contradict each other.
In particular, there is one of the 4 pieces of the lower
 sheet that receives color
$(a*b)*c$ or $(a*c)*(b*c)$ depending on what path was followed
 to compute the color. 
Since these values agree in the quandle, there is no contradiction.
$\Box$

\begin{sect} {\bf Definition.\/} 
\label{triplepointsign}
{\rm
Note that when three sheets form a triple point, they have relative
positions {\it top, middle, bottom} 
with respect to the projection
direction of $p: {\bf R}^4 \rightarrow {\bf R}^3$.
The {\it sign of a triple point} 
is positive
if the normals of top, middle, bottom sheets in this order 
match the 
right-handed 
orientation of the $3$-space. Otherwise 
the sign is negative.
 
} \end{sect}

\begin{sect} {\bf Definition.\/} {\rm
A (Boltzmann) weight at a triple point, $\tau$, 
is defined as follows.
Let $R$ be the octant from which all normal vectors of the 
three sheets point outwards; let
a  coloring ${\cal C}$ be given.
Let $p$, $q$, $r$ be colors of the 
bottom, middle, and top 
sheets respectively, that bound the region $R$. 
Let $\epsilon (\tau) $ be the sign of the triple point,
and $\theta$ be a 
quandle 
$3$-cocycle.
Then the Boltzmann 
weight $ B( \tau, {\cal C})$ 
assigned to $\tau$ with respect to ${\cal C}$
is defined to be $\theta (p,q,r) ^{ \epsilon (\tau) }$
where $p$, $q$, $r$ are colors described above.
} \end{sect}

\begin{sect} {\bf Definition.\/} {\rm 
The {\it partition function}, or a {\it state-sum}, 
is the expression 
$$ \sum_{\cal C}  \prod_{\tau}  B( \tau, {\cal C} )$$ 
where $B ( \tau, {\cal C})$ 
is
the Boltzmann weight assigned to $\tau$.
As in the classical case, 
we take the coefficient of the cohomology to be 
the group ring ${\bf Z}[A]$ where $A$ is the coefficient 
group written multiplicatively.
} \end{sect}

\begin{figure}
\begin{center}
\mbox{
\epsfxsize=5in
\epsfbox{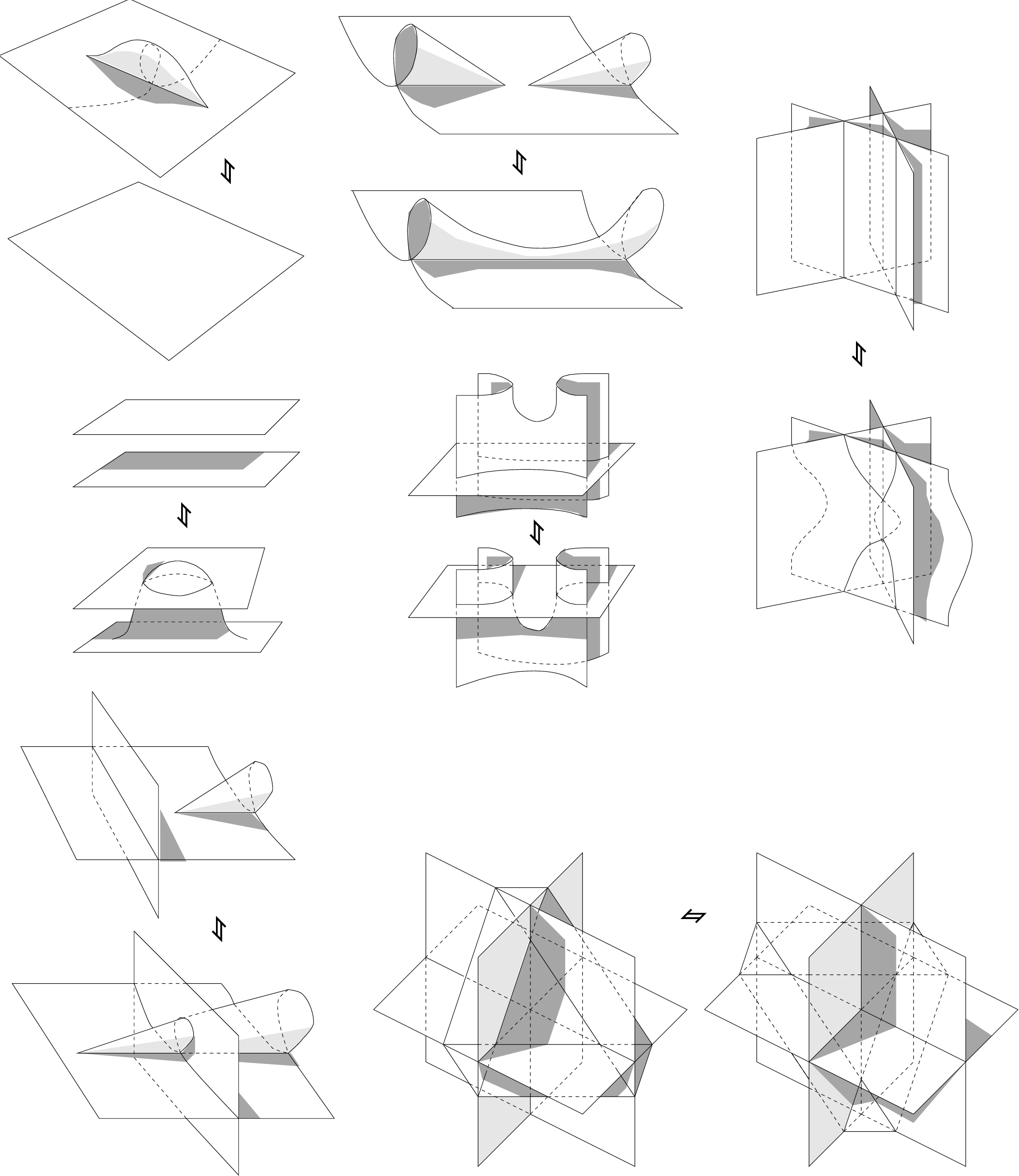}
}
\end{center}
\caption{Roseman moves for knotted surface diagrams  }
\label{rose}
\end{figure}

\begin{sect} {\bf Theorem.\/} 
The partition function does not depend on the choice of 
knotted surface diagram.
Thus it is an invariant of knotted surfaces $F$,
and denoted by $\Phi (F)$
(or $\Phi_{\theta} (F) $ to specify the $3$-cocycle $\theta$ used).
\end{sect} 
{\it Proof Sketch.\/}
Roseman generalized Reidemeister moves to knotted surfaces, 
and 
their projections 
are depicted in Fig.~\ref{rose} \cite{CS:book,Rose}.
Thus two knotted surface diagrams represent isotopic knotted surface 
if and only if the diagrams are related by a finite sequence 
of moves,
called {\it Roseman moves,}
taken from this list.   
The well-definedness of the state-sum is proved by showing that it
remains invariant under the Roseman moves.
In particular, Figs.~\ref{tetraL} and \ref{tetraR}
represent movie descriptions of the terahedral move,
which involves four general position planes (right-bottom of Fig.~\ref{rose}).
Thus these figures show that the state-sum is invariant under this move,
for a specific choice of orientations.
Other cases are checked to prove the well-definedness.
 $\Box$

As in the classical dimension, we can show the following.

\begin{sect} {\bf Proposition.\/} 
If $\Phi_{\theta}$ and $\Phi_{\theta '}$ denote the state-sum invariants 
defined from cohomologous cocycles  $\theta$ and $\theta '$
then $\Phi_{\theta} =\Phi_{\theta '} $ 
(so that $\Phi_{\theta} (K)=\Phi_{\theta '}(K)$ for any knotted surface $K$).
In particular,
if $\theta$ is a $3$-coboundary,
then the state-sum defined above
 is equal to the number of
colorings.

\end{sect}

\section{Non-trivial Cocycles}
\label{cocysec}

Here we define a few exemplary quandles and give some 
quandle cocycles that are not coboundaries.

\begin{sect}{\bf Presentation of the Cohomology Groups.\/}
{\rm Suppose that the coefficient group $A$ is either a cyclic group,
${\bf Z},$ ${\bf Z}_n$, or the rational numbers, ${\bf Q}$.
Define a {\it characteristic function}
$$\chi_x(y) = \left\{ \begin{array}{lr} 1 & {\mbox{\rm if}} \  \ x=y
\\
                                               0   & {\mbox{\rm if}} \ \  x\ne
y \end{array} \right.$$ 
from the free abelian group generated by $X^n$ to the group $A.$
The set $\{ \chi_x: x \in X^n \ \}$ of such
functions
spans the group $C^n_{\mbox{\rm rack}}(X,A)$ of cochains.
Thus if  $f \in C^n_{\mbox{\rm rack}}(X,A)$ is a cochain, then
$$f = \sum_{x \in X^n}  
C_x
\chi_x.$$
We are interested in those $f\/$s 
in $P^n$; {\it i.e.\/} those
homomorphisms that vanish on
$S= \{ (x_1, \ldots , x_n) \in X^n: x_j = x_{j+1} \ 
{\mbox{\rm for some }} \ 
j \}$.
So we can write
$$f = \sum_{x \in X^n \setminus S} \ \ \  
C_x \chi_x.$$

If $\delta f=0$, then $f$ vanishes on expressions of the form
$$\sum_j (-1)^{j+1}(x_0, \ldots , \hat{x}_j, \ldots, x_n)
+
\sum_k (-1)^k(x_0 * x_k, \ldots, x_{k-1}* x_k, x_{k+1}, \ldots , x_n).$$

In computing the cohomology we consider all such expressions as
$(x_0, \ldots x_n)$ ranges over all $(n+1)$-tuples for which
each consecutive pair of elements is distinct. 
By evaluating linear
combinations of 
characteristic functions on these expressions, we determine
those functions that are cocycles.
Similarly, we compute the coboundary on each of the characteristic
functions in the previous dimension, to determine which linear
combinations of characteristic functions are coboundaries.
}\end{sect}

\begin{sect} {\bf Definition \cite{FR}.\/} {\rm
A rack is called {\it trivial} if $x*y=x$ for any $x,y$.

The {\it dihedral quandle} $R_n$ of order $n$ is the quandle consisting of
reflections of the regular $n$-gon with the conjugation as operation.
The dihedral group $D_{2n}$  has a presentation 
$$ \langle x, y | x^2 = 1 = y^n, xyx=y^{-1} \rangle $$
where $x$ is a reflection and $y$ is a rotation of a regular $n$-gon.
The set of reflections $R_n$ in this presentation 
is $\{ a_i = xy^i : i=0, \cdots, n-1 \} $
where we use the subscripts from ${\bf Z}_n$ in the following computations. 
The operation is 
$$a_i * a_j = a_j^{-1} a_i a_j =x y^{j} x y^i x y^j 
= x y^j y^{-i} y^j = a_{2j-i}.$$
Compare with 
the well known $n$-coloring of
 knot diagrams \cite{FoxTrip}.

The quandles with three elements are classified in \cite{FR}
and consist of three isomorphic classes, the trivial one,
$R_3$, and $P_3$.

Let $S_4$ denote the quandle 
with  four elements,  
denoted by  $0,1,2,3,$ 
with the relations
$$\begin{array}{ccccccccc}0 &=& 0 * 0 &= & 1*2 &=& 2* 3 &=& 3*1 \\
                        1  &=& 0*3 &=& 1*1 &=& 2* 0 &=& 3*2 \\
                         2 &=& 0*1 &=& 1*3&=& 2*2 &=& 3*0 \\
                         3 &=& 0*2 &=& 1*0&=&  2*1&=& 3*3. \end{array}$$
This quandle is the following set of 3-cycles in the permutation group on 4 
elements: 
$\{ 0=(243), 1=(134), 2=(142), 3=(123) \}$ with conjugation as the operation.
} \end{sect}

\begin{sect} {\bf Lemma.\/}
Any 
cochain
on a  trivial quandle
is a cocycle. Only the zero map is a coboundary.
\end{sect}
{\it Proof.} This follows from the definitions. $\Box$

\vspace{5mm}

In \cite{CJKLS,FRS1,FRS2,Flower,Greene},  cohomology groups
are computed for some quandles. 
A {\sc Maple} program is found in \cite{MSMaple}.
The techniques that we used in \cite{CJKLS} are applied to give the following:

\begin{itemize}
\item $H^2_Q(R_3, A) = 0$  for any $A$,
\item $H^2_Q(R_4, {\bf Z}_2) = ({\bf Z}_2)^4$,
\item $H^2_Q(R_4, A) = A \times A$  for any $A$ without order 2 elements,
\item $H^2_Q(R_5, A) = 0$  for any $A$,
\item $H^2_Q(R_6, A) = A \times A$  for any $A$,
\item $H^2_Q(S_4, {\bf Z}_2) = {\bf Z}_2$,
\item $H^2_Q(S_4, A) = 0$  for any $A$ without order 2 elements
\end{itemize}

For the third cohomology, we have
\begin{itemize}
\item $H^3_Q(P_3, A) = A \times A$ for any $A$,
\item $H^3_Q(R_3, {\bf Z}_3) = {\bf Z}_3$,
\item $H^3_Q(R_3, A) = 0$  for any $A$ without order 3 elements,
\item $H^3_Q(R_4, {\bf Z}_2) = ({\bf Z}_2)^8$,
\item $H^3_Q(R_4, {\bf Z}_q) = {\bf Z}_q \times {\bf Z}_q$ for any odd prime $q$,
\item $H^3_Q(R_4, {\bf Z}) = {\bf Z} \times {\bf Z} \times {\bf Z}_2 \times {\bf Z}_2$,
\item $H^3_Q(R_4, {\bf Q}) = {\bf Q} \times {\bf Q}$,
\item $H^3_Q(R_5, {\bf Z}_5) = {\bf Z}_5$,
\item $H^3_Q(R_5, A) = 0$  for any $A$ without order 5 elements,
\item $H^3_Q(S_4, {\bf Z}_2) = ({\bf Z}_2)^3$,
\item $H^3_Q(S_4, {\bf Z}_4) = ({\bf Z}_2)^2 \times {\bf Z}_4$,
\item $H^3_Q(S_4, {\bf Z}_q) = 0$ for any odd prime $q$,
\item $H^3_Q(S_4, {\bf Z}) = {\bf Z}_2$,
\item $H^3_Q(S_4, {\bf Q}) = 0$.
\end{itemize}

An important class of quandles are Alexander quandles.

\begin{sect} {\bf Definition \cite{FR,K&P}.\/} {\rm
Let $\Lambda = {\bf Z}[T, T^{-1}]$ be the Laurent polynomial ring
over the integers. Then any $\Lambda$-module $M$
 has a quandle structure defined by
$a*b= Ta + (1-T) b$ for $a, b \in M$. 
} \end{sect}

There are many  $\Lambda$-modules such that 
$L \otimes_{\bf Z} {\bf Z}_n$ is a finite quandle. When
${\bf Z}_n[ T, T^{-1} ]  / (h(T)) $  is a finite quandle,
it is called a {\it (mod $n$)-Alexander quandles.}
Let us point out an exceptional case. Suppose that $\gcd{(a,n)}>1$, and
consider
$
\Lambda_{n,a}
={\bf Z}_n[ T, T^{-1} ]  / (T-a) $. Then 
$\Lambda_{n,a}$ 
is not a quandle because 
axiom II fails.
If $\gcd{(a,n)}=1$, then $\Lambda_{n,a}$ 
is a quandle.

Some of the quandles we have already seen are related to Alexander quandles.

\begin{itemize}
\item
${\bf Z}_n[T, T^{-1}]/(T + 1) \cong R_n$,
\item
${\bf Z}_2[T, T^{-1}]/(T^2 -1) \cong R_4$,
\item
${\bf Z}_2[T, T^{-1}]/(T^2 +T +1) \cong S_4$ 
(The correspondence is 
 $0 \leftrightarrow 0$,
      $1 \leftrightarrow  1$,
      $2 \leftrightarrow  1+T, $ and
      $3 \leftrightarrow T$).
\end{itemize}

\begin{sect} {\bf Table.\/} \label{cdimtable} {\rm 
Extending these results, we give 
the following  table (next page) 
of  cohomological dimensions
 for some quandles with mod $p$ coefficients
(where $p$ is a prime), computed by {\sc Maple}. 
The orders $q$ of the coefficient groups $A={\bf Z}_q$ 
are indicated in the table.
The quandles are chosen as follows. 
First, the programs require 
time beyond our patience for 
quandles of larger than ten elements.
Second, Burau matrices enable us to write a program to compute
invariants for knots in the knot table (which will be presented
in the next section) for Alexander quandles, and all examples
of quandles we have dealt with are Alexander quandles 
(dihedral quandles are the case $T=-1$ in the Alexander quandles). 
Hence, we computed for Alexander quandles of less than 10 elements,
and we considered the Alexander quandles of the form  
${\bf Z}_p[T,T^{-1}]/(h(T))$, where 
$h$ is a polynomial whose leading and constant terms are invertible 
in ${\bf Z}_p$ for the quandle to be finite.
By multiplying by a unit, 
we can assume that such polynomials are monic. 
In the case deg$(h(T))=3$, typical elements are of the form 
$a+bT+cT^2$ and hence the quandle has the order $p^3$, 
and for this order to be less than $10$, we only have 
the choice $p=2$, and the choices of $h(T)$ are 
$T^3+1$, $T^3+T^2+1$, and $T^3+T+1$ that are listed 
in the table below as the last $3$ entrees. 
The cases for smaller degrees are as shown in the table. 
} \end{sect} 


\begin{center}
\begin{tabular}{|l|c||c|c|c|c|c|c|c|c||c|c|c|c|c|} \hline
 &  & \multicolumn{8}{c|}{ $q=$ $\setminus \; \; $ $dim H^2( Q, {\bf Z}_q)$}
 & \multicolumn{5}{c|}{$dim H^3( Q, {\bf Z}_q)$} \\
\hline 
quandle $Q$ & order     & 2&3&5&7&11&13&17&19 &  2&3&5&7&11 \\
\hline \hline 
$R_3$ &3&  0&0&0&0&0&0&0&0&  0&1&0&0&0  \\ \hline 
$R_4$ &4&  4&2&2&2&2&2&2&2&  8&2&2&2&2  \\ \hline      
$R_5$ &5&  0&0&0&0&0&0&0&0&  0&0&1&0&0  \\ \hline 
$R_6$ &6& 2&2&2&2&2&2&2&2&  2&4&2&2&2 \\ \hline    
$R_7$ &7& 0&0&0&0&0&0&0&0&   & & & & \\ \hline   
$R_8$ &8& 4&2&2&2&2&2&2&2& & & & & \\ \hline 
$R_9$ &9& 0&0&0&0&0&0&0&0& & & & & \\ \hline 
${\bf Z}_5[T, T^{-1}]/(T-2)$ &5&  0&0&0&0&0&0&0&0&  0&0&0&0&0 \\ \hline 
${\bf Z}_5[T, T^{-1}]/(T-3)$ &5&  0&0&0&0&0&0&0&0&  0&0&0&0&0 \\ \hline 
${\bf Z}_7[T, T^{-1}]/(T-2)$ &7&  0&0&0&0&0&0&0&0& & & & & \\ \hline 
${\bf Z}_7[T, T^{-1}]/(T-3)$ &7&  0&0&0&0&0&0&0&0& & & & & \\ \hline 
${\bf Z}_7[T, T^{-1}]/(T-4)$ &7&  0&0&0&0&0&0&0&0& & & & & \\ \hline 
${\bf Z}_7[T, T^{-1}]/(T-5)$ &7&  0&0&0&0&0&0&0&0& & & & & \\ \hline 
${\bf Z}_8[T, T^{-1}]/(T-3)$ &8&  4&2&2&2&2&2&2&2& & & & & \\ \hline 
${\bf Z}_8[T, T^{-1}]/(T-5)$ &8&  16&12&12&12&12&12&12&12& & & & & \\ \hline
${\bf Z}_9[T, T^{-1}]/(T-2)$ &9&  0&0&0&0&0&0&0&0& & & & & \\ \hline 
${\bf Z}_9[T, T^{-1}]/(T-4)$ &9&  6&9&6&6&6&6&6&6& & & & &  \\ \hline 
${\bf Z}_9[T, T^{-1}]/(T-5)$ &9&  0&0&0&0&0&0&0&0& & & &  & \\ \hline 
\hline 
${\bf Z}_2[T, T^{-1}]/(T^2+1)$ &4&   2&0&0&0&0&0&0&0& 8&2&2&2&2 \\ \hline
${\bf Z}_2[T, T^{-1}]/(T^2+T+1)$ &4&  1&0&0&0&0&0&0&0& 3&0&0&0&0  \\ \hline
${\bf Z}_3[T, T^{-1}]/(T^2+1)$ &9&    0&1&0&0&0&0&0&0& & & &  & \\ \hline
${\bf Z}_3[T, T^{-1}]/(T^2-1)$ &9&    6&6&6&6&6&6&6&6& & & &  & \\ \hline
${\bf Z}_3[T, T^{-1}]/(T^2+T+1)$ &9&  3&6&3&3&3&3&3&3& & & & &  \\ \hline
${\bf Z}_3[T, T^{-1}]/(T^2-T+1)$ &9&  0&0&0&0&0&0&0&0& & & & &  \\ \hline
${\bf Z}_3[T, T^{-1}]/(T^2+T-1)$ &9&  0&0&0&0&0&0&0&0& & & &  & \\ \hline
\hline
${\bf Z}_2[T, T^{-1}]/(T^3+1)$ & 8&   4&2&2&2&2&2&2&2& & &  & &  \\ \hline
${\bf Z}_2[T, T^{-1}]/(T^3+T^2+1)$ &8& 0&0&0&0&0&0&0&0&  & & &  & \\ \hline
${\bf Z}_2[T, T^{-1}]/(T^3+T+1)$ & 8& 0&0&0&0&0&0&0&0& & & &  & \\ \hline
\end{tabular} \\[5mm]
{Table~\ref{cdimtable} : Cohomological dimensions of  Alexander quandles } 
 \end{center} 


\begin{sect} {\bf Remark.\/} {\rm
We conjecture that the dimension with $A={\bf Z}$ is the dimension
for prime $q$'s when they have the same value for most of them we 
computed. For example, we already know \cite{CJKLS} that 
$H^2(R_4, {\bf Z})={\bf Z}^2$ and 
in Table~\ref{cdimtable} 
we have dimensions $2$ for all $q$ but $q=2$. 
So it is natural to conjecture the same pattern for other quandles 
when it happens.
} \end{sect}

\begin{sect}{\bf Remark.\/} {\rm
The blank entries in Table~\ref{cdimtable} means that {\sc Maple} did  not
finish computations of the first item of the program 
within  $24$ hours. 
} \end{sect}

\begin{sect}{\bf Remark.\/} {\rm
In Table~\ref{cdimtable}, ${\bf Z}_9[T, T^{-1}]/(T-7)$ is omitted since it is isomorphic to
${\bf Z}_9[T, T^{-1}]/(T-4)$; the mapping  
$$f(0)=0, f(1)=1, f(2)=2, f(3)=6, f(4)=7, f(5)=8, f(6)=3, f(7)=4, f(8)=5$$
gives an isomorphism.

Two quandles $Q$ and $R$ are said to be {\it dual quandles}
 if there is a one-to-one correspondence $\gamma: Q \rightarrow R$ 
between their elements, and
$\gamma (a) * \gamma (b) = \gamma (a \bar{*} b )$ where 
$c= a \bar{*} b$ is the unique element $c \in Q$ such that $a=c*b$.
In general,
if $ab \equiv 1 \bmod p$, then  
${\bf Z}_p[T, T^{-1}]/(T-a)$ and ${\bf Z}_p[T, T^{-1}]/(T-b)$
are dual quandles. 
The above quandles, 
${\bf Z}_9[T, T^{-1}]/(T-4)$ and ${\bf Z}_9[T, T^{-1}]/(T-7)$,
are not only isomorphic but also dual of each other. More generally,
we have the following. 
} \end{sect}

\begin{sect}{\bf Lemma.\/} 
The quandles 
${\bf Z}_p[T, T^{-1}]/(T-a)$ and ${\bf Z}_p[T, T^{-1}]/(T-b)$ are 
dual   to 
each other if $ab \equiv 1 \ {\rm mod} \ p$. 
(We assume $\gcd(p,a)=
\gcd(p,b)=1$ 
so that the quandles have $p$ elements.
In this case we denote ${\bf Z}_p[T, T^{-1}]/(T-a)$ by $\Lambda_{p,a}$.)
\end{sect}
{\it Proof.\/}
Consider the identity map ${\bf Z}_p \to {\bf Z}_p$.
$x *_1 y = ax + (1-a)y$ and $x *_2 y = bx + (1-b)y$.
The dual of ${\bf Z}_p[T, T^{-1}]/(T-b)$ has the operation
$x \overline{*_2} y = (x + (b-1)y)b^{-1}$.  However
if $ab \equiv 1$ (hence $b^{-1}=a$), then
 $x \overline{*_2} y = (x +(b-1)y)a  = ax + (1-a)y = x *_1 y$. 
 Hence $\Lambda_{p,a}$
is  dual to $\Lambda_{p,b}$.
$\Box$

In a subsequent paper, 
we will prove that dual quandles have isomorphic cohomology groups.

\begin{sect}{\bf Remark.\/} {\rm 
Recall that  $\Lambda_{9,4}$ and $\Lambda_{9,7}$ are
isomorphic.  This implies that $\Lambda_{9,4}$ is ``self-dual",
{\it i.e.}, 
isomorphic to its dual.
In general, $\Lambda_{p,a}$ is not self-dual, and hence $\Lambda_{p,a}$
is not isomorphic to $\Lambda_{p,b}$ even if 
$ab \equiv 1 \ {\rm mod} \ p$.

We have that $\Lambda_{5,2}$ is not isomorphic to $\Lambda_{5,3}$,
$\Lambda_{7,3}$ is not isomorphic to $\Lambda_{7,5}$.
} \end{sect} 

\begin{sect}{\bf Question.\/} {\rm 
When are $\Lambda_{p,a}$  
is isomorphic to $\Lambda_{p,b}$ with $ab \equiv 1 \ {\rm mod}  \ p $?
Equivalently, when is $\Lambda_{p,a}$ self-dual?
} \end{sect}

\section{Computations  for Knot Table} 
\label{tablesec}

For Alexander quandles, closed braids and Burau representations
can be used to compute the cocycle invariants. 
In this section we present 
computational results obtained by 
{\sc Maple}. 
The {\sc Maple} program we used in this section can be found in 
\cite{MSMaple}. 
Here is a sketch of an algorithm we used.

\begin{enumerate}
\item 
We used the same program that we used in the preceding section
to compute cocycle groups for Alexander quandles. 
Then we made specific choices of cocycles. 
The choices are made after some experiments.

\item 
If
$B= \left[ \begin{array}{cc} 0&T \\1&1-T \end{array} \right]$  
is the Burau matrix, and $v$ is a row vector of $\Lambda^2$,
where $\Lambda$ is an  Alexander quandle, representing the colors assigned to
the top two strings at a positive crossing, then the 
colors at the bottom 
strings  
are represented by a vector $vB$. 
Furthermore, the cocycle contribution
at this crossing is $\phi(v_1, v_2)$ where $v= (v_1, v_2)$. 
The actual vectors have larger dimensions depending on the braid index
of the given closed braid, but it is simply shifted to the 
position of the corresponding braid generator.
\item
If $B(w)$ is the Burau representation of
a braid word $w$, 
then a vector $v$ colors the closed braid
if and only if $vB(w)=v$. The {\sc Maple} program searches for 
all such colors, and then evaluates the state-sum contributions.
\end{enumerate}
 
We used the closed braid form for knots in the table 
given in \cite{Jones} up to (including) 
 $9$ crossings
to obtain the following results.
We list some of the  quandles with non-trivial cohomology computed in the
preceding section. 
We chose $q$ with highest dimensions among all non-zero dimensions.
We conjecture that the invariants are positive integers for
knots with $R_n$ if $n$ is even. 
Note that some quandles have non-trivial values for links \cite{CJKLS},
and could define interesting invariants for links even if they 
are trivial for knots.

\begin{enumerate}
\item 
For ${\bf Z}_2[T, T^{-1}]/(T^2+1)$
with the coefficient $A={\bf Z}_2$ 
all knots in the table up to $9$ crossings have the 
trivial invariant (have value $4$).
Thus none colors  non-trivially.  

For 
 ${\bf Z}_9[T, T^{-1}]/(T-4)$ and  ${\bf Z}_3[T, T^{-1}]/(T^2+T+1)$ 
with 
 $A={\bf Z}_3$, 
all knots of the first half ($42$ out of $84$) of the table up to $9$
crossings have the trivial invariant (have value $9$).
We stopped the program  after $5$ days.

\item
For $S_4={\bf Z}_2[T, T^{-1}]/(T^2+T+1)$ with the coefficient $A={\bf Z}_2$, 
we used the cocycle $\phi=\prod \chi_{(a,b)}$ where 
the product is taken over all pairs $(a,b)$ such that 
$a,b \in \{0,1,T+1 \}$ and $a \neq b$. The invariants take the following values.
\begin{itemize} 
\item  
$4(1+3t)$ for $3_1$, $4_1$, $7_2$, $7_3$, $8_1$, $8_4$, $8_{11}$,
$8_{13}$, $9_1$, $9_6$, $9_{12}$, $9_{13}$, $9_{14}$, $9_{21}$, 
$9_{23}$, $9_{35}$, $9_{37}$.
\item 
$16(1+3t)$ for $8_{18}$, $9_{40}$.
\item
$16$ for $8_5$, $8_{10}$, $8_{15}$, $8_{19}-8_{21}$,
$9_{16}$, $9_{22}$, $9_{24}$, $9_{25}$, $9_{28}-9_{30}$,
$9_{36}$, $9_{38}$, $9_{39}$, $9_{41}-9_{45}$, $9_{49}$. 
\item
$4$ otherwise.
\end{itemize}

\item
For ${\bf Z}_3[T, T^{-1}]/(T^2+1)$ with the coefficient $A={\bf Z}_3$, 
we used the cocycle 
\begin{eqnarray*}
\phi &=& {\chi _{2\,T, \,2}} +
 2\,{\chi _{2\,T, \,T}} +
 2\,{\chi _{2\,T, \,1 + T}} +
 2\,{\chi _{2\,T, \,2 + T}} +
 2\,{\chi _{1 + 2\,T, \,0}} +
 2\,{\chi _{1 + 2\,T, \,1}} +
 2\,{\chi _{1 + 2\,T, \,2}} +
 {\chi _{1 + 2\,T, \,T}}
\\ & &  + 
{\chi _{1 + 2\,T, \,1 + T}} + 
{\chi _{1 + 2\,T, \,2 + T}} +
{\chi _{2\,T, \,1}} + 
{\chi _{2\,T, \,0}}  +
 2\,{\chi _{0, \,2\,T}} +
 2\,{\chi _{0, \,2 + T}} +
{\chi _{0, \,T}} + 
{\chi _{0, \,2}}
\\ & &  +
 2\,{\chi _{0, \,1}}  +
 2\,{\chi _{1, \,2\,T}} +
{\chi _{1, \,2 + T}} +
 2\,{\chi _{1, \,1 + T}} +
2\,{\chi _{1, \,2}} +
\,{\chi _{1, \,0}} +
{\chi _{0, \,1 + 2\,T}} +
{\chi _{2, \,1}}
\\ & &  + 
2\,{\chi _{2, \,0}} +
{\chi _{1, \,1 + 2\,T}} +
{\chi _{2, \,1 + 2\,T}} +
 2\,{\chi _{2, \,2\,T}} +
{\chi _{2, \,1 + T}} 
+ 2\,{\chi _{2, \,T}} +
{\chi _{T, \,2}} +
2\,{\chi _{T, \,0}}
\\ & &  +
{\chi _{T, \,1 + T}} +
 2\,{\chi _{T, \,2 + T}} +
 {\chi _{T, \,2\,T}} +
 2\,{\chi _{T, \,1 + 2\,T}} +
{\chi _{1 + T, \,1}} +
 2\,{\chi _{1 + T, \,T}}  +
 2\,{\chi _{1 + T, \,2}}+
{\chi _{1 + T, \,2 + T}} 
\\ & &  + 
2\,{\chi _{1 + T, \,1 + 2\,T}} +
{\chi _{1 + T, \,2\,T}} +
{\chi _{2 + T, \,0}} +
 2\,{\chi _{2 + T, \,1+ T}} + 
{\chi _{2 + T, \,T}} +
{\chi _{2 + T, \,2\,T}} +
 2\,{\chi _{2+ T, \,1 + 2\,T}}.
\end{eqnarray*}
 The invariants 
take the following values.
\begin{itemize}
\item 
$9(1+4t+4t^2)$ for $4_1$, $5_2$, $8_3$, $8_{17}$, $8_{18}$, $8_{21}$,
$9_6$, $9_7$, $9_{11}$, $9_{24}$, $9_{26}$, $9_{37}-9_{39}$, 
$9_{47}$.
\item 
$297+216t+216t^2$ for $9_{40}$.
\item
$81$ for $6_3$, $8_2$, $8_{19}$, $8_{24}$, $9_{12}$, $9_{13}$, 
$9_{46}$.
\item 
$9$ otherwise.
 \end{itemize}

\end{enumerate}

\begin{sect} {\bf Remark.\/} {\rm
With $S_4$, we do not know any knot with  the invariant
 not equal to  an integer or $k(4+12t)$ for an integer $k$. 
However, the torus link $T(5, 15)$ has the invariant
$544 + 480t$ (with $S_4$ and $A={\bf Z}_2$, with the same
cocycle as above). 
} \end{sect}

\section{Computations for Torus Knots}
\label{torussec}

The Burau matrices for torus knots have  periodicities with some quandles.
In this section we give such periodicities and use them to compute
 the invariants for 
torus knots.  
Throughout the section 
let $Q$ be a finite  Alexander quandle ($Q={\bf Z}_m[T, T^{-1}]/ (h(T))$
for some positive integer $m$ and a Laurent polynomial $h(T)$).

\begin{figure}
\begin{center}
\mbox{
\epsfxsize=1in
\epsfbox{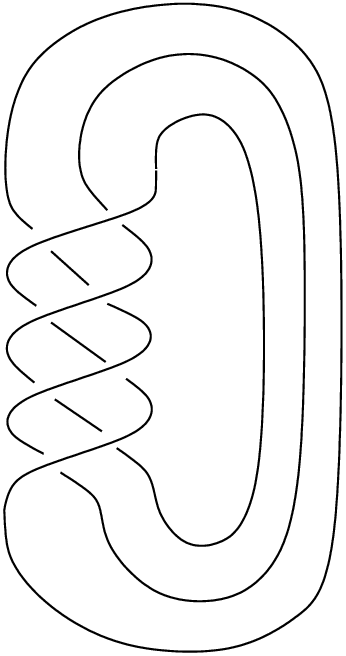}
}
\end{center}
\caption{The torus knot $T(3,4)$  }
\label{torus34}  
\end{figure}

Consider a $(n,k)$-torus knot or link $L=T(n,k)$ (Fig.~\ref{torus34}
depicts $T(3,4)$).
Then  $L$ may
 be represented as 
the closure of 
an $n$-string braid $\beta$ of the form
$(\sigma_{n-1} \sigma_{n-2} \ldots \sigma_{1})^k$.  If we color
this braid by the quandle $Q$, any possible color of the closed
braid $\hat{\beta}$ can be uniquely determined by a choice of
colors on the top segment of the $n$-strands expressed as a
vector, say 
$[ a_1, a_2, \cdots, a_n]$
where $a_i \in Q$ for $i=1,2, \ldots n$,
which we call a color vector.
 Note that 
$[ a_1, a_2, \cdots, a_n]=\sum^{n}_{i=1}a_ie_i$, where $e_i$
is a basic unit vector, 
with all $0$ except a $1$ in the $i$-th position.

The braid word $\sigma_{n-1} \sigma_{n-2} \ldots \sigma_{1}$ on
strands colored 
$[ a_1, a_2, \cdots, a_n]$ on the top ends alters this $n$-tuple of colors
to the color vector $[a_n, a_1 * a_n, \cdots, a_{n-1}*a_n ]$
assigned on the strings below the word
 $\sigma_{n-1} \sigma_{n-2} \ldots \sigma_{1}$.
  This operation can be represented by
multiplying the color vector on the right by the $n \times n$
matrix $$A= \left[
\begin{array}{cccc} 0& T& \ldots &0\\ \vdots & 0 & \ddots \\ 0& & & T\\
1 & 1-T & \ldots & 1-T \end{array} \right] $$
which is a product of Burau matrices. 

So for any $j \in {\bf N} \cup \{0\}$ with $ 0 \leq j \leq k$, the
color vectors after
the braid word 
$(\sigma_{n-1} \sigma_{n-2} \ldots \sigma_{1} )^j $  is
$[ a_1, a_2, \cdots, a_n] A^j=
(\sum_{i=1}^{n}a_ie_i)A^j= \sum_{i=1}^{n}a_i(e_iA^j). $ 
Thus any
choice of color vectors induces nontrivial color of $L$ if
$\sum_{i=1}^{n}a_i(e_iA^j) = \sum^{n}_{i=1}a_ie_i$, which occurs
if and only if $e_iA^j=e_i$ for every $i=1,2,\ldots, n.$  Note
that this is equivalent to $A^j=I.$

\begin{sect} {\bf Definition.\/} {\rm The {\em color period} of a quandle
$Q$
for the family ${\cal T}(n)=\{T(n,m):m \in {\bf Z}\}$ is the
 minimum positive integer k such that $A^k=I.$ } \end{sect}

\begin{sect} {\bf Lemma.\/} If $e_nA^{kn}=e_n$ for some 
$k \in {\bf Z}$, then $e_jA^{kn}=e_j$ for all $j=1,2,\ldots ,n.$
\end{sect}
{\it Proof.\/} Let $e_nA^{kn}=e_n$.  Then $e_nA^{(k-1)n+1}=Ke_1$ for
some $K \in Q.$  Since $e_jA=Te_{j+1}$ for
$j=1,2,\ldots,n-1,$
$$e_nA^{kn}=e_nA^{(k-1)n+1}A^{n-1}=Ke_1A^{n-1}=KT^{n-1}e_n=e_n.$$
Hence, $KT^{n-1}=1.$

Now consider $e_jA^{kn}$ for $j \in \{1,2,\ldots ,n\}.$ Then,
$$e_jA^{kn}=(e_jA^{n-j})A^{(k-1)n+j}=T^{n-j}(e_nA^{(k-1)n+1})A^{j-1}
=T^{n-j}Ke_1A^{j-1}=T^{n-j+j-1}Ke_j=e_j. \Box $$

\begin{sect} {\bf Proposition.\/}
 For $A$ as above, the color periods, $p$,
of the following quandles $Q$
 for the family ${\cal T}(n)=\{T(n,m):m \in {\bf Z}\}$
are:
\begin{itemize}
    \item $Q={\bf Z}_8[T,T^{-1}]/(T-3)$, $p=\left\{
     \begin{array}{lll} 2n &  {\mbox{\rm if}} \  \ n = 2k+1, & k \in
     {\bf N}, \\ 4n & {\mbox{\rm if}} \  \ n  =2k, & k \in {\bf N}.
\end{array}
      \right.$

    \item $Q={\bf Z}_3[T,T^{-1}]/(T^2+1)$, $p= \left \{
    \begin{array}{lll} 2n & {\mbox{\rm if}} \ \ n =4k+2, &  k \in
     {\bf N}, \\ 3n & {\mbox{\rm if}} \ \ n=4k, & k \in
     {\bf N}, \\    4n & {\mbox{\rm otherwise}}. \end{array}
     \right.$

    \item $Q={\bf Z}_3[T,T^{-1}]/(T^2-1)$, $p= \left \{
    \begin{array}{lll} 2n & {\mbox{\rm if}} \ \  n=2k+1, &  k \in
     {\bf N},\\ 3n & {\mbox{\rm if}} \ \ n=2k, &  k \in
     {\bf N}. \end{array} \right.$

 \item $Q=  {\bf Z}_2[T,T^{-1}]/(T^3-1)$, 
           $p= \left \{  \begin{array}{lll} 2n  & {\mbox{\rm if}} \ \ n =3k, &
              k \in {\bf Z}, \\
             3n & {\mbox{\rm otherwise}}.
\end{array} \right.$

 \item $Q={\bf Z}_2[T,T^{-1}]/(T^2+T+1)=S_4$, 
 $p= \left \{  \begin{array}{lll} 2n  & {\mbox{\rm if}} \ \ n =3k, & k \in
     {\bf Z}, \\
3 & {\mbox{\rm if}} \ \ n=2,  \\ 3n & {\mbox{\rm otherwise}}.
\end{array} \right.$

     \item $Q=R_j$ for $ j=2k, k \in {\bf N}$, $p= \left \{
    \begin{array}{lll} 2n & {\mbox{\rm if}} \ \ n =2k+1, &  k \in
     {\bf N}, \\ kn & {\mbox{\rm if}} \ \ n=2k, & k \in
     {\bf N}. \end{array} \right.$

    \item $Q= \left \{ \begin{array}{l} {\bf
    Z}_8[T,T^{-1}]/(T-5), \;  {\mbox{\rm or}} \\ {\bf
    Z}_2[T,T^{-1}]/(T^2+1) \end{array} \right.$,  $p=2n.$

    \item $Q= \left \{ \begin{array}{l} {\bf
    Z}_9[T,T^{-1}]/(T-4), \;  {\mbox{\rm or}}\\ {\bf
    Z}_9[T,T^{-1}]/(T-7), \;  {\mbox{\rm or}}\\ {\bf
    Z}_3[T,T^{-1}]/(T^2+T+1) \end{array} \right.$ , $p=3n.$

\end{itemize}
\end{sect}
{\it Proof.\/} We prove the cases $Q={\bf Z}_2[T,T^{-1}]/(T^2+T+1)=S_4$,
 $Q={\bf Z}_3[T,T^{-1}]/(T^2+1)$,
and $Q=R_j$ for $ j \in {\bf N}$.  The rest are similar.

First consider $Q={\bf Z}_2[T,T^{-1}]/(T^2+T+1)=S_4$, and let
$n=3k$ for some $k \in {\bf N}$.

Consider $e_nA=[  1, 1-T, \cdots , 1-T] $. 
 Denote by $q_b(a)$ the operation $a*b$ in $S_4$.
Note that $$[ a_1, a_2,  \cdots,  a_n]A 
=[  a_n, q_{a_n}(a_1),  q_{a_n}(a_2), \cdots , q_{a_n}(a_{n-1})] $$
 and that $q_b^3(a)=(q_b \circ q_b \circ q_b)(a)=a$ for
all $a,b \in S_4$.  
Hence  $$e_nA^n=[  1-T, \cdots ,  1-T, q_{1-T}^{n-1}(1)] 
=[ 1-T,  \cdots ,  1-T,  T]. $$ 
Furthermore,
$$e_nA^{n+1}=[  T, q_T(1-T), \cdots , q_T(1-T)]
 = [T, 0, \cdots,  0], $$ 
and
 $$e_nA^{2n}=[  0, \cdots ,  0, q_0^{n-1}(T)] 
 = [   0,  \cdots ,  0,  1] =e_n.$$
Hence, since $e_nA^p \neq e_n$ for any $0 < p < 2n$ and
$e_iA^{2n}=e_i$ for all $i=1,2, \ldots n$
by the preceding Lemma, 
  the color period  is $2n$.

Now consider $n=3k+1$ for some $k \in {\bf N} \cup \{0\}$.
  Again,
consider $e_nA=[ 1, 1-T, \cdots ,  1-T] $. 
 Note that 
 $$e_nA^n=[  1-T, \cdots ,  1-T,  q_{1-T}^{n-1}(1)] 
=[  1-T,   \cdots ,  1-T,  1] ,$$ 
and so $e_nA^{n+1}=[  1,  T,  \cdots ,  T,] .$ 
 Thus, $e_nA^{2n}=[  T,  \cdots,  T, 1] $, 
which leads to $e_nA^{2n+1}=e_1$.  From
this, it can be seen that
$$ 
e_nA^{3n}=T^{n-1}e_n=T^{3k+1-1}e_n$$
 for some $k \in {\bf Z}$.  Thus $e_nA^{3n}=T^{3k}e_n=e_n$.
Thus 
 the color period is $3n.$

For
$n=3k+2$ 
for some $k \in {\bf N} \cup \{0\}$.  Again,
consider $e_nA=[  1,  1-T,  \cdots ,  1-T] $. 
 Note $$e_nA^n=[ 1-T, \cdots ,  1-T,  q_{1-T}^{n-1}(1)]
 =[  1-T,  \cdots , 1-T,  0] ,$$ 
and so $e_nA^{n+1}=[ 0,  1,  \cdots , 1] , $ 
which is $e_n$ when $n=2$. 
 Thus, $e_nA^{2n}=[  1, \cdots  , 1,  T+1] $,
which leads to $e_nA^{2n+1}=(T+1)e_1=T^2e_1$. From this, it can be
seen that 
$$e_nA^{3n}= 
T^{2+(n-1)}e_n=T^{n+1}e_n=T^{3k+2+1}e_n=e_n$$
since $n=3k+2$ for some $k \in {\bf Z}.$
Thus
 the color period  here is $3n.$

For $Q={\bf Z}_3[T,T^{-1}]/(T^2+1)$, we proceed as with $S_4$. If
$n \equiv 0$  (mod 4), then 
we have the following computations.
\begin{eqnarray*}
e_nA & = & [ 1, 1-T,  \cdots , 1-T] \\
 &  = & [  1,  1+2T,  \cdots , 1+2T] \\
e_nA^n & = & [ 2T+1, \cdots , 2T+1, 2T+2] \\
e_nA^{n+1} & = & [ 2T+2, T+2, \cdots , T+2] \\
e_nA^{2n} & = & [ T+2, \cdots ,  T+2,  T] \\
e_nA^{2n+1} & = & [T, 0, \cdots , 0] .
\end{eqnarray*}
  Thus, $e_nA^{3n}=e_n.$
 So for $n \equiv 0$ (mod 4),
 the color period  of a $n$-strand torus link is $3n$.
Now let $n \equiv 2$ (mod 4).  Then $e_nA^n=[ 2T+1, 2T+1, \cdots , 2T]$,
 and so $e_nA^{n+1}=2Te_1$. Thus $e_nA^{2n}=e_n$.
 So, the color period is $2n$.
For $n \equiv 1$ (mod 4), one computes 
\begin{eqnarray*} 
e_nA^n &= & [ 2T+1,  \cdots , 2T+1, 1] \\
e_nA^{n+1} &=&[  1, 2, \cdots , 2] \\
e_nA^{2n} &=& [ 2, \cdots , 2, 1] \\
e_nA^{2n+1} &=& [ 1, T+1,  \cdots , T+1] \\
e_nA^{3n} & = & [ T+1, \cdots , T+1, 1] \\
e_nA^{3n+1} &= & [ 1, 0, \cdots , 0] .
\end{eqnarray*}
  Hence, $e_nA^{4n}=e_n$.
For $n \equiv 3$ (mod 4), one computes 
\begin{eqnarray*}
e_nA^n &=& [ 2T+1,  \cdots , 2T+1, T+1] \\
e_nA^{n+1} &=&[ T+1, T, \cdots , T] \\
e_nA^{2n} & =& [ T, \cdots , T, T+2] \\
e_nA^{2n+1} &=& [ T+2, 2T+2, \cdots , 2T+2] \\
e_nA^{3n} & = & [2T+2, \cdots , 2T+2, 2] \\
e_nA^{3n+1}&=&[ 2, 0, \cdots , 0] 
\end{eqnarray*}
  Hence, $e_nA^{4n}=e_n$.
Thus the color periods  are $4n$.

Finally, for $Q=R_j$ for $ j=2k, k \in {\bf N}$, note that $Q$ is
isomorphic to ${\bf Z}_j[T,T^{-1}]/(T+1)$.  So, $T=-1$, which
makes $1-T=2$, so
 $e_nA=[ 1, 2, \cdots , 2] .$ 
 If $n$ is odd, $e_nA^n=[ 2,  \cdots ,  2, 1] $,
 so $e_nA^{n+1}=[1, 0, \cdots , 0] $. 
 This means that $e_nA^{2n}=e_n$. 
  Thus the color period is $2n$.
If $n$ is even $e_nA^n=[ 2, \cdots , 2, 3] $. 
Noting that 
$n*(n+1)=n+2$, for $i=1,2, \ldots , k$ we will get
\begin{eqnarray*}
e_nA^{in} &=& [ 2i, \cdots , 2i, 2i+1] , \\
e_nA^{in+1} &= & [ 2i+1, 2i+2, \cdots , 2i+2] .
\end{eqnarray*}
  When $i=k-1, 2i+2=2(k-1)+2=2k=j \equiv 0$ (mod $j$),
 and so $e_nA^{nk}=e_n$. $\Box$

\begin{sect}{\bf Proposition.\/} Let $p$ be the color period of a
quandle $Q$ for the family ${\cal T}(n)$.  Let $\theta \in Z^2(Q,{\bf Z}_q).$
  Then $\Phi_{\theta}(T(n,k))$ is periodic with
respect to $k$ with period at most $pq$. \end{sect}
{\it Proof.\/} Let $Q$ be a quandle that has color period $p$ with
respect  to ${\cal T}(n)$.   Then any torus
knot or link $T(n,s)$ with $s>p$, $s=ph+r \; (0 \leq r <  p)$, 
 can be thought of as 
the closure of a finite number $(h)$ of
the  block  $(\sigma_{n-1} \sigma_{n-2}\ldots \sigma_{1})^p$ 
and a remainder block $(\sigma_{n-1} \sigma_{n-2}\ldots \sigma_{1})^r$.  
It has been shown that the
color vector after the $p\/$th block is the same as the initial
(at the top of the braid representation) color vector. Hence
the $(ph+r)$-th color vector is the same as the $r$-th color vector
for $r=1,2,\ldots,p-1.$  
Let $v$ be an initial color vector 
which colors $T(n,s)$. 
Let $t^{\alpha}$ and $t^{\beta}$  be the contribution of
 to the state-sum
 of  $(\sigma_{n-1} \sigma_{n-2}\ldots \sigma_{1})^p$ 
and $(\sigma_{n-1} \sigma_{n-2}\ldots \sigma_{1})^r$ respectively,
for the color with the initial color vector $v$.
Then the contribution of $(\sigma_{n-1} \sigma_{n-2}\ldots \sigma_{1})^s$
is $t^{\alpha h + \beta}$ since   $s=ph+r$. 
Hence if $h=q$ (with $A={\bf Z}_q$), we obtain 
$t^{\alpha h + \beta} = t^{\beta}$, and therefore 
$\Phi(T(n,ph+r))=\Phi(T(n,r))$. $\Box$

\begin{sect}{\bf Example.\/}  {\rm 
The above Proposition implies that the invariants 
$\Phi(T(n,r))$ for $0 \leq r < pq$ determine all the rest
of $\Phi(T(n,s))$. Here we present the list of 
$\Phi(T(2, s))$, for a  cocycle in  
$Z^2( {\bf Z}_8[T,T^{-1}]/(T-5), {\bf Z}_2 )$ 
(as this quandle displays 
a variety 
of polynomial values). 
In this case the period $pq$ is $8$. In the following table,
 $k$ represents any integer.
The cocycle used is 
$$\phi= {\chi _{0, \,1}} + {\chi
_{0, \,5}} + {\chi _{1, \,5}} + {\chi _{ 2, \,1}} + {\chi _{2,
\,5}} + {\chi _{3, \,5}} + {\chi _{5, \,1}}
 + {\chi _{7, \,1}}.$$

\begin{center}
\begin{tabular}{|l|c|} \hline
 torus link    & invariant
\\
\hline \hline 
$T(2,8k)$  & $64$ 
\\ \hline 
 $T(2,1 + 8k)$ & $8$ \\
\hline $T(2,2+8k)$ & $28 + 4t$ \\
\hline $T(2,3+8k)$ & $8$ \\
\hline $T(2,4+8k)$  & $48 + 16t$ \\
\hline $T(2,5+8k) $  & $8 $ \\
\hline $T(2,6+8k)$  & $28 + 4t$ \\
\hline $T(2,7+8k)$  & $8$ \\
  \hline
\end{tabular} 
 \end{center}

} \end{sect}

\begin{sect}{\bf Table.\/} \label{torustable} {\rm 
The following table presents 
calculations
of $\Phi(T(n,k))$ with some quandles for some small values of $k$. 
The quandles are Alexander quandles with $A={\bf Z}_q$.
Although they  do not give complete
set of initial conditions (as the periods $pq$ are sometimes too large),
they give useful  information together with the proposition above.
The table lists  $\Phi(T(n,k))$ for $n \leq k$ since 
$T(n,k)=T(k,n)$.
} \end{sect}

We used the following cocycles: 

\begin{eqnarray*}
\Theta_1 &=& {\chi _{0, \,2}} + {\chi _{1, \,0}} + {\chi _{1,
\,2}} + {\chi _{ 2, \,0}}
\\
\Theta_2 &= & {\chi _{0, \,3}} + {\chi _{4, \,1}} + {\chi _{2, \,5}}
+ {\chi _{ 2, \,3}} + {\chi _{2, \,1}} + {\chi _{1, \,5}} + {\chi
_{1, \,3}}
 + {\chi _{3, \,5}} + {\chi _{3, \,1}} 
\\ 
\Theta_3 & = & {\chi _{0, \,4}} + {\chi _{2, \,4}} +
{\chi _{2, \,0}}
 + {\chi _{4, \,0}} + {\chi _{3, \,4}} + {\chi _{3, \,2}} + {\chi
 _{4, \,2}} + {\chi _{5, \,0}}
\\
& & \mbox{}  + {\chi _{2, \,6}} + {\chi _{3, \,
0}} + {\chi _{3, \,6}} + {\chi _{4, \,6}} + {\chi _{5, \,2}} + {
\chi _{6, \,2}}  + {\chi _{1, \,6}} + {\chi _{1, \,4}}
\\
\Theta_4 & = & {\chi _{1, \,0}} + {\chi_{T, \,0}}
\\
\Theta_5 & = & {\chi _{0, \,1}} +
{\chi _{0, \,T}} + {\chi _{1, \,0}} + {\chi _{ 1, \,T}} + {\chi
_{T, \,0}} + {\chi _{T, \,1}}
\\
\Theta_6 & = & {\chi _{1, \,0}} + 2\,{\chi
_{0, \,2}} + {\chi _{0, \,1} } + 2\,{\chi _{2, \,1}} + {\chi _{1,
\,2}} + 2\,{\chi _{0, \,T}}
 + {\chi _{1, \,T + 1}} + {\chi _{T, \,0}} + 2\,{\chi _{T, \,1}}
 + {\chi _{T + 1, \,T}} \\
 & & \mbox{} + {\chi _{2 + T, \,T + 1}} + 2\,{\chi _{2 + T, \,T}}
 + 2\,{\chi _{2\,T, \,0}} + 2\,{\chi _{2 + T, \,1 + 2\,T}} + {
\chi _{2 + T, \,2\,T}} + {\chi _{2\,T, \,T + 1}} \\
 & & \mbox{} + {\chi _{2\,T, \,2}} + 2\,{\chi _{2\,T, \,1}} + {
\chi _{1 + 2\,T, \,0}} + 2\,{\chi _{0, \,2 + 2\,T}} + {\chi _{2\,
T, \,1 + 2\,T}} + {\chi _{0, \,2\,T}} + 2\,{\chi _{0, \,2 + T}}
 \\
 & & \mbox{} + 2\,{\chi _{0, \,T + 1}} + 2\,{\chi _{2\,T, \,2 + T
}} + {\chi _{1, \,2\,T}} + 2\,{\chi _{1 + 2\,T, \,T + 1}} + {\chi
 _{1 + 2\,T, \,T}} + {\chi _{1 + 2\,T, \,2}} \\
 & & \mbox{} + {\chi _{1 + 2\,T, \,1}} + 2\,{\chi _{1, \,2 + 2\,T
}} + {\chi _{2, \,2\,T}} + {\chi _{2, \,2 + T}} + 2\,{\chi _{2 +
2\,T, \,2}} + {\chi _{2, \,T + 1}} + {\chi _{2, \,T}} \\
 & & \mbox{} + {\chi _{2 + 2\,T, \,2 + T}} + 2\,{\chi _{T, \,2}}
 + 2\,{\chi _{T, \,2 + T}} + 2\,{\chi _{T, \,T + 1}} + 2\,{\chi
_{T, \,1 + 2\,T}} + {\chi _{T, \,2\,T}} + {\chi _{T + 1, \,1}} \\
 & & \mbox{} + {\chi _{T + 1, \,2}} + 2\,{\chi _{T + 1, \,2 + T}}
 + 2\,{\chi _{T + 1, \,2\,T}} + 2\,{\chi _{T + 1, \,1 + 2\,T}} +
2\,{\chi _{2 + T, \,0}} + {\chi _{2 + T, \,1}}
\\
\Theta_7& =& {\chi _{0, \,1}} + {\chi _{2,
\,1}} + {\chi _{0, \,T}}
 + {\chi _{1, \,T}} + {\chi _{2 + T, \,T}} + {\chi _{2 + T, \,2
 + 2\,T}} + {\chi _{2\,T, \,1}} + {\chi _{1, \,2 + 2\,T}} + {\chi
 _{2 + 2\,T, \,T}} + {\chi _{2, \,T}} \\
 & & \mbox{} + {\chi _{2 + 2\,T, \,1}} + {\chi _{2, \,2 + 2\,T}}
 + {\chi _{T + 1, \,2 + 2\,T}}
\\
 \Theta_8 & = &{\chi _{0, \,1}} + {\chi _{0, \,T}} + {\chi _{0, \,2 + 2
\,T}} + 2\,{\chi _{1, \,T}} + {\chi _{1, \,2 + 2\,T}} + 2\,{\chi
_{2, \,1}} + {\chi _{2, \,2 + 2\,T}} + {\chi _{T, \,1}} + 2\,{
\chi _{T, \,2 + 2\,T}} \\
 & & \mbox{} + {\chi _{T + 1, \,T}} + 2\,{\chi _{T + 1, \,2 + 2\,
T}} + 2\,{\chi _{2 + T, \,1}} + 2\,{\chi _{2 + T, \,T}} + 2\,{
\chi _{2 + T, \,2 + 2\,T}} + {\chi _{2\,T, \,1}} \\
 & & \mbox{} + 2\,{\chi _{2\,T, \,T}} + 2\,{\chi _{2 + 2\,T, \,1}
} + {\chi _{2 + 2\,T, \,T}} 
\\
\Theta_9 &=& {\chi _{0, \,1}} + {\chi _{0, \,5}} + {\chi _{1,
\,5}} + {\chi _{ 2, \,1}} + {\chi _{2, \,5}} + {\chi _{3, \,5}} +
{\chi _{5, \,1}}
 + {\chi _{7, \,1}}
\\
\Theta_{10} & = & {\chi _{0, \,1}} + {\chi _{0, \,4}} + {\chi _{0, \,7}}
+ {\chi _{ 3, \,1}} + {\chi _{3, \,4}} + {\chi _{3, \,7}} + {\chi
_{6, \,1}}
 + {\chi _{6, \,4}} + {\chi _{6, \,7}}
\end{eqnarray*}


\begin{center}
\begin{tabular}{|l|c|c||c|c|} \hline
quandle $Q$ & cocycle     & A & T(n,k) & invariant \\ \hline
\hline $R_4$ & $\Theta_1$& ${\bf Z}_2$ & $T(2,4)$ & $8+8t$
\\ \cline{4-5}
 & &  & $T(4,6)$ & $8+8t$ 
\\ \cline{4-5}
 & &  & $T(4,8)$ & $128 + 128t$
\\ \hline
$R_6$ & $\Theta_2$ & ${\bf Z}_2$ & $T(2,2)$  & $6 +6t$
\\ \cline{4-5}
 & &  & $T(2,6)$ & $18 + 18t$ 
\\ \cline{4-5}
 & &  & $T(4,4)$ &  $216 + 216t$ 
\\ \cline{4-5}
 & &  & $T(4,12)$ & $648 +648t$
 \\ \hline
 $R_8$ & $\Theta_3$ & ${\bf Z}_2$ & $T(2,8)$ & $32+32t$ 
\\ \cline{4-5}
 & &  & $T(4,16)$ & $2048 + 2048t$  
\\ \hline
 ${\bf Z}_2[T, T^{-1}]/(T^2+1)$ &  $\Theta_4$
 & ${\bf Z}_2$ & $T(2,4)$ & $8+8t$ 
\\ \cline{4-5}
 & &  & $T(4,6)$ & $8+8t$ 
\\ \cline{4-5}
 & &  & $T(4,8)$ & $128+128t$ 
\\ \hline 
${\bf Z}_2[T, T^{-1}]/(T^2+T+1)$ &  $\Theta_5$
 & ${\bf Z}_2$ & $T(2,3)$ & $4 +12t$ 
\\ \cline{4-5}
 & &  & $T(3,3)$ & $4 +12t$
\\ \cline{4-5}
 & &  & $T(3,6)$ & $16 +48t$ 
\\ \hline 
${\bf Z}_3[T, T^{-1}]/(T^2+1)$& $\Theta_6 $ & ${\bf Z}_3$ &
$T(2,4)$ & $9 + 36t+36t^2$ 
\\ \cline{4-5}
& &  & $T(4,4)$ & $297 +216t +216t^2$ 
\\ \cline{4-5}
 & &  & $T(4,8)$ & $297 +216t +216t^2$ 
\\ \cline{4-5}
 & & & $T(4,10)$ & $9 +36t +36t^2$ 
\\ \hline  ${\bf Z}_3[T, T^{-1}]/(T^2-1)$ & $\Theta_7$ & ${\bf Z}_2$
& $T(2,2)$ & $21 + 6t$ 
\\ \cline{4-5}
 & & & $T(2,6)$ & $63 + 18t$ 
\\ \cline{4-5}
 & & & $T(3,3)$ & $63 + 18t$
 \\ \cline{4-5}
 & & & $T(4,4)$ & $1647 + 540t$ 
\\ \cline{4-5}
 & & & $T(4,12)$ & $4941 + 1620t$ 
\\ \hline  
${\bf Z}_3[T, T^{-1}]/(T^2+T+1)$ & $\Theta_8$ & ${\bf
Z}_3$ & $T(3,3)$ & $135 +54t + 54t^2$ 
\\ \cline{4-5}
 & &  & $T(3,6)$ & $135 +54t + 54t^2$ 
\\  \hline  
 ${\bf Z}_8[T, T^{-1}]/(T-5)$ & $\Theta_9$ & ${\bf Z}_2$ &
$T(2,2)$ & $28 + 4t$ 
\\ \cline{4-5}
 & & & $T(2,4)$ & $48 + 16t$ 
\\ \cline{4-5}
 & & &  $T(2,6)$ & $28 + 4t$ 
\\ \cline{4-5}
 & & & $T(3,3)$ & $104 + 24t$ 
\\ \cline{4-5}
& & & $T(4,4)$ & $1600 + 448t$ 
\\ \cline{4-5}
 & & & $T(4,6)$ & $48 + 16t$ 
\\ \cline{4-5}
 & & & $T(4,6)$ & $3072 + 1024t$ 
\\ \hline  
${\bf Z}_9[T, T^{-1}]/(T-7)$ & $\Theta_{10}$ & ${\bf Z}_2$ &
$T(2,6)$ & $ 63+18t$ 
\\ \cline{4-5}
 & & & $T(3,3)$ & $ 81+162t$
 \\ \cline{4-5}
 & & & $T(3,9)$ & $567+162t$ 
\\ \cline{4-5}
 & & & $T(4,12)$ & $ 4941+1620t$ 
\\
\hline
 \end{tabular} \\[5mm]
{Table~\ref{torustable} : Nontrivial invariants of torus links}
\end{center}  


\section{Computations with Movies}
\label{defate}

\begin{figure}
\begin{center}
\mbox{
\epsfxsize=6in
\epsfbox{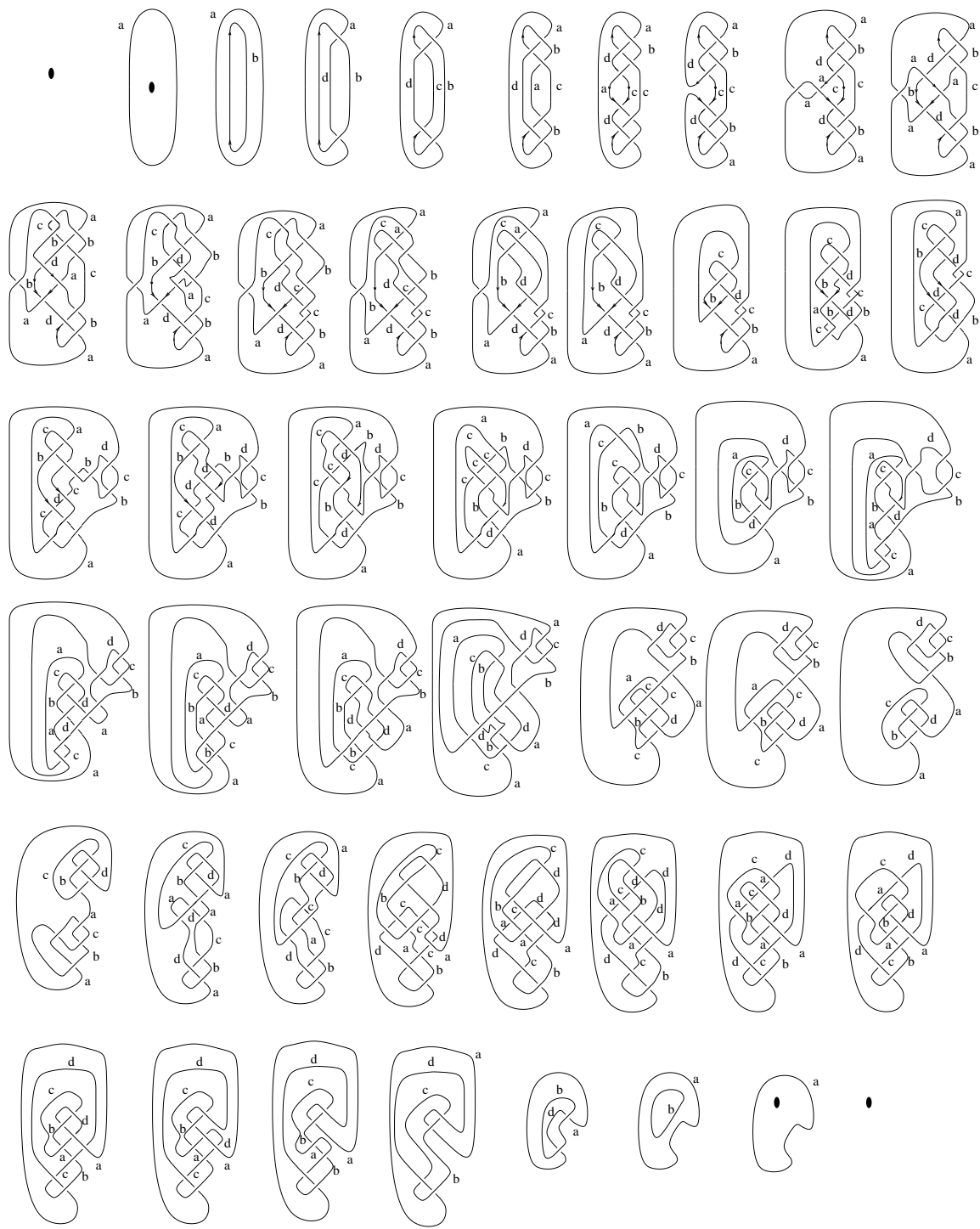}
} \end{center}
\caption{A deform spun figure eight knot     
}
\label{deffig8label}
\end{figure}

In this section we 
present a method of computing the cocycle invariants
using movies, with a sample calculation 
for a deform-spun figure-8 knot, 
the movie of which is illustrated in Fig.~\ref{deffig8label}.

First, we describe the knotted 
sphere that the movie illustrates. The 
sequence of {\it stills} consist of classical knot diagrams. 
Successive stills differ by 
a Reidemeister move of type  II or III, 
a critical point of the surface, or a planar 
isotopy of the underlying diagram.
Type I moves are not used in this particular example, though 
they may appear in general. 
Exceptions are made for maximal/minimal points (or the birth/death 
of small circles) where dots are depicted in some stills in the figure
to make clear where they occurred. 
In this case, the exact moment
of critical points are depicted in a still, instead of between 
successive stills.
The entries in the figure are referred to by a pair of 
numbers that indicate the row and column in the figure. 
Thus still (1,8) is a standard picture of 
the connected sum of a pair of figure-8 knots.

Stills (1,1), (1,2), (6,7), and (6,8)  
represent critical levels in which simple closed 
curves are born or die.
In particular, the still (1,1) consists of a single dot 
where a birth occurs. The still (1,2) consists of one circle
and a dot.

 These curves are colored $a$ 
and $b$ with the curve colored $b$ nested inside 
the curve colored $a$. 
 Other critical points (saddle points) occur 
between stills (1,7) and (1,8) and between 
(6,1) and (6,2). 
Since these are all the critical points the resulting surface depicted 
is a sphere. 

The labels indicate a coloring of the 
diagram by $S_4$.
It is seen that for any  
$a,b\in S_4$, if 
 $d=a*b$ and $c=b*d$, then such a choice gives a coloring.
In particular, there is a non-trivial coloring of this knotted surface
by $S_4$.
The quandle rule holds at each crossing 
point of the stills. As the crossings move 
from still to still, they sweep out double 
point arcs, and the quandle coloring rule 
holds along these double point arcs. 
The non-trivial  
colorings show that the quandle of the 
embedded sphere is non-trivial, and so 
the sphere is indeed knotted. 
Observe that the proof of knotting does not depend 
on an {\it a priori} knowledge that the sphere is 
deform-spun 
knot. 

Next we describe the story-board 
of the movie. 
Roughly the second, third, fourth rows each involve 
pushing the top figure-8 tangle past an arc of the 
bottom figure-8 tangle. The figure is arranged to fit 
on a page, and so at the end of the first and 
second rows the top knot has begun its next migration. 
Each of the four crossings of the top figure-8 tangle 
must cross under or over one of the four arcs of the bottom figure-8 tangle.
 When the crossings past these arcs triple points in the projection occur. 
The triple points appear as Reidemeister type III 
moves between the stills. There are 16 such triple 
points they occur in the scenes listed in the table below.

Most of the rest of the changes are 
type II Reidemeister moves.
The exceptions are $(2,6)\mapsto(2,7)$, 
where a rescaling has taken 
 place, 
and $(4,7) \mapsto (5,1) $,
where the bottom tangle has rotated clockwise 
around the left most arc to become the top tangle in $(5,1)$.

One can use the
movie to obtain the following presentation of the quandle of the
 knotted sphere.

$$\langle x,y: y*( x *\bar{y})=x*y; y*(y*(x* \bar{y})))=x \rangle.$$

\begin{figure}
\begin{center}
\mbox{
\epsfxsize=5in
\epsfbox{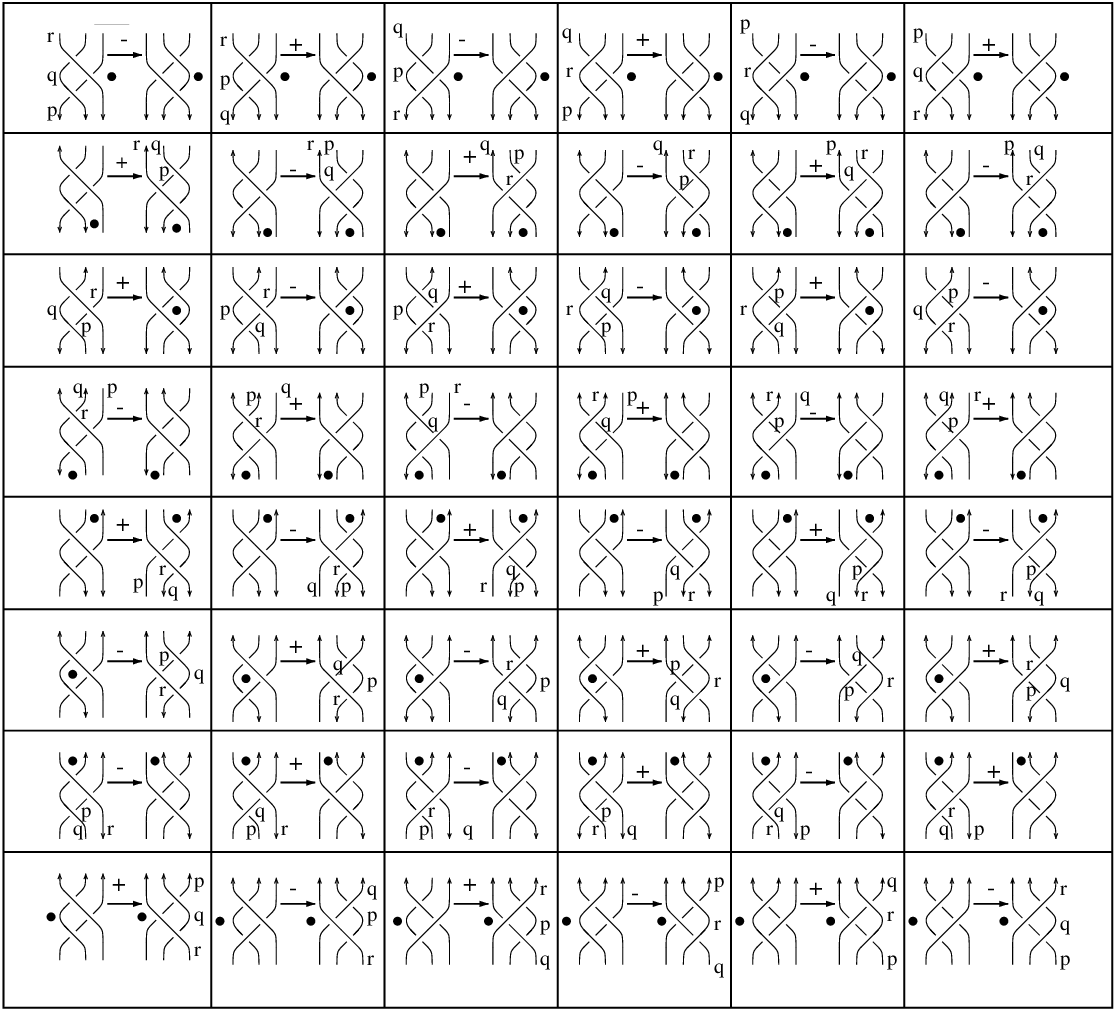}
}
\end{center}
\caption{All possibilities of oriented type III moves }
\label{alltypeIII}
\end{figure}

In order to evaluate the cocycle invariant for cocycles in $Z^3(S_4,A)$, 
we need to examine the signs of the triple points 
and the colors of the regions away from which the normal point.
This examination is easily made by comparing the given scenes in the movie
 with the  
Reidemeister type III moves that are illustrated in Fig.~\ref{alltypeIII}.
 This figure 
illustrates the 96 possible type III moves
 (48 in the direction of the arrows, 48 in the opposite direction). 
The black dots indicate the region into which all normals point when
 the normal points to the left of the oriented arc. 
The labels $(p,q,r)$ indicate the colors on the bottom, middle, 
and top (resp.) sheets that are the arguments of the quandle cocycles.
The entries in Fig.~\ref{alltypeIII}
 are indicated in a manner similar to those in the previous figure. 
 Thus in entry $(1,6)$ the braid move
$\sigma_1 \sigma_2 \sigma_1 
\mapsto \sigma_2 \sigma_1 \sigma_2$ 
is
indicated.

In the following table we indicate the scenes 
that correspond to triple points, the entry 
in Fig.~\ref{alltypeIII} to which the scenes correspond. 
The third column indicates whether the entry in the figure  
is read forward (left to right) or backward (right to left). 
If it is backward then the sign of the 
triple point is opposite to that indicated in the figure. 
 The last column indicates the Boltzmann weight associated to the triple point in the movie Fig.~\ref{deffig8label}. 

\begin{center}
\begin{tabular}{||c|c|c|c||} \hline \hline
Scene in movie        & type III & direction & Weight \\ \hline \hline
$(1,9) \mapsto (1,10)$& $(1,2)$ & $\leftarrow$ &  $\theta(b,a,a)^{-1}$  \\ \hline
$(2,1)\mapsto (2,2)$  & $(2,5)$ & $\rightarrow$ & $\theta(c,b,a)$ 
\\ \hline
$(2,2) \mapsto (2,3)$ & $(2,4)$ & $\rightarrow$ & $\theta(c,d,a)^{-1}$ 
\\ \hline
$(2,3)\mapsto (2,4)$  & $(1,4)$ & $\rightarrow$ &  $\theta(b,c,a)$ 
\\ \hline
$(2,8) \mapsto (2,9)$ & $(5,4)$ & $\rightarrow$ & $\theta(c,b,a)^{-1}$
\\ \hline
$(3,1) \mapsto (3,2)$ & $(6,1)$ & $\rightarrow$ & $\theta(d,c,d)^{-1}$ 
\\ \hline
$(3,2) \mapsto (3,3)$ & $(6,4)$ & $\rightarrow$ & $\theta(c,c,b)$ 
\\ \hline
 $(3,3) \mapsto (3,4)$& $(5,1)$ & $\rightarrow$ & $\theta(c,b,c)$ 
\\ \hline
$(4,1) \mapsto (4,2)$ & $(1,6)$ & $\rightarrow$ & $\theta(b,a,b)$ 
\\ \hline
$(4,2)\mapsto (4,3)$  & $(2,3)$ & $\rightarrow$ & $\theta(c,d,b)$ 
\\ \hline
$(4,3)\mapsto (4,4)$  & $(2,6)$ & $\rightarrow$ & $\theta(c,b,b)^{-1}$ 
\\ \hline
$(4,4)\mapsto (4,5)$  & $(1,3)$ & $\rightarrow$ & $\theta(b,c,b)^{-1}$ 
\\ \hline
$(5,2)\mapsto(5,3)$   & $(2,6)$ & $\leftarrow$ & $\theta(a,b,c)$ 
\\ \hline
$(5,4) \mapsto (5,5)$& $(4,5)$ & $\leftarrow$ & $\theta(d,d,a)$ 
\\ \hline
$(5,5)\mapsto(5,6)$   & $(2,4)$ & $\rightarrow$ & $\theta(c,d,b)^{-1}$ 
\\ \hline
$ (5,6)\mapsto (5,7)$ & $(1,1)$ & $\rightarrow$ & $\theta(c,b,d)^{-1}$ 
\\ \hline
\hline 
\end{tabular}
\end{center}

Recall from \cite{CJKLS} that the following are non-trivial 3-cocycles for $S_4$
with various coefficients.

\begin{eqnarray*}
\eta_{1}&=&
+\chi_{(0,1,0)}+\chi_{(0,2,1)}+\chi_{(0,2,3)}+\chi_{(0,3,0)}
+\chi_{(0,3,1)}+\chi_{(0,3,2)}+\chi_{(1,0,1)} \\ & & 
+\chi_{(1,0,3)}
+\chi_{(1,2,0)}+\chi_{(1,3,1)}+\chi_{(2,0,3)}+\chi_{(2,1,0)}
+\chi_{(2,1,3)}+\chi_{(2,3,2)}; \\
\eta_{2}&=&
+\chi_{(0,1,2)}-\chi_{(0,1,3)}-\chi_{(0,2,1)}+\chi_{(0,3,0)}
+\chi_{(0,3,1)}-\chi_{(0,3,2)}+2\chi_{(1,0,1)} \\
& & +\chi_{(1,0,2)} 
+\chi_{(1,0,3)}-\chi_{(1,2,0)} 
+\chi_{(1,3,2)}+\chi_{(2,0,1)}
+\chi_{(2,0,2)}+\chi_{(2,0,3)} \\ & & 
+\chi_{(2,1,3)}+\chi_{(3,0,1)}
+\chi_{(3,0,2)}+\chi_{(3,0,3)}+\chi_{(3,1,3)}; \\	
\eta_{11}&=&
-\chi_{(0,1,0)}-\chi_{(0,1,3)}+\chi_{(0,3,1)}
+\chi_{(0,3,2)}-\chi_{(1,0,1)}-\chi_{(1,0,2)}-\chi_{(1,0,3)} \\ & & 
+\chi_{(1,2,0)}-\chi_{(1,2,1)}+\chi_{(1,3,0)}+\chi_{(1,3,1)}
+\chi_{(1,3,2)}+\chi_{(2,0,3)}-\chi_{(2,1,0)}-\chi_{(3,0,2)}
+\chi_{(3,2,3)} 
\end{eqnarray*}

Thus we have
\begin{sect}{\bf Theorem.} The state sum invariants for the deform-spun figure-8 knot that is depicted in Fig.~\ref{deffig8label} is
$$\Phi_\theta =\sum  \theta(b,c,a)\theta(c,b,c)\theta(b,a,b) 
\theta(a,b,c)\theta(c,d,a)^{-1}\theta(d,c,d)^{-1}\theta(b,c,b)^{-1}\theta(c,b,d)^{-1}$$
where the sum is taken over all pairs of elements $a,b\in S_4$
and $a=d*b$ while $c=b*d$.
We have the values
\begin{itemize} 
\item 
$\Phi_{\eta_{11}} = 16.$
\item
$\Phi_{\eta_1} = 4 +12 t$ for $t$ a generator of ${\bf Z}_2.$
\item
$\Phi_{2\eta_1}= 4+12t^2$ for $t$ a generator of ${\bf Z}_4.$
\item
$\Phi_{\eta_2} = 4 +12 t$ for $t$ a generator of ${\bf Z}_2.$
\item
$\Phi_{\eta_2}= 4+12t$ 
 for $t$ a generator of ${\bf Z}_4.$
\end{itemize}
\end{sect}

\section{
Twist-spun Torus Knots} 
\label{moviesec}

Let $\tau^k T(p,q)$ denote the $k$-twist spun $(p,q)$-torus knot (or link).
The construction of twist spun knots 
was 
first defined by Zeeman \cite{Zeeman},
see also \cite{Rolf}.
In this section we 
use movies to 
compute the state-sum expressions 
of $\tau^k T(2,m)$ with 
cocycles over 
dihedral quandles. 
We prove 
the invariants are periodic  
with respect to $k$. 
Then we evaluate the state-sum for some  quandle cocycles with various 
coefficients.
First we establish some notation.

\begin{sect}{\bf Notation.\/} {\rm
Let $x,y$ denote elements of an Alexander quandle 
$\Lambda = {\bf Z}_p [T, T^{-1}]/ (f(T))  $; let $s$ denote an integer.
Define quandle elements $G(s)= G(x,y,s,T)$ recursively by
\begin{eqnarray*} 
G(-1) &=& x \\
G(0) & =& y  \\
G(s+1) &=& T G(s-1) + (1-T) G(s) \end{eqnarray*}
Then $G(-2) = T^{-1} y +(1-T^{-1})x$, and for $s\ge 0,$
$$G(s) = x \sum_{j=1}^s (-1)^{j+1} T^j + y \sum_{j=0}^s (-1)^j T^j.$$

Define $h(x,y,0)=y$, 
 and 
$h(x,y,n)= T^{-n}y + (1-T^{-n})x$.
Then $h$ satisfies the relation
$T h(x,y,n+1)+(1-T)x = h(x,y,n).$

If $\theta$ is a $3$-cocycle  for the quandle $\Lambda$, then define
$$\Theta^m_0(x,y, T) = \prod_{j=0}^{m-1} 
\theta(G(x,y,-2,T), G(x,y,j-1,T), G(x,y,j,T))^{-1}.$$
and 
$$\Theta^m_1(x,y,T) = \prod_{j=0}^{m-1}
 \theta(G(x,y,j-2,T), G(x,y,j-1,T), G(x,y,-2,T)).$$
}\end{sect}

\begin{figure}
\begin{center}
\mbox{
\epsfxsize=3in
\epsfbox{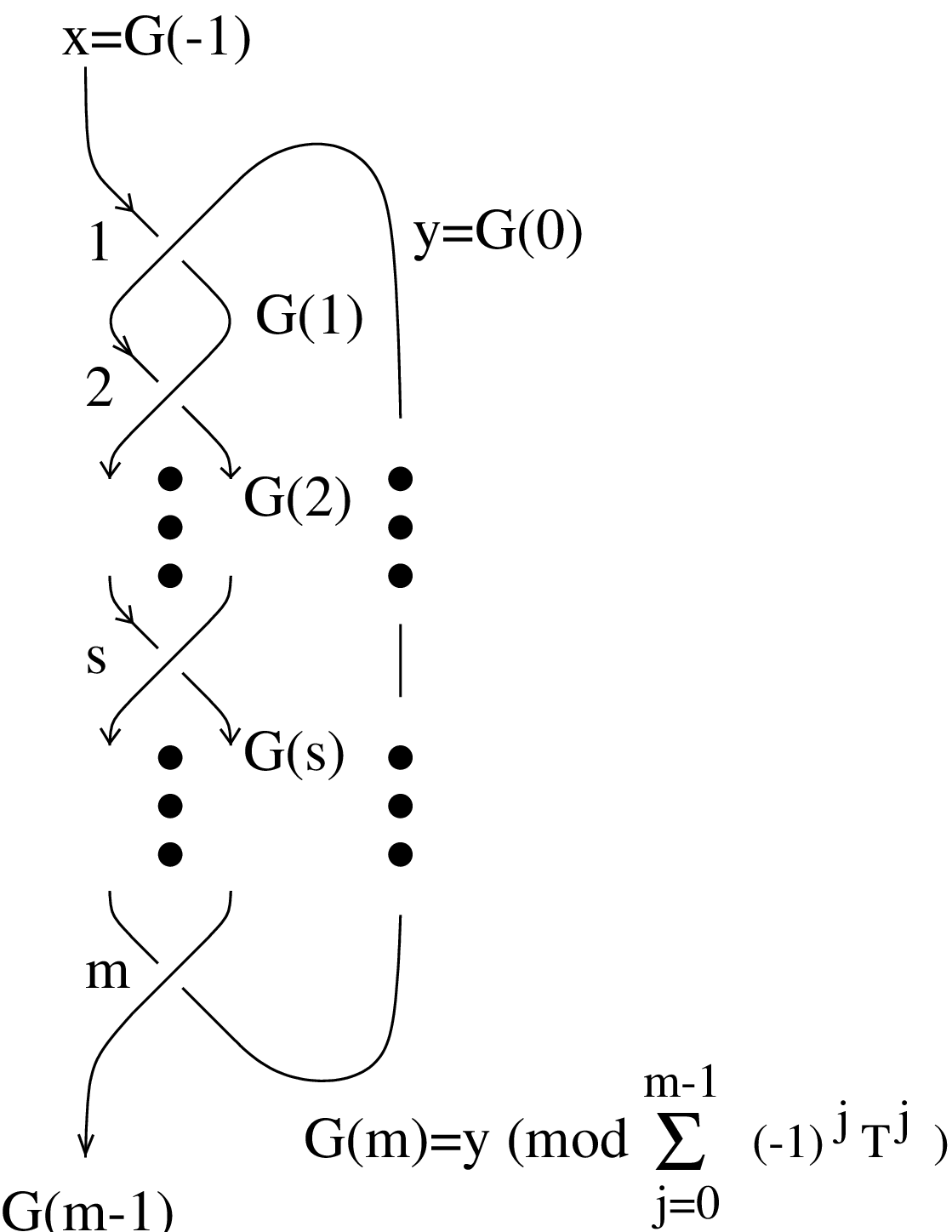}
}
\end{center}
\caption{Alexander colors for $\tau^k T(2,m)$ }
\label{twist1LX}
\end{figure}

\begin{figure}
\begin{center}
\mbox{
\epsfxsize=3in
\epsfbox{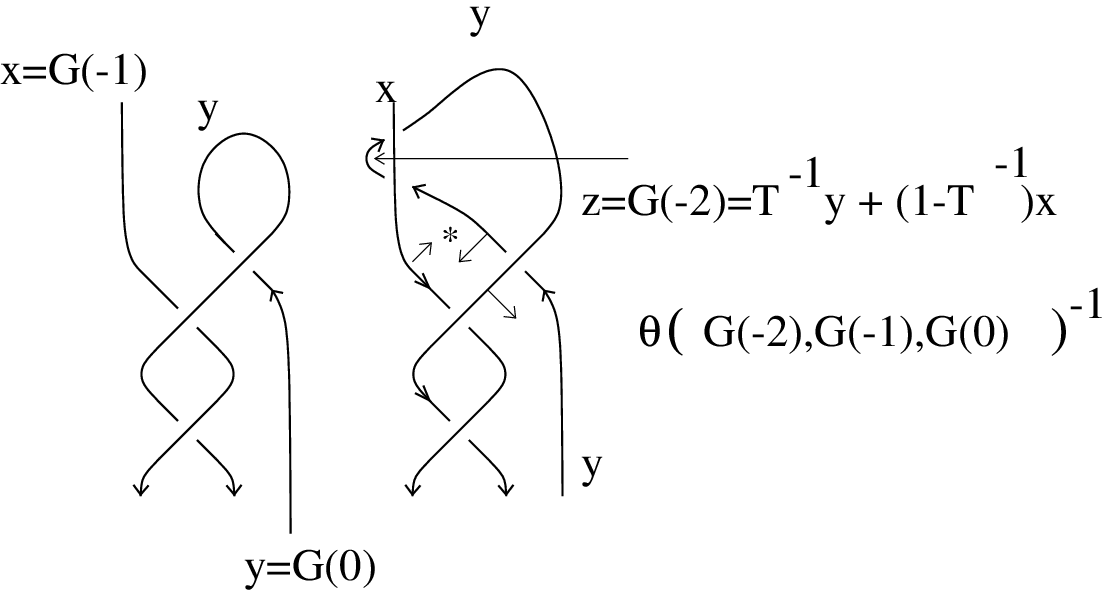}
}
\end{center}
\caption{Start twisting on top }
\label{twist2LX}
\end{figure}

\begin{figure}
\begin{center}
\mbox{
\epsfxsize=4.5in
\epsfbox{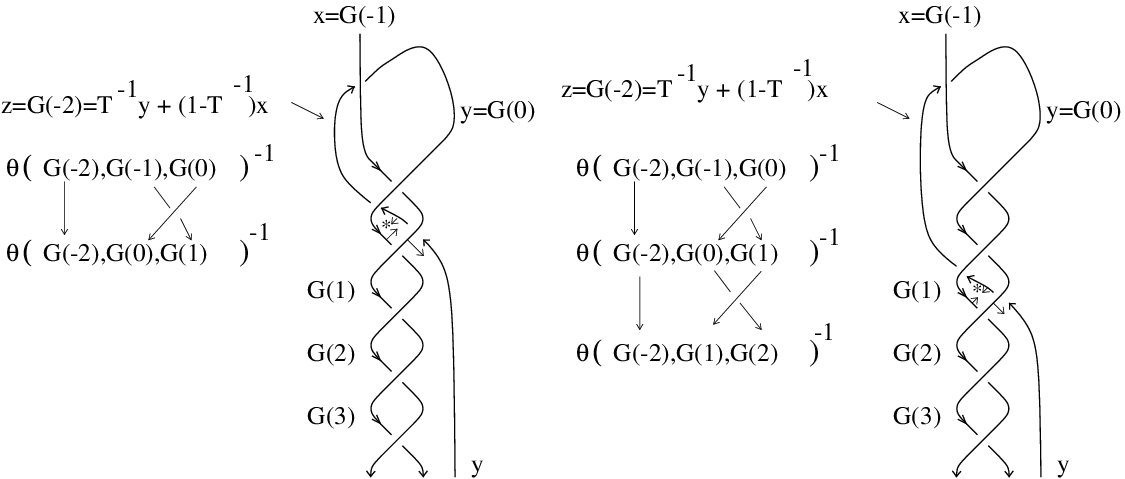}
}
\end{center}
\caption{Continue twisting by Reidemeister III moves }
\label{twist3LX}
\end{figure}

\begin{figure}
\begin{center}
\mbox{
\epsfxsize=2in
\epsfbox{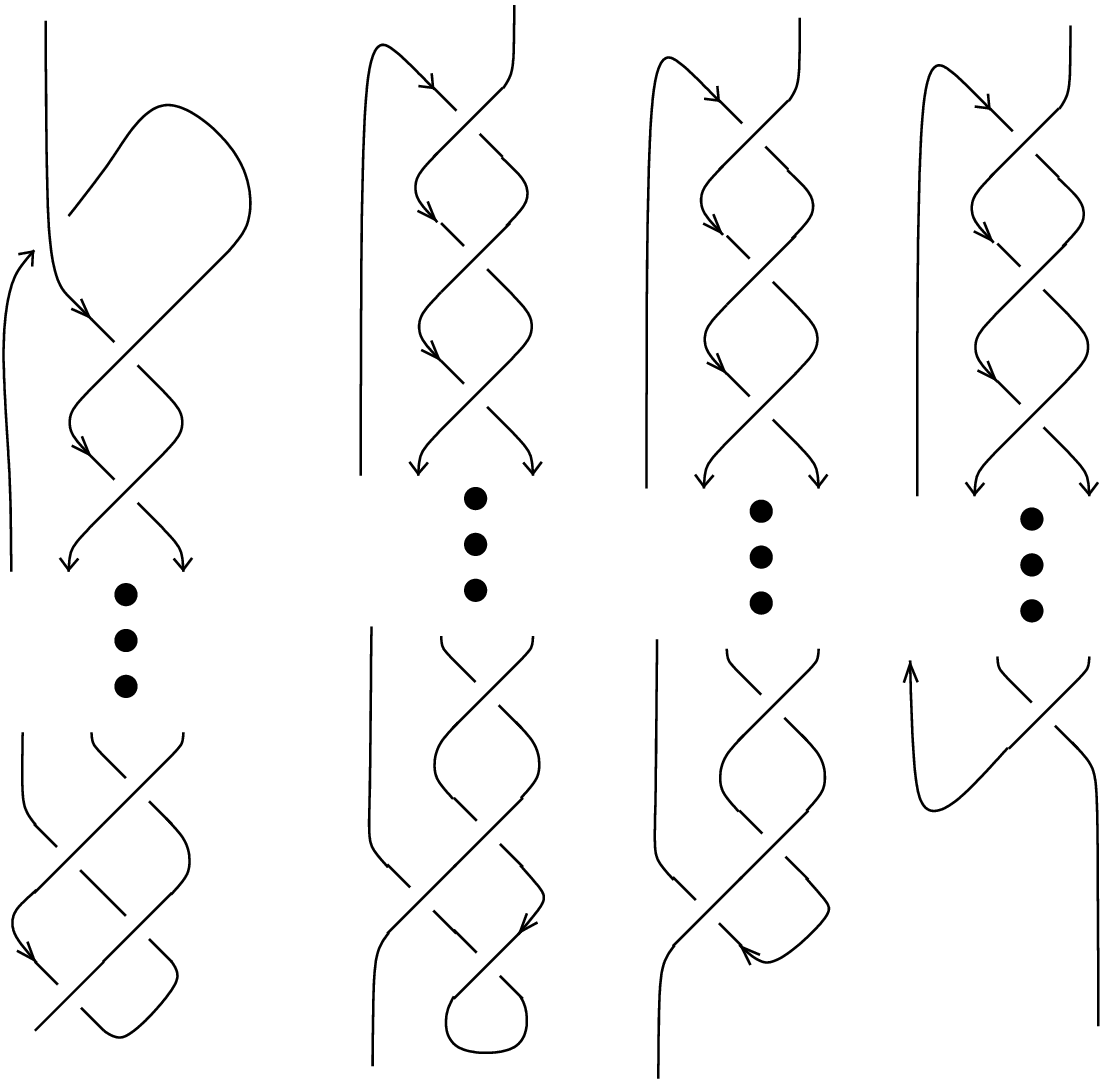}
}
\end{center}
\caption{Complete half a twist}
\label{twist4LX}
\end{figure}

\begin{figure}
\begin{center}
\mbox{
\epsfxsize=3in
\epsfbox{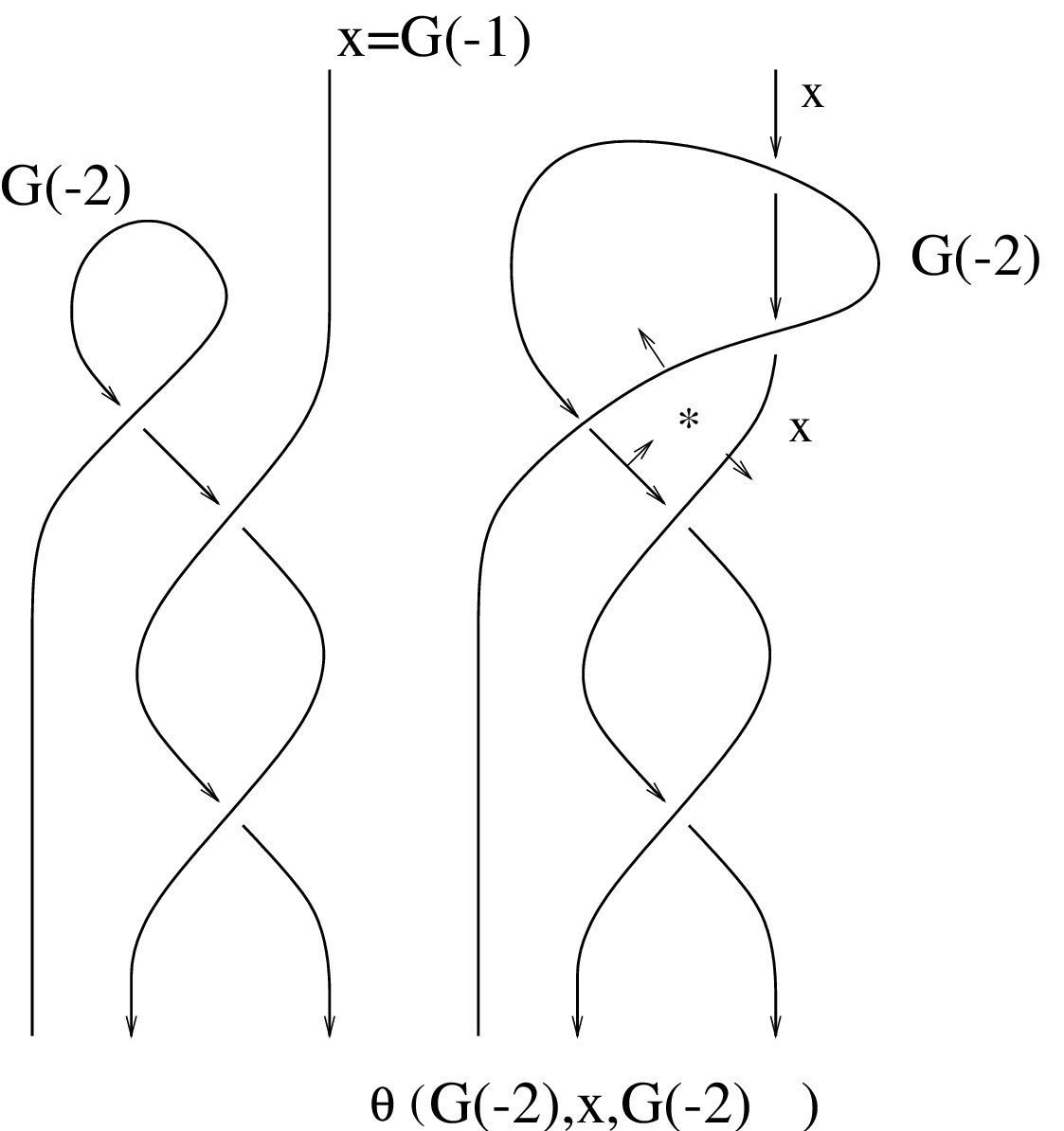}
}
\end{center}
\caption{Start twising again on top }
\label{twist5LX}
\end{figure}

\begin{figure}
\begin{center}
\mbox{
\epsfxsize=4.5in
\epsfbox{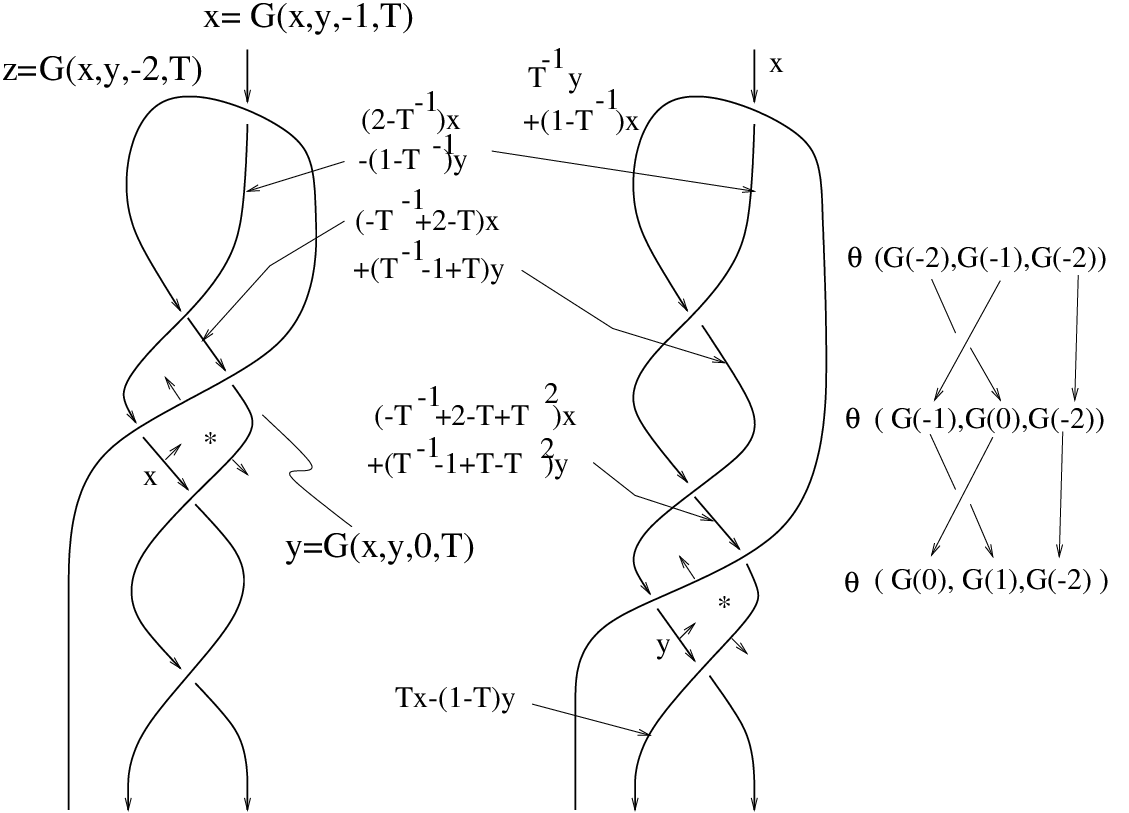}
}
\end{center}
\caption{Keep twisting by Reidemeister III moves }
\label{twist6LX}
\end{figure}

\begin{figure}
\begin{center}
\mbox{
\epsfxsize=3in
\epsfbox{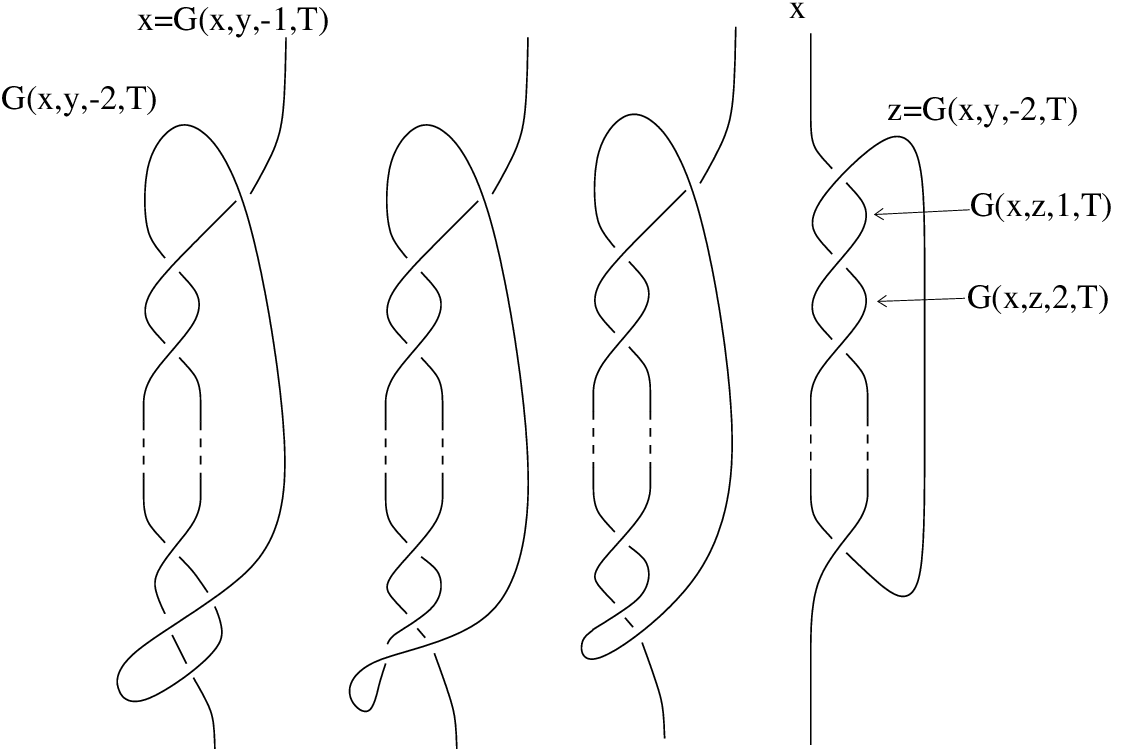}
}
\end{center}
\caption{One full twist completed }
\label{twist7LX}
\end{figure}

Next we describe the important scenes in the movie
 of the twist-spun torus knots,
$\tau^kT(2,m)$. The tangle depicted in Fig.~\ref{twist1LX} is
 the tangle of $T(2,m)$ upon
which we will apply the twisting. 
The initial stages of the first half-twist are indicated in 
Figs.~\ref{twist1LX} and \ref{twist2LX}.
 The final stages of the half-twist are indicated in
Fig.~\ref{twist3LX}. In Fig.~\ref{twist1LX} 
quandle elements
 are indicated on the center arcs
of the torus tangle. After the $s\/$th crossing the arc 
on the right of the crossing 
is colored $G(s)$. 
The normal direction to the surface is chosen to be a right pointing
 arrow on arcs that run downward. 
On the right-hand side of  Fig.~\ref{twist2LX}
 a new arc with color $G(-2)$ is born via a type II Reidemeister move.
 The asterisk in the triangle indicates the location of a type III
 Reidemeister move that is about to take place.
 Compare this to Fig.~\ref{alltypeIII} entry (2,6). 
This and the subsequent type III Reidemesiter moves during the
 first half-twist are all 
of this type; the signs are all negative. 
The arguments for the cocycles are indicated in Fig.~\ref{twist2LX}
where the relationship between the arguments on subsequent type III moves
 is indicated by the crossed arrows. The second half-twist is indicated
 in Figs.~\ref{twist4LX} through \ref{twist7LX}.
 Here the type III move appears in Fig.~\ref{alltypeIII} as entry (6,6).

\begin{sect} {\bf Theorem.\/}
The surface $\tau^k T(2,m)$, with the orientation shown in Fig.~\ref{twist1LX},
colors nontrivially with 
the Alexander  rack $\Lambda = {\bf Z}_p [T, T^{-1}]/ (f(T))  $
 if and only if
$1-T+ \cdots \pm T^{m-1} = 0$  and $T^k =1$ in $\Lambda$.

The state-sum  invariant
of $\tau^k T(2,m)$ with a cocycle $\theta$ 
of the Alexander quandle $\Lambda$  is 
$$\sum_{x,y} \prod_{n=0}^{k-1} \Theta^m_0 (x, h(x,y,n),T)
 \Theta^m_1(x,h(x,y,n),
T)$$
\end{sect}
{\it Proof.\/} As discussed above,
Fig.~\ref{twist1LX} through Fig.~\ref{twist7LX} depict
one full twist of a diagram of $T(2,m)$
by a sequence of Reidemeister moves. 
First, Fig.~\ref{twist1LX} shows that $T(2,m)$ is non-trivially
colored by $\Lambda$ if and only if 
$G(s)=G(x,y,m,T)= x \sum_{j=1}^m (-1)^{j+1} T^j + y \sum_{j=0}^m (-1)^j T^j=y$ if and only if
$\sum_{j=1}^m (-1)^{j+1} T^j = 0$
in $\Lambda$. 
After one full twist, the color on the right-hand arc of the figures changes
 to $h(x,y,1)= T^{-1}y+(1-T^{-1})x$. By induction, the color on 
the right-hand arc 
after $k$ full twists is $h(x,y,k)=T^{-k}y+(1-T^{-k})x$
Hence $\tau^k T(2,m)$ is non-trivially colored if and only if $T^k=1$ 
in $\Lambda$.
Note also that the colors on $\tau^k T(2,m)$ are periodic with 
period $k$ if $T^k=1$. 
 

%

%

%

After one full twist the Reidemeister type III moves have
 contributed factors of
$\Theta^m_0(x,y,T)\Theta^m_1(x,y,T)$ to a term in the state-sum, and the color 
$y=h(x,y,0)$ on the right-hand arc has changed to $h(x,y,1)$ when
 the twist is complete.
The result follows. $\Box$.

\begin{sect} {\bf Corollary.\/}
Suppose $1-T+ \cdots \pm T^{m-1}=0 $ and $T^n=1$ in the given 
Alexander rack $\Lambda$. Let $\theta$ be a $3$-cocycle
in $Z^3(\Lambda, {\bf Z}_q)$. 
Then 
 $\Phi_{\theta}(\tau^k T(2,m))$
is periodic 
with respect to $k$ with period  $nq$. 
\end{sect}
{\it Proof.\/} 
The colors have period $n$ and the state-sum takes values in ${\bf Z}_q$,
so the state-sum has period $nq$.
$\Box$

\begin{sect} {\bf Corollary.\/}
The surface $\tau^k T(2,m)$ colors nontrivially with 
the dihedral rack $R_h$ if and only if $h=m$ and $2|k$.
The state-sum invariant of 
$\tau^{2w} T(2,m)$ with a cocycle $\theta$ 
of a dihedral rack $R_m$ is
\begin{eqnarray*} 
 \Phi_{\theta}(\tau^{2w} T(2,m)) & = & \sum_{x,y}
\prod_{n=0}^{k-1} \Theta^m_0 (x, h(x,y,n),-1)
 \Theta^m_1(x,h(x,y,n),
-1) \\
& =& \sum_{x,y} \left[ \Theta_0^m (x,y,-1) \Theta_1^m (x,y,-1) \right] 
\left[ \Theta_0^m (x,2x-y,-1) \Theta_1^m (x,2x-y,-1) \right] \\
& \times &\left[ \Theta_0^m (x,y,-1) \Theta_1^m (x,y,-1) \right] 
\left[ \Theta_0^m (x,2x-y,-1) \Theta_1^m (x,2x-y,-1) \right] \\
& \times & \cdots {\mbox{\rm a total of $w$ factors}} \cdots \\
& \times& \left[ \Theta_0^m (x,y,-1) \Theta_1^m (x,y,-1) \right] 
\left[ \Theta_0^m (x,2x-y,-1) \Theta_1^m (x,2x-y,-1) \right] 
\end{eqnarray*}

Furthermore, if the coefficient of the cohomology is $A={\bf Z}_q$,
i.e. if $\theta \in C^2_Q(R_m,{\bf Z}_q)$,
then $\Phi_{\theta}(\tau^k T(2,m))$ 
is periodic 
with respect to $k$ with period  $2q$. 
\end{sect}
{\it Proof.\/} 
Evaluate $h(x,y,n)$ at $T=-1$ to obtain 
$$(-1)^{2j} y + (1-(-1)^{2j}) x = y,$$
$$(-1)^{2j+1} y + (1-(-1)^{2j+1}) x = 2x-y.$$
The terms of the state-sum have the form
$\Theta^w$
so that it has period $q$. Since the color has period $2$,
the invariant has period $2q$. 
$\Box$

\begin{sect} {\bf Table.\/} \label{tautable} {\rm                 
For the rest of the section we present {\sc Maple} computations
of the state-sum invariants for twist spun $T(2, m)$ torus knots
using the above formulas. 
Computations are summarized in the following table. 
The Alexander ideal of $\tau^k T(2,m)$ is
$(\Delta_m=T^{m-1}- \cdots \pm 1, T^k-1)$. Somewhat simplified ideals 
are shown in the table. 
Since $\Delta_m$ divides $T^{2m}-1$ if $m$ is odd and 
$T^{m}-1$ if $m$ is even, the Alexander modules are periodic 
with respect to $k$. Thus we listed up to the smallest period in the table.

The quandles we used are those whose $3$-cocycles were computed 
in Section~\ref{cocysec} with non-trivial cohomological dimension.
Often there is a pattern that the dimensions are the same for 
several prime numbers of $q$. In that case, the smallest among such 
is chosen. 

Since the dihedral quandles $R_m$ 
and other quandles of the same order appear often as quandles
that color $\tau^k T(2,m)$ non-trivially, we listed knots with $m<7$,
as {\sc Maple} does not finish computations for larger quandles
in predictable time. 

The blank entries means that either only trivial quandles color, 
or no quandles color. In these cases the invariant column is also 
left blank. 

For $\tau^4 T(2,m)$, the invariant follows from 
periodicity, since $\tau^2 T(2,m)$ colors with the same quandle 
$R_m$ 
In this case, the same cocycle was used for both knotted surfaces.
When the same quandles are used for various $k$, we applied the 
periodicity in this way. Such cases are marked by $(p)$ in the table.
The first few examples of such values are provided in the table,
but the others are left for the reader.

When there are more than one choice of cocycles that give different 
values of the invariant, some
different  choices are made by experiments.
In particular, we do not know those listed exhaust all possibilities.
The cocycles of choice are presented with the table. 

If, for example, $3+6t$ is a value of the invariant with $A={\bf Z}_3$,
then $3+6t^2$ appears also as a value, by taking the negative of
the cocycles of the former value. Such cases are listed but cocycles 
are not given in duplicate. 
} \end{sect}

The choices of cocycles are as follows.
The indices in front of cocycles represent indices in the table:
$(m-k-A-b)$ represents the invariant for $\tau^k T(2, m)$, 
the first (represented by $(-A-)$) choice of quandle/coefficient,
the second (represented by $(-b)$) choice of a cocycle.

\begin{eqnarray*}
(3-2-A) & &
 2\,{\chi _{0, \,1, \,2}} + 2\,{\chi _{0, \,2, \,1}} + 2\,{\chi _{
1, \,0, \,1}} + 2\,{\chi _{1, \,0, \,2}} + 2\,{\chi _{2, \,0, \,1
}} +2 {\chi _{2, \,0, \,2}}
 \\
(3-3-A-a) & & 
{\chi _{0, \,T, \,1 + T}} 
 +{\chi _{0, \,1 + T, \,1}}
+ {\chi _{1, \,0,\,1}} + {\chi _{1, \,0, \,T}} +{\chi _{1, \,0,\,1 + T}}
+  {\chi _{1, \,T, \,0}} + {\chi _{1, \,T, \,1}} +
{\chi _{1, \,1 + T, \,0}}
\\ & & \mbox{}
 + {\chi _{T, \,0, \,1}}  +
{\chi _{T, \,0, \,T}} + {\chi _{T, \,1, \,0}} + {\chi _{T, \,1 + T, \,0}}
\\
(3-3-A-b) & & 
 {\chi _{1, \,T, \,1 + T}} + {\chi _{1, \,T, \,0}} + {\chi _{1, \,T, \,1}}
 +{\chi _{1, \,0, \,1}}
 +  {\chi _{0, \,1 + T, \,0}} +{\chi _{0, \,T, \,0}} +
{\chi _{0, \,T, \,1}} + {\chi _{0, \,1, \,T}}
\\ & & \mbox{}
 + {\chi _{T, \,1 + T, \,0}} +
{\chi _{T, \,0, \,T}} + {\chi _{T, \,0, \,1 + T}} + 
 {\chi _{1 + T, \,T, \,1 + T}} + {\chi _{1 + T, \,0, \,T}} +
 {\chi _{1 + T, \,0, \,1 + T}} + {\chi _{1 + T, \,0, \,1}}
\\
(4-2-A-a) & & 
{\chi _{0, \,1, \,0}}  +  {\chi _{0, \,1, \,3}} +
{\chi _{0, \,2, \,1}} +{\chi _{0, \,2, \,3}} +{\chi _{0, \,3, \,0}} 
 + {\chi _{0, \,3, \,1}}
\\
(4-2-A-b) & & 
 {\chi _{0, \,3, \,1}} + {\chi _{0, \,2, \,3}} + {\chi _{0, \,3, \,2}}
 + {\chi _{0, \,2, \,0}} + {\chi _{2, \,0, \,1}} 
 + {\chi _{2, \,1, \,3}} + {\chi _{2, \,1, \,0}}
 + {\chi _{2, \,0, \,2}} 
\\
(4-2-B-a) & &
{\chi _{0, \,1, \,0}} + {\chi _{0, \,1,\,T}} +{\chi _{0, \,T, \,0}}
 + {\chi _{0, \,T, \,1}}
 + {\chi _{0, \,1 + T, \,1}} + {\chi _{0, \,1 + T, \,T}} 
 \\
  (4-2-B-b) & &
{\chi _{0, \,1, \,0}} + {\chi _{0, \,1,\,T}} +{\chi _{0, \,T, \,0}}
 + {\chi _{0, \,T, \,1}}
  + {\chi _{0, \,1 + T, \,1}} +
{\chi _{0, \,1 + T, \,T}} +{\chi _{1, \,0, \,1}}
 + {\chi _{1, \,0, \,1 + T}}
\\ & & \mbox{}
 + {\chi _{1, \,T, \,0}}
 + {\chi _{1, \,T, \,1 + T}} + {\chi _{1, \,1 + T, \,0}}
 + {\chi _{1, \,1 + T, \,1}} 
\\
(5-2-A) & & 
3\,{\chi _{3, \,0, \,2}} + 
4\,{\chi _{3, \,0, \,1}} +
 3 \,{\chi _{3, \,0, \,3}} +
{\chi _{3, \,0, \,4}}  +
{\chi _{3, \,4, \,0}} + 
 4 \,{\chi _{3, \,4, \,3}} 
 + 
2\,{\chi _{3, \,1, \,2}} +
{\chi _{3, \,1, \,3}} 
\\ & &  +
 2 \,{\chi _{4, \,0, \,1}}
+ 3\,{\chi _{4, \,0,\,2}} +
2\,{\chi _{3, \,2, \,1}} +
 3\,{\chi _{4, \,0, \,3}} 
+ {\chi _{4, \,0, \,4}}
  +2\,{\chi _{0, \,1, \,0}}  +
 4\,{\chi _{0, \,1, \,2}} +
 2 \,{\chi _{0, \,1, \,3}} 
\\ & & 
 + 2\,{\chi _{0, \,1, \,4}} +
 2 \,{\chi _{0, \,2, \,0}} 
+
 3\,{\chi _{1, \,0, \,1}} +
{\chi _{0, \,3, \,1}} +{\chi _{0, \,3, \,2}} +
3\,{\chi _{0, \,4, \,0}} +{\chi _{0, \,4, \,1}}
4\,{\chi _{1, \,0, \,2}}
\\ & & 
 + 
3\,{\chi _{1, \,0, \,3}} +
 4\,{\chi _{3, \,2, \,4}} + 
4\,{\chi _{0, \,4, \,3}} + 
3 \,{\chi _{0, \,4, \,2}} +
 2\,{\chi _{1, \,0, \,4}} +
2\,{\chi _{1, \,3, \,2}} +
 2 \,{\chi _{1, \,3, \,1}} +
2 \,{\chi _{1, \,2, \,4}}
\\ & & 
 +
3\,{\chi _{1, \,3, \,0}} +
3\,{\chi _{1, \,2, \,1}} +
2\,{\chi _{1, \,2, \,3}} +
3\,{\chi _{1, \,2, \,0}} +
4\,{\chi _{1, \,3, \,4}} +
3\,{\chi _{1, \,4, \,0}} +
4\,{\chi _{1, \,4, \,1}} +
{\chi _{2, \,0, \,2}}
\\ & & 
 + {\chi _{2, \,0, \,3}} 
 +3 \,{\chi _{2, \,1, \,2}}+
2\,{\chi _{2, \,1, \,3}} +
4 \,{\chi _{2, \,1, \,0}}+
2\,{\chi _{2, \,3,\,0}} +
2\,{\chi _{2, \,1, \,4}} +
 2 \,{\chi _{2, \,3, \,4}} +
2\,{\chi _{2, \,4, \,0}}
\\ & & 
 +
3\,{\chi _{2, \,4,\,1}} +
3\,{\chi _{2, \,3, \,2}} +
4\,{\chi _{2, \,3, \,1}} + {\chi _{2, \,4, \,2}} + {\chi _{4, \,2, \,1}}
\\
(6-2-A) & & 
{\chi _{0, \,1, \,2}} +
{\chi _{0, \,2, \,1}} +
{\chi _{1, \,0, \,1}} + 
{\chi _{1, \,0, \,2}} +
{\chi _{2, \,0, \,1}} +
 {\chi _{2, \,0, \,2}} 
\\
(6-2-B-a) & & 
{\chi _{2, \,3, \,1}}
 +{\chi _{2, \,3, \,4}}
 +2\,{\chi _{2, \,5, \,0}}
 +2\,{\chi _{2, \,5, \,3}} +
 {\chi _{4, \,0, \,2}} +
 {\chi _{4, \,0, \,5}} +
 {\chi _{4, \,1, \,2}}
 + {\chi _{4, \,1, \,5}}
\\ & & 
 +{\chi _{4, \,1, \,3}} +
{\chi _{4, \,2, \,0}}+
 {\chi _{4, \,2,\,3}} +
{\chi _{4,\,1, \,0}}
{\chi _{2, \,1, \,3}} +
 {\chi _{2, \,1, \,5}} 
{\chi _{2, \,1, \,2}} +
2\,{\chi _{0, \,1, \,2}}
\\ & & 
+ 2\,{\chi _{0, \,1, \,5}} +
 2\,{\chi _{0, \,3, \,2}}
 2\,{\chi _{0, \,3, \,5}}
 +2\,{\chi _{0, \,2, \,4}}
+2\,{\chi _{0, \,2, \,1}}
+ {\chi _{0, \,4, \,2}}
+2\,{\chi _{0, \,4, \,3}} + 
2\,{\chi _{0, \,4, \,0}}
\\ & &  +
 2\,{\chi _{0, \,5, \,0}}
+2\,{\chi _{0, \,5, \,1}}
+ {\chi _{0, \,4, \,5}} +
 2\,{\chi _{0, \,5, \,4}}
+2\,{\chi _{0, \,5, \,3}} + 
{\chi _{2, \,0, \,1}} 
+{\chi _{2, \,0, \,2}} +
{\chi _{2, \,0, \,5}}
\\ & & +
{\chi _{2, \,0, \,4}} 
+{\chi _{2,\,1, \,0}} 
\\
(6-2-B-b) & & 
{\chi _{5, \,2, \,3}} +
2\,{\chi _{1, \,2, \,3}} +
 2\,{\chi _{1, \,4, \,0}} +
 2\,{\chi_{1, \,4, \,3}} 
+{\chi _{1, \,3, \,2}} + 
{\chi _{1, \,3, \,5}} +
2\,{\chi _{1, \,4, \,2}} +
 2\,{\chi _{1, \,4, \,5}}
\\ & & +
{\chi _{1, \,5, \,0}} +
{\chi_{1, \,5, \,3}} +
{\chi _{2, \,0, \,1}} +
{\chi _{2, \,0, \,2}} +
{\chi _{2, \,0, \,5}} +
{\chi _{2, \,0, \,4}} +
{\chi _{2, \,1, \,0}}
 +{\chi _{2, \,1, \,2}}
\\ & & 
+ {\chi _{2, \,1, \,3}} +
{\chi _{2, \,1, \,5}} +
{\chi _{2, \,3, \,1}} +
{\chi _{2, \,3, \,4}} +
2\,{\chi _{2, \,5, \,0}} +
 2\,{\chi _{2, \,5, \,3}} +
{\chi _{3, \,0, \,2}} +
{\chi _{3, \,1, \,0}}
\\ & &  +
{\chi _{3, \,0, \,5}} +
{\chi _{3, \,1, \,3}}+
{\chi _{3, \,1, \,2}} +
{\chi _{3, \,1, \,5}} +
{\chi _{3, \,2, \,0}} +
{\chi _{3, \,2, \,3}} +
2\,{\chi _{1, \,0, \,5}} +
2\,{\chi _{1, \,0, \,2}}
\\ & &  + 
2\,{\chi _{0, \,5, \,1}} +
2\,{\chi _{0, \,5, \,3}} +
{\chi _{4, \,0, \,2}}+
2\,{\chi _{0, \,5, \,4}}+
{\chi _{0, \,4, \,2}}+
2\,{\chi _{0, \,4, \,3}}+
2\,{\chi _{0, \,4, \,0}} +
{\chi _{0, \,4, \,5}}
\\ & &  +
 2\,{\chi _{0, \,5, \,0}} +
\,{\chi_{4, \,0, \,5}}+
 {\chi _{4, \,1, \,2}} + 
{\chi _{4, \,1, \,5}} +
{\chi _{4, \,2, \,0}} +
{\chi _{4, \,1, \,0}} +
{\chi _{4, \,1, \,3}}+
 {\chi _{4, \,2, \,3}}
\\ & &  +
{\chi _{5,\,0, \,2}} +
{\chi _{5, \,1, \,3}}+
{\chi _{5, \,1, \,2}} +
 {\chi _{5, \,1, \,5}} +
 {\chi _{5, \,2, \,0}} + 
{\chi _{5, \,1, \,0}} +
{\chi _{5, \,0, \,5}} +
2\,{\chi _{1, \,2, \,0}}
\\ & & +
2\,{\chi _{0, \,1,\,2}} +
2\,{\chi_{0, \,1, \,5}} +
2\,{\chi _{0, \,2, \,1}} +
2\,{\chi_{0, \,2, \,4}} +
 2\,{\chi _{0, \,3, \,2}}
\\
(6-2-B-c) & & 
2\,{\chi _{5, \,0, \,4}} +
 2\,{\chi _{5, \,1, \,0}}+
{\chi _{5, \,0, \,5}} +
{\chi _{5, \,1, \,2}} +
{\chi _{5, \,1, \,5}}+
2\,{\chi _{5, \,1, \,3}} + 
{\chi _{5, \,2,\,1}} +
{\chi _{5, \,2, \,4}}
\\ & &  +
 {\chi _{5, \,3, \,1}} +
{\chi _{5, \,0, \,2}} +
2\,{\chi _{0, \,1,\,2}} +
2\,{\chi_{0, \,1, \,5}}+
2\,{\chi _{0, \,2, \,1}} +
2\,{\chi _{0, \,2, \,4}}+
2\,{\chi _{0, \,4, \,0}}+
 2\,{\chi _{0,\,3, \,2}}
\\ & &  +
 2\,{\chi _{0, \,3, \,5}} +
{\chi _{0, \,4, \,5}} +
2\,{\chi _{0, \,5, \,0}}+ 
2\,{\chi _{0, \,4, \,3}}+
{\chi _{0, \,4, \,2}} +
2\,{\chi _{0, \,5, \,1}} +
2\,{\chi _{0, \,5, \,3}}+
2\,{\chi _{0, \,5, \,4}}
\\ & &  +
{\chi _{1, \,0, \,2}}+
{\chi _{1, \,0, \,1}} +
{\chi _{1, \,0, \,5}} +
2\,{\chi _{1, \,0, \,4}} +
2\,{\chi _{1, \,2, \,0}} +
2\,{\chi _{1, \,2, \,3}} +
{\chi _{1, \,2,\,4}}+
2\,{\chi _{1, \,3, \,1}}
\\ & &  + 
{\chi _{1, \,4,\,2}} +
{\chi _{1, \,4, \,5}} +
 2\,{\chi _{1, \,5, \,1}} + 
{\chi _{2, \,0, \,1}} + 
{\chi _{2, \,0, \,4}} +
{\chi _{2, \,0, \,2}}+ 
{\chi _{2, \,0, \,5}} +
{\chi _{2, \,1, \,0}}
\\ & & + 
{\chi _{2, \,1, \,3}} +
{\chi _{2, \,1, \,2}}+
{\chi _{2, \,1, \,5}}+
{\chi _{2, \,3, \,1}} +
{\chi _{2, \,3, \,4}} + 
2\,{\chi _{2, \,5, \,0}} +
 2\,{\chi _{2, \,5, \,3}} +
{\chi _{3, \,0, \,2}}
\\ & &  +
{\chi _{3, \,0, \,1}} +
{\chi _{3, \,0, \,5}} +
2\,{\chi _{3, \,0, \,4}} +
{\chi _{3, \,1, \,0}}+
2\,{\chi _{3, \,1, \,2}} +
{\chi _{3, \,1, \,3}} +
2\,{\chi _{3, \,2, \,0}} +
 2\,{\chi _{3, \,1, \,5}}
\\ & &  +
2\,{\chi _{3, \,2, \,3}} +
 2\,{\chi _{3, \,2, \,1}}+
{\chi _{3, \,2, \,4}} +
2\,{\chi _{3, \,4, \,0}}+
{\chi _{3, \,4, \,2}} +
2\,{\chi _{3, \,4, \,3}} +
{\chi _{3, \,4, \,5}} +
{\chi _{3, \,5, \,1}}
\\ & &  +
{\chi _{4, \,0, \,2}} +
 {\chi_{4, \,1, \,2}} + 
{\chi _{4, \,0, \,5}} +
{\chi _{4, \,1, \,0}} +
{\chi _{4, \,1, \,3}} +
 {\chi _{4, \,1, \,5}}+
 {\chi _{4, \,2, \,3}} +
{\chi _{4, \,2, \,0}} 
\end{eqnarray*}

\vspace*{-1cm}
\begin{center}
{Table~\ref{tautable} :   Cocycle invariants of twist spun torus knots $\tau^k T(2,m)$}
\begin{tabular}{|l|l||l|l|l|} \hline 
Torus & $k$ & Alexander & Color & Invariants \\
links &     & ideals    & quandles &        \\  \hline \hline 
$T(2,3)$ & 2 & $(T+1, 3) (\cong R_3)$ 
		 & $Z^3(R_3, {\bf Z}_3)$ & $3+6t$ $\; (3-2-A)$
\\ \cline{5-5}
& & & & $3+6t^2 $ 
\\ \cline{2-5} 
& 3 & $(T^2-T+1, 2)$
    & $Z^3({\bf Z}_2[T, T^{-1}]/(T^2+T+1),{\bf Z}_2) $ &
 $ 10 + 6 t  \; (3-3-A-a)$
\\ \cline{5-5}
& & & &  $8+8t \; (3-3-A-b)$
\\ \cline{2-5} 
& 4 & $(T+1, 3) (\cong R_3)$
	 & $Z^3(R_3, {\bf Z}_3)$ &  $3+6t^2, \; 3+6t \;  (p) $ 
\\ \cline{2-5} 
& 5 & $(1)$ &  &  
\\ \cline{2-5} 
& 6 & $(T^2-T+1)$ &  $Z^3(R_3, {\bf Z}_3)$ & $9$ for any cocycle $ \; (p) $  
\\ \cline{4-5} 
&   &             & $Z^3({\bf Z}_2[T, T^{-1}]/(T^2+T+1),{\bf Z}_2) $ &
 $4$ for any cocycle $  \; (p) $ 
\\ \hline  
$T(2,4)$ & 2 & $(T^2-1, 2T-2)$ &   $Z^3(R_4, {\bf Z}_2)$ & 
 $ 12+4t \; (4-2-A-a)$ 
 \\ \cline{5-5}
& & & & $ 8 + 8t \;  (4-2-A-b)$ 
\\ \cline{4-5} 
& & &  $Z^3(R_4, {\bf Z}_3)$ & 16 
\\ \cline{4-5} 
& & & $Z^3({\bf Z}_2[T, T^{-1}]/(T^2+1),{\bf Z}_2) $ & 
 $12+4t \; (4-2-B-a)  $ 
\\ \cline{5-5}
& & & & $8+8t \; (4-2-B-b)  $ 
\\ \cline{4-5} 
& & & $Z^3({\bf Z}_2[T, T^{-1}]/(T^2+1),{\bf Z}_3) $ & 
$16 $ for any cocycles
\\ \cline{2-5} 
&3& $(1)$ & &  
\\ \cline{2-5} 
&4& $(\Delta_4)$ & & 
\\ \hline 
$T(2,5)$ & 2 & $(T^2-1, 2T-3)$&  $Z^3(R_5, {\bf Z}_5)$ & 
$5+10t+10t^4 \;  (5-2-A) $ 
\\ \cline{5-5}
& & & & $5+10t^2+10t^3 $ 
\\ \cline{2-5} 
&3& $(1)$  & &  
\\ \cline{2-5} 
&4& $(T+1, 5) (\cong R_5)$ &  $Z^3(R_5, {\bf Z}_5)$ & (p)
\\ \cline{2-5} 
&5& $(\Delta_5, 2)$ &  & 
\\ \cline{2-5} 
&6&  $(T+1, 5) (\cong R_5)$ &  $Z^3(R_5, {\bf Z}_5)$ & (p)
\\ \cline{2-5} 
&7&  $(1)$  & &  
\\ \cline{2-5} 
&8& $(T^2+5, 2T-3)$ & $Z^3(R_5, {\bf Z}_5)$ & (p)
\\ \cline{2-5} 
&9&  $(1)$  & &  
\\ \cline{2-5} 
&10& $(T^2-1, 2T-3)$ &  $Z^3(R_5, {\bf Z}_5)$ &  (p)
\\ \hline 
$T(2,6)$ & 2& $(T^2-1, 3T-3)$ & 
 $Z^3(R_3, {\bf Z}_3)$ & 
$3+6t \; (6-2-A)$ 
\\ \cline{5-5}
& & & & $3+6t^2 $
\\ \cline{4-5} 
& & &  $Z^3(R_6, {\bf Z}_2)$ & $36$ for any cocycle 
\\ \cline{4-5} 
& & &  $Z^3(R_6, {\bf Z}_3)$ & 
$24 + 12t \; (6-2-B-a) $ 
\\ \cline{5-5}
 & & & & $ 24+12t^2 $ 
\\ \cline{5-5}
 & & & & $ 12+24t \; (6-2-B-b) $
 \\ \cline{5-5}
 & & & & $ 12+24t^2 $ 
\\   \cline{5-5}
 & & & & $ 12+12t+ 12 t^2 \; (6-2-B-c) $  
\\ \cline{2-5} 
&3& $(T^3-1)$ & $Z^3({\bf Z}_2[T, T^{-1}]/(T^2+T+1),{\bf Z}_2) $ & $16$ for any cocycle 
\\ \cline{2-5}
&4& $(T^2-1, 3T-3)$  & 
 $Z^3(R_3, {\bf Z}_3)$ & (p) 
\\ \cline{4-5} 
& & &  $Z^3(R_6, {\bf Z}_2)$ & (p) 
\\ \cline{4-5} 
& & &  $Z^3(R_6, {\bf Z}_3)$ &  (p)
\\ \cline{2-5}
&5&  $(1)$  & &  
\\ \cline{2-5}
&6& $(\Delta_6)$ &  $Z^3(R_6, {\bf Z}_3)$ & (p) 
\\ \cline{4-5} 
& & &  $Z^3({\bf Z}_2[T, T^{-1}]/(T^2+T+1),{\bf Z}_2) $ & (p)  
\\ \hline 
\end{tabular}
\end{center}

\clearpage

\section{Computations with Surface Braids}
\label{sfcebrsec}

In this section, we use surface braid theory
 to provide computations for $\tau^2 T(2,m)$.
 Some of the results here coincide with those of
 the previous section, but there are advantages to 
the current approach that we outline. 
First, by making computations in more than one context, 
we are able to cross-check our results. 
Second, the braid chart provides a snap-shot of the entire 
knotted surface whereas in a movie description more work is
 needed to reconstruct 
the diagram of the knotting. 
Third, presentations for the fundamental quandle and the form of
 the partition function can be read directly from the braid chart.  
Thus the techniques we present for the current example are applicable
 in general to surfaces that are given in braid form.

We begin with a  brief review of the theory of surface braids;
 see \cite{Kam:top,Kam:ribbon,CS:book} for more details.

\begin{sect}{\bf Definition.\/} {\rm 
Let  $D^2$  and  $D$  be 2-disks and  $X_m$  a
fixed set of $m$  interior points of  $D^2$.
By  $pr_1 : D^2 \times D \to D^2$  and
$pr_2 : D^2 \times D \to D$, we mean the
projections to the first factor and to the
second factor.

A {\em surface braid\/} 
(\cite{Kam:ribbon}, \cite{Rud})
 of degree  $m$ is a
compact, oriented surface  $S$  properly
embedded in  $D^2 \times D$  such that the
restriction of  $pr_2$  to  $S$  is a degree-$m$
simple branched covering map and  $\partial S
= X_m \times \partial D \subset D^2 \times \partial D$.
A degree-$m$ branched covering map $f : S \to D$
is {\em simple\/} if
$| f^{-1}(y) | = m$ or $m-1$ for $y \in D$.
In this case, the branch points are simple ($z \mapsto z^2$).

A surface braid  $S$  of degree  $m$  is extended
to a closed surface  $\widehat{S}$  in  $D^2 \times
S^2$  such that
$\widehat{S} \cap (D^2 \times D) = S$ and
$\widehat{S} \cap (D^2 \times \overline{D}) =
X_m \times \overline{D}$,
where  $S^2$  is the 2-sphere obtained from  $D^2$
by attaching a 2-disk  $\overline{D}$  along
the boundary.
    By identifying  $D^2 \times S^2$  with the
tubular neighborhood of a standard 2-sphere
in ${\bf R}^4$, we assume that  $\widehat{S}$  is a
closed oriented surface embedded in  ${\bf R}^4$.
We call it the {\em closure\/}  of $S$  in ${\bf R}^4$.
It is proved in 
\cite{Kam:top} 
that every closed oriented surface embedded in ${\bf R}^4$ is
ambient isotopic to the closure of a surface braid.

Two surface braids  $S$  and  $S'$  in $D^2 \times D$
are said to be {\em equivalent\/}  if there is an
isotopy  $\{ h_t \}$  of  $D^2 \times D$  such that
\begin{enumerate}
\item
$h_0 = {\rm id}$, $h_1(S) = S'$,
\item
for each $t \in [0,1]$, $h_t$ is fiber-preserving;
that is, there is a homeomorphism
$\underline{h}_t : D \to D$ with
$\underline{h}_t \circ pr_2 = pr_2 \circ h_t$, and
\item
for each $t \in [0,1]$, $h_t |_{D^2 \times \partial D}
= {\rm id}$.
\end{enumerate}

Let  $C_m$  be the configuration space of unordered
$m$  interior points of  $D^2$.  We identify the
fundamental group  $\pi_1(C_m, X_m)$  of
$C_m$  with base point  $X_m$  with the  
braid group  $B_m$  
on $m$ strings. 
    Let  $S$  
denote 
a surface braid and $\Sigma(S) \subset D$  
the
branch point set of the branched covering map
$S \to D$.
For a path
$ a : [0,1] \to D - \Sigma(S)$,  we define a path
$$ \rho_S (a) : [0,1] \to C_m  $$
by
$$ \rho_S (a) (t) =
pr_1( S \cap (D^2 \times \{ a(t) \})). $$
    If  $pr_1( S \cap (D^2 \times \{ a(0) \})) =
pr_1( S \cap (D^2 \times \{ a(1) \})) = X_m$,
then the path  $\rho_S (a)$ represents an element
of  $\pi_1(C_m, X_m) = B_m$.
    Take a point  $y_0$ in $\partial D$.
The {\em braid monodromy\/}  of  $S$  is the
homomorphism
$$ \rho_S : \pi_1( D - \Sigma(S), y_0 ) \to B_m $$
such that  $\rho_S ([a]) = [ \rho_S (a) ]$
for any loop  $a$  in  $D - \Sigma(S)$  with
base point  $y_0$. 
}\end{sect}

\begin{sect}{\bf Definition.\/} {\rm
Let $\Sigma (S) = \{y_1, \ldots, y_n \}$. Take a regular neighborhood 
$N(\Sigma (S)) = N(y_1)\cup \cdots \cup N(y_n)$ in $D$.
A {\em Hurwitz arc system\/}  
${\cal A} = ( \alpha_1, \dots, \alpha_n )$  for $\Sigma(S)$
is an $n$-tuple of simple arcs in 
$E(\Sigma (S)) = {\mbox{\rm cl }}(D \setminus N(\Sigma (S))$ such that each 
$\alpha_i$ starts from a point of $\partial N (y_i)$ and ends at $y_0,$ and 
$\alpha_i \cap \alpha_j = \{ y_0 \}$  for
$ i \neq j $,  and $\alpha_1, \ldots, \alpha_n$
appear in this order around  $y_0$.

Let $\gamma_i$ ($i= 1, \dots, n$)  be 
the 
loop
$\alpha_i^{-1} \cdot \partial N(y_i)\cdot \alpha_i$ 
in $D - \Sigma(S)$  
with base point  $y_0$
which goes along $\alpha_i$, turns around the
endpoint of  $\alpha_i$ in the positive direction,
and returns along  $\alpha_i$.
The {\em braid system\/}  of $S$  associated with
${\cal A}$ is an  $n$-tuple of  $m$-braids
$$
( \rho_S ([\gamma_1]), \rho_S ([\gamma_2]), \dots,
\rho_S ([\gamma_n]) ). $$
} \end{sect}

\begin{sect}{\bf Definition.\/} {\rm
An {\it $m$-chart} \cite{Kam:ribbon} 
is oriented, labelled graph $\Gamma$ in $D$, which may
be empty or have closed edges without vertices called
{\em hoops\/}, satisfying the following conditions:
\begin{enumerate}
\item
Every vertex has degree one, four or six.
\item
The labels of edges are in $\{ 1, 2, \dots, m-1 \}$.
\item
For each degree-six vertex, three consective edges are
oriented inward and the other three are outward,
and these six edges are labelled  $i$  and  $i+1$
alternately for some  $i$.
\item
For each degree-four vertex, diagonal edges have the
same label and are oriented coherently, and
the labels  $i$  and  $j$  of the diagonals satisfy
$| i-j | > 1$.
\end{enumerate}
We call 
a 
degree 1 (resp. degree 6) vertex a
{\em black\/} (resp. {\em white\/}) vertex. 
A degree 4 vertex is called a {\it crossing point of 
the chart}.

For an $m$-chart $\Gamma$, we consider a surface braid
$S$ of degree $m$ satisfying the following conditions:
\begin{enumerate}
\item
For a regular neighborhood  $N(\Gamma)$  of $\Gamma$ 
in $D$ and 
for any  $y \in D - {\rm int}N(\Gamma)$
the projection is                       
 $pr_1( S \cap (D^2 \times \{ y \}) ) = X_m$,  
where $X_m$ denotes the $m$ fixed interior points of $D^2$. 

\item
The branch point set of  $S$  coincides with the set
of the black vertices of  $\Gamma$.
\item
For a path  $\alpha: [0,1] \to D$  which is in general
position with respect to $\Gamma$  and  $\alpha(0)$,
$\alpha(1)$ are in  $D - {\rm int}N(\Gamma)$,
the $m$-braid determined by  $\rho_S (\alpha)$  is
the  $m$-braid presented by the intersection braid
word  $w_\Gamma (a)$.
\end{enumerate}
Then we call  $S$  a {\em surface braid described
by $\Gamma$\/}.
}\end{sect}

\begin{figure}
\begin{center}
\mbox{
\epsfxsize=5in
\epsfbox{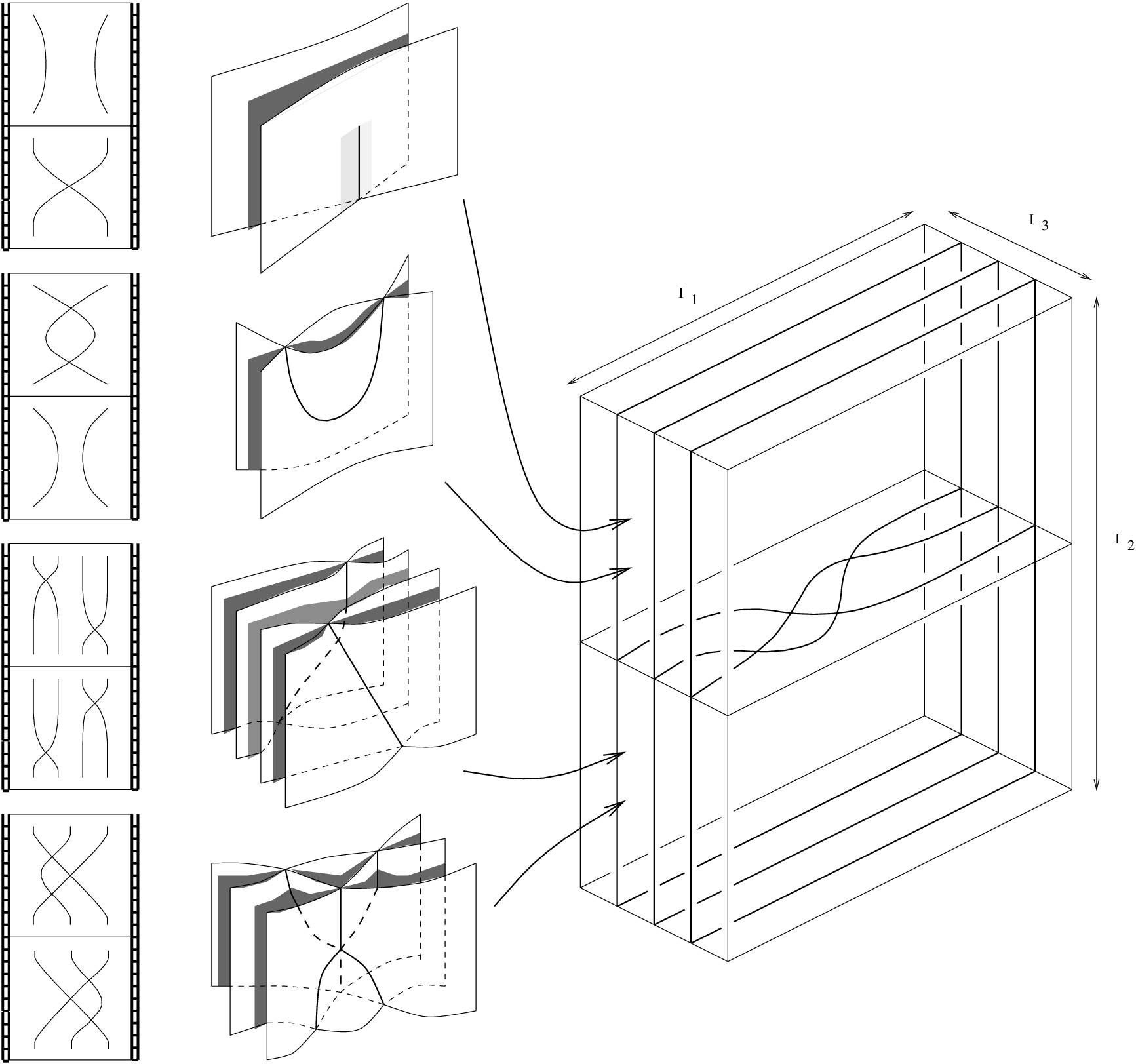}
} \end{center}
\caption{Projections and charts         
}
\label{ebcproj}
\end{figure}

In general, the singularity set of the image of $S$ by the
projection  $I_1 \times I_2 \times I_3 \times I_4
\to I_2 \times I_3 \times I_4$  is identified 
naturally
with the chart $\Gamma$ 
in the sense of \cite{Kam:Nato} and \cite{CS:book}.
The white vertices are in one-to-one correspondence to
the triple points and the black vertices 
are to
the branch points.  Figure~\ref{ebcproj} shows the 
relationship schematically, see
\cite{Kam:Nato} and \cite{CS:book} 
 for details.

\begin{figure}
\begin{center}
\mbox{
\epsfxsize=5in
\epsfbox{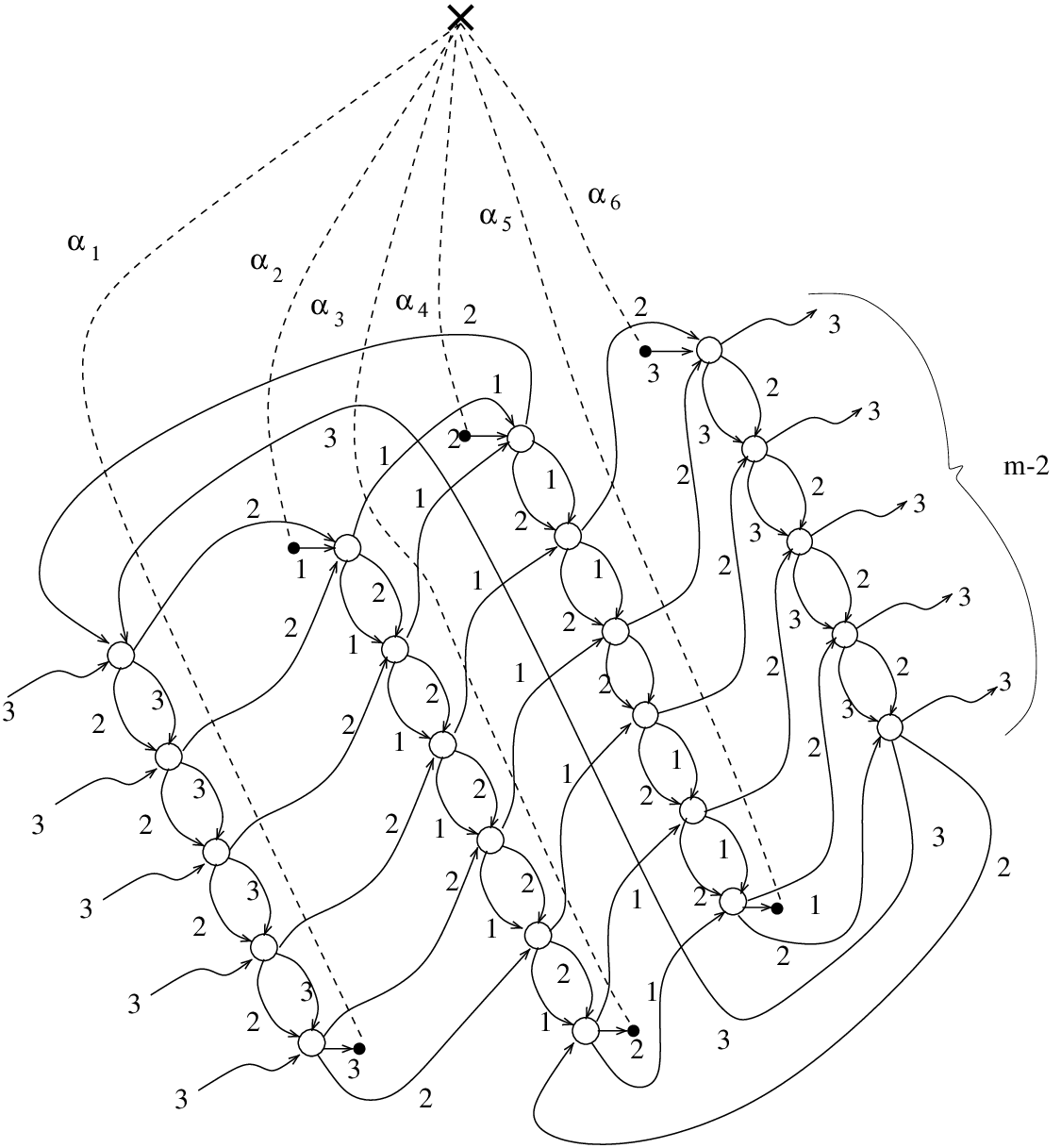}
}
\end{center}
\caption{A braid chart and Hurwitz system for 
the 2-twist spun $T(2,m)$: Here $m=7$}
\label{taut2malph}
\end{figure}

Let $m$ be an integer with $m \geq 3$ and let
$\Gamma$ be the $4$-chart illustrated in
Fig.~\ref{taut2malph}. 
 Let $S$ be the
surface braid of degree $4$ described by $\Gamma$.
The closure of $S$  in $R^4$ is the surface link $\tau^2(T(2,m))$
which is obtained from
$T(2,m)$, the torus knot/link of type $(2,m)$,
by Zeeman's 2-twist spinning.
(If $m$ is odd, then it is a 2-knot.  If $m$ is even,
then it consists of 2 components; one is a 2-sphere and
the other is a torus.  If $m=0$, $1$, or $2$, then $\tau^2(T(2,m))$
is a ribbon surface link and every quandle cocycle invariant of it is
trivial.  So we assume $m \geq 3$ in this paper.)

\begin{sect} {\bf Theorem.\/} 
For a quandle 3-cocycle $\theta
\in Z^3_Q(X;A)$, the state sum
$\Phi_\theta ( \tau^2(T(2,m)) )$
is
\begin{eqnarray*}
  \sum_{y_1,y_2}
&&
\prod_{j : {\rm odd} ;  1\leq j \leq m-2}
\theta (
y_2 \ast (y_1 y_2)^{(j-1)/2}y_1^{-1},
y_1 \ast (y_2 y_1)^{(j+1)/2}y_1^{-2}, 
y_1)^{-1} \\
&\times&
\prod_{j : {\rm even} ;  1\leq j \leq m-2}
\theta (
y_1 \ast (y_2 y_1)^{j/2}y_1^{-2},  
y_2 \ast (y_1 y_2)^{j/2}y_1^{-1},  
y_1)^{-1} \\
&\times&
\prod_{j : {\rm odd} ;  1\leq j \leq m-1}
\theta (
y_1 \ast (y_2 y_1)^{(j-5)/2}y_2,
y_2 \ast (y_1 y_2)^{(j-3)/2},
y_1 \ast y_2 y_1^{-1})^{-1} \\ 
&\times&
\prod_{j : {\rm even} ;  1\leq j \leq m-1}
\theta (
y_2 \ast (y_1 y_2)^{(j-4)/2},
y_1 \ast (y_2 y_1)^{(j-4)/2}y_2,
y_1 \ast y_2 y_1^{-1})^{-1} \\ 
&\times&
\prod_{j : {\rm odd} ;  1\leq j \leq m-1}
\theta (
y_1 \ast y_2 y_1^{-1},  
y_2 \ast (y_1 y_2)^{(j-3)/2},
y_1 \ast (y_2 y_1)^{(j-3)/2}y_2)^{+1} \\
&\times&
\prod_{j : {\rm even} ;  1\leq j \leq m-1}
\theta (
y_1 \ast y_2 y_1^{-1},
y_1 \ast (y_2 y_1)^{(j-4)/2}y_2,
y_2 \ast (y_1 y_2)^{(j-2)/2})^{+1} \\
&\times&
\prod_{j : {\rm odd} ;  1\leq j \leq m-2}
\theta (
y_1,
y_2 \ast (y_1 y_2)^{(j-1)/2},
y_1 \ast (y_2 y_1)^{(j-1)/2}y_2)^{+1} \\
&\times&
\prod_{j : {\rm even} ;  1\leq j \leq m-2}
\theta (
y_1,
y_1 \ast (y_2 y_1)^{(j-2)/2}y_2,
y_2 \ast (y_1 y_2)^{j/2})^{+1},
\end{eqnarray*}
where  $y_1, y_2$ run over all elements of  $X$ satisfying
$y_2 \ast (y_1 y_2)^{(m-1)/2} = y_1$ (if $m$ is odd)
or $y_1 \ast (y_2 y_1)^{m/2} = y_1$ (if $m$ is even), and
$y_1 \ast y_2^2 = y_1$.
\end{sect}

The formula is restated as follows.

\begin{sect} {\bf Theorem.\/} 
Let $m=2n+1$ (resp. $m=2n$).
For a quandle 3-cocycle $\theta
\in Z^3_Q(X;A)$, the state sum
$\Phi_\theta ( \tau^2(T(2,m)) )$
is

\begin{eqnarray*}
  \sum_{y_1,y_2}
&&
\prod_{k=1}^{n ({\rm resp.}~ n-1)}
\theta (
y_2 \ast (y_1 y_2)^{k-1}y_1^{-1},
y_1 \ast (y_2 y_1)^{k}y_1^{-2},  
y_1)^{-1} \\
&\times&
\prod_{k=1}^{~~~~n-1~~~~ }
\theta (
y_1 \ast (y_2 y_1)^{k}y_1^{-2}, 
y_2 \ast (y_1 y_2)^{k}y_1^{-1},  
y_1)^{-1} \\
&\times&
\prod_{k=1}^{~~~~~n~~~~~ }
\theta (
y_1 \ast (y_2 y_1)^{k-3}y_2,
y_2 \ast (y_1 y_2)^{k-2},
y_1 \ast y_2 y_1^{-1})^{-1} \\ 
&\times&
\prod_{k=1}^{n ({\rm resp.}~ n-1)}
\theta (
y_2 \ast (y_1 y_2)^{k-2},
y_1 \ast (y_2 y_1)^{k-2}y_2,
y_1 \ast y_2 y_1^{-1})^{-1} \\
&\times&
\prod_{k=1}^{~~~~~n~~~~~ }
\theta (
y_1 \ast y_2 y_1^{-1},  
y_2 \ast (y_1 y_2)^{k-2},
y_1 \ast (y_2 y_1)^{k-2}y_2)^{+1} \\
&\times&
\prod_{k=1}^{n ({\rm resp.}~ n-1)}
\theta (
y_1 \ast y_2 y_1^{-1}, 
y_1 \ast (y_2 y_1)^{k-2}y_2,
y_2 \ast (y_1 y_2)^{k-1})^{+1} \\
&\times&
\prod_{k=1}^{n ({\rm resp.}~ n-1)}
\theta (
y_1,
y_2 \ast (y_1 y_2)^{k-1},
y_1 \ast (y_2 y_1)^{k-1}y_2)^{+1} \\
&\times&
\prod_{k=1}^{~~~~n-1~~~~ }
\theta (
y_1,
y_1 \ast (y_2 y_1)^{k-1}y_2,
y_2 \ast (y_1 y_2)^{k})^{+1},
\end{eqnarray*}
where  $y_1, y_2$ run over all elements of  $X$ satisfying
$y_2 \ast (y_1 y_2)^{n} = y_1$
(resp. $y_1 \ast (y_2 y_1)^{n} = y_1$)
and
$y_1 \ast y_2^2 = y_1$.
\end{sect} 
{\it Proof.\/} 
For the Hurwitz arc system $(\alpha_1, \dots, \alpha_6)$
illustrated in Fig.~\ref{taut2malph}, 
the braid system
$$(
w_1^{-1} \sigma_{k_1}^{\epsilon_1} w_1,
w_2^{-1} \sigma_{k_2}^{\epsilon_2} w_2,
\dots,
w_6^{-1} \sigma_{k_6}^{\epsilon_6} w_6) $$
of  $S$ is given by
$$\begin{array}{ll}

w_1 = \sigma_2^{-(m-2)}\sigma_3\sigma_2, \quad &
\sigma_{k_1}^{\epsilon_1} = \sigma_3, \\

w_2 = \sigma_2^{-1}\sigma_3\sigma_2, \quad &
\sigma_{k_2}^{\epsilon_2} = \sigma_1^{-1}, \\

w_3 = \sigma_1^{-(m-1)}\sigma_3\sigma_2, \quad &
\sigma_{k_3}^{\epsilon_3} = \sigma_2, \\

w_4 = \sigma_1^{-1}\sigma_2, \quad &
\sigma_{k_4}^{\epsilon_4} = \sigma_2^{-1}, \\

w_5 = \sigma_2^{-(m-2)}, \quad &
\sigma_{k_5}^{\epsilon_5} = \sigma_1, \\

w_6 = \sigma_2^{-1}, \quad &
\sigma_{k_6}^{\epsilon_6} = \sigma_3^{-1}.

\end{array}$$

\begin{figure}
\begin{center}
\mbox{
\epsfxsize=5in
\epsfbox{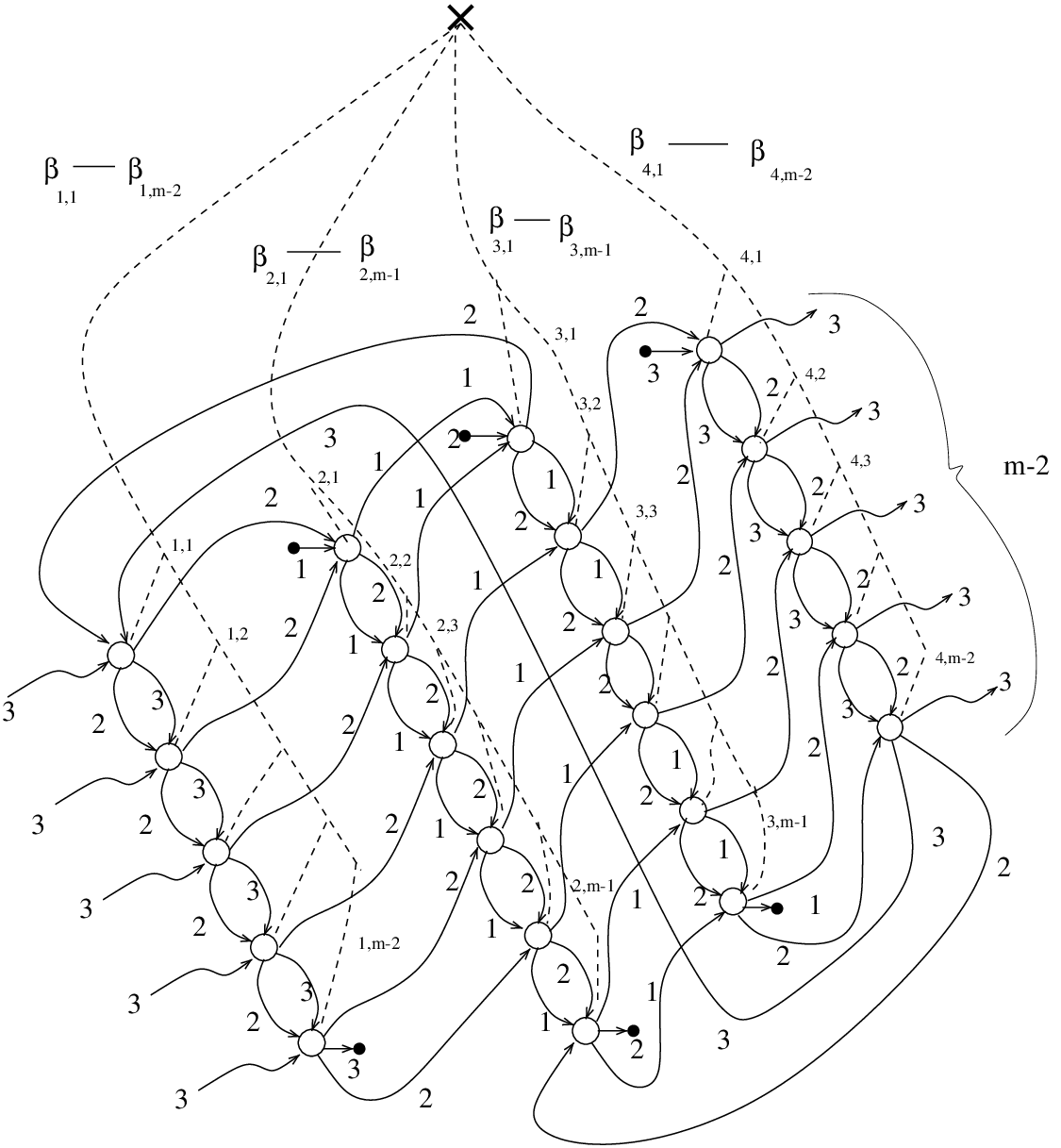}
}
\end{center}
\caption{Paths to the preferred regions near triple points}
\label{taut2mbet}
\end{figure}

For arcs $\{\beta_{1,1}, \dots, \beta_{1,m-2},
\beta_{2,1}, \dots, \beta_{2,m-1},
\beta_{3,1}, \dots, \beta_{3,m-1},
\beta_{4,1}, \dots, \beta_{4,m-2}\}$ as in
Fig.~\ref{taut2mbet}, 
the intersection word $w_{\Gamma}(\beta_{i, j})$
$(i \in\{1,2,3,4\})$ are as follows  (sgn means the sign
of the white vertex $W_{i, j}$ near the starting point of
$\beta_{i, j}$,
and labels mean the labels around $W_{i, j}$):

$$\begin{array}{llll}

w_{\Gamma}(\beta_{1, j}) = \sigma_2^{-(j-1)}\sigma_3\sigma_2,
\quad & {\rm sgn} =-1,  &
{\rm labels} = \{2,3\},  &
(j= 1, \dots, m-2); \\

w_{\Gamma}(\beta_{2, j}) = \sigma_1^{-(j-1)}\sigma_3\sigma_2,
\quad & {\rm sgn} =-1,  &
{\rm labels} = \{1,2\},  &
(j= 1, \dots, m-1); \\

w_{\Gamma}(\beta_{3, j}) = \sigma_2^{-(j-2)},
\quad & {\rm sgn} =+1,  &
{\rm labels} = \{1,2\},  &
(j= 1, \dots, m-1); \\

w_{\Gamma}(\beta_{4, j}) = \sigma_3^{-(j-1)},
\quad & {\rm sgn} =+1,  &
{\rm labels} = \{2,3\},  &
(j= 1, \dots, m-2); \\

\end{array}$$

Suppose that $m$ is odd.  Then
the quandle automorphisms $Q(w_i)$ $(i=1,\dots,6)$ of
$F_Q \langle x_1, \dots, x_4 \rangle$ are as follows:

\begin{eqnarray*}
Q(\sigma_2^{-(m-2)}\sigma_3\sigma_2) &:&
  \left\{
  \begin{array}{ll}
  x_1 &\to x_1, \\
  x_2 &\to x_4 \ast (x_3 x_4)^{(m-3)/2}x_2^{-1}, \\
  x_3 &\to x_3 \ast (x_4 x_3)^{(m-3)/2}x_4 x_2^{-1}, \\
  x_4 &\to x_2,
 \end{array} \right. \\
Q(\sigma_2^{-1}\sigma_3\sigma_2) &:&
  \left\{
  \begin{array}{ll}
  x_1 &\to x_1, \\
  x_2 &\to x_4 \ast x_2^{-1}, \\
  x_3 &\to x_3 \ast x_4 x_2^{-1}, \\
  x_4 &\to x_2,
\end{array} \right. \\
Q(\sigma_1^{-(m-1)}\sigma_3\sigma_2) &:&
  \left\{
  \begin{array}{ll}
  x_1 &\to x_1 \ast (x_2 x_3 x_2^{-1} x_1)^{(m-3)/2} x_2 x_3 x_2^{-1},
\\
  x_2 &\to x_3 \ast (x_2^{-1} x_1 x_2 x_3)^{(m-1)/2} x_2^{-1}, \\
  x_3 &\to x_4 \ast x_2^{-1}, \\
  x_4 &\to x_2,
\end{array} \right. \\
Q(\sigma_1^{-1}\sigma_2) &:&
  \left\{
  \begin{array}{ll}
  x_1 &\to x_3 \ast x_2^{-1}, \\
  x_2 &\to x_1 \ast x_2 x_3 x_2^{-1}, \\
  x_3 &\to x_2, \\
  x_4 &\to x_4,
\end{array} \right. \\
Q(\sigma_2^{-(m-2)}) &:&
  \left\{
  \begin{array}{ll}
  x_1 &\to x_1, \\
  x_2 &\to x_3 \ast (x_2 x_3)^{(m-3)/2}, \\
  x_3 &\to x_2 \ast (x_3 x_2)^{(m-3)/2} x_3, \\
  x_4 &\to x_4,
\end{array} \right. \\
Q(\sigma_2^{-1}) &:&
  \left\{
  \begin{array}{ll}
  x_1 &\to x_1, \\
  x_2 &\to x_3, \\
  x_3 &\to x_2 \ast x_3, \\
  x_4 &\to x_4,
\end{array} \right.
\end{eqnarray*}

Hence the defining relations
$Q(w_i)( x_{k_i} ) = Q(w_i)( x_{k_i +1} )$ ($i=1,\dots, 6$) of
$Q(S)$  are
\begin{eqnarray*}
  x_3 \ast (x_4 x_3)^{(m-3)/2}x_4 x_2^{-1}  &=&
  x_2, \\
  x_1   &=&
  x_4 \ast x_2^{-1}, \\
  x_3 \ast (x_2^{-1} x_1 x_2 x_3)^{(m-1)/2} x_2^{-1} &=&
  x_4 \ast x_2^{-1}, \\
  x_1 \ast x_2 x_3 x_2^{-1}  &=&
  x_2, \\
  x_1 &=&
  x_3 \ast (x_2 x_3)^{(m-3)/2}, \\
  x_2 \ast x_3 &=&
  x_4.
\end{eqnarray*}

Thus the quandle $Q(S)$, for odd $m$, is
\begin{eqnarray*}
\langle x_1, \dots, x_4 |
& & x_3 \ast (x_2 x_3)^{(m-1)/2} = x_2, \\
& & x_2 \ast x_3^2 = x_2, \\
& & x_1 = x_2 \ast x_3 x_2^{-1}, \\
& & x_4 = x_2 \ast x_3
\rangle \\
=
\langle x_2, x_3 |
& & x_3 \ast (x_2 x_3)^{(m-1)/2} = x_2, \\
& & x_2 \ast x_3^2 = x_2
\rangle.
\end{eqnarray*}

This is isomorphic to the dihedral quandle $R_m$.

Suppose that $m$ is even.  Then
the quandle automorphisms $Q(w_i)$ $(i=1,\dots,6)$ of
$F_Q \langle x_1, \dots, x_4 \rangle$ are as follows:

\begin{eqnarray*}
Q(\sigma_2^{-(m-2)}\sigma_3\sigma_2) &:&
  \left\{
  \begin{array}{ll}
  x_1 &\to x_1, \\
  x_2 &\to x_3 \ast (x_4 x_3)^{(m-4)/2}x_4 x_2^{-1}, \\
  x_3 &\to x_4 \ast (x_3 x_4)^{(m-2)/2}x_2^{-1}, \\
  x_4 &\to x_2,
 \end{array} \right. \\
Q(\sigma_2^{-1}\sigma_3\sigma_2) &:&
  \left\{
  \begin{array}{ll}
  x_1 &\to x_1, \\
  x_2 &\to x_4 \ast x_2^{-1}, \\
  x_3 &\to x_3 \ast x_4 x_2^{-1}, \\
  x_4 &\to x_2,
\end{array} \right. \\
Q(\sigma_1^{-(m-1)}\sigma_3\sigma_2) &:&
  \left\{
  \begin{array}{ll}
  x_1 &\to x_3 \ast (x_2^{-1} x_1 x_2 x_3)^{(m-2)/2} x_2^{-1}, \\
  x_2 &\to x_1 \ast (x_2 x_3 x_2^{-1} x_1)^{(m-2)/2} x_2 x_3 x_2^{-1},
\\
  x_3 &\to x_4 \ast x_2^{-1}, \\
  x_4 &\to x_2,
\end{array} \right. \\
Q(\sigma_1^{-1}\sigma_2) &:&
  \left\{
  \begin{array}{ll}
  x_1 &\to x_3 \ast x_2^{-1}, \\
  x_2 &\to x_1 \ast x_2 x_3 x_2^{-1}, \\
  x_3 &\to x_2, \\
  x_4 &\to x_4,
\end{array} \right. \\
Q(\sigma_2^{-(m-2)}) &:&
  \left\{
  \begin{array}{ll}
  x_1 &\to x_1, \\
  x_2 &\to x_2 \ast (x_3 x_2)^{(m-4)/2} x_3, \\
  x_3 &\to x_3 \ast (x_2 x_3)^{(m-2)/2}, \\
  x_4 &\to x_4,
\end{array} \right. \\
Q(\sigma_2^{-1}) &:&
  \left\{
  \begin{array}{ll}
  x_1 &\to x_1, \\
  x_2 &\to x_3, \\
  x_3 &\to x_2 \ast x_3, \\
  x_4 &\to x_4,
\end{array} \right.
\end{eqnarray*}

Hence the defining relations
$Q(w_i)( x_{k_i} ) = Q(w_i)( x_{k_i +1} )$ ($i=1,\dots, 6$) of
$Q(S)$  are
\begin{eqnarray*}
  x_4 \ast (x_3 x_4)^{(m-2)/2}x_2^{-1}  &=&
  x_2, \\
  x_1   &=&
  x_4 \ast x_2^{-1}, \\
  x_1 \ast (x_2 x_3 x_2^{-1} x_1)^{(m-2)/2} x_2 x_3 x_2^{-1} &=&
  x_4 \ast x_2^{-1}, \\
  x_1 \ast x_2 x_3 x_2^{-1}  &=&
  x_2, \\
  x_1 &=&
  x_2 \ast (x_3 x_2)^{(m-4)/2} x_3, \\
  x_2 \ast x_3 &=&
  x_4.
\end{eqnarray*}

Thus the quandle $Q(S)$, for even $m$, is
\begin{eqnarray*}
\langle x_1, \dots, x_4 |
& & x_2 \ast (x_3 x_2)^{m/2} = x_2, \\
& & x_2 \ast x_3^2 = x_2, \\
& & x_1 = x_2 \ast x_3 x_2^{-1}, \\
& & x_4 = x_2 \ast x_3
\rangle \\
=
\langle x_2, x_3 |
& & x_2 \ast (x_3 x_2)^{m/2} = x_2, \\
& & x_2 \ast x_3^2 = x_2
\rangle.
\end{eqnarray*}

Now we calculate the Boltzmann weight for each white
vertex.  Let $\iota: \langle x_1, x_2, x_3, x_4 \rangle
\to Q(S)$ be the natural projection map (that is $i_\ast$
in \cite{CJKLS}).  Let $\theta $ be a quandle 3-cocycle of a
finite quandle, and $c$ a coloring.

\begin{enumerate} \item 
$w_{\Gamma}(\beta_{1, j}) = \sigma_2^{-(j-1)}\sigma_3\sigma_2,
\quad  {\rm sgn} =-1,
{\rm labels} = \{2,3\},
(j= 1, \dots, m-2).$         
\begin{enumerate}
\item
For odd $j$ ($1\leq j \leq m-2$),  the composition of
$Q(\sigma_2^{-(j-1)}\sigma_3\sigma_2)$ and $\iota$ maps
  \begin{eqnarray*}  \left\{  \begin{array}{lll}
  x_1 &\to x_1
      &\to x_2 \ast x_3 x_2^{-1} \\  
  x_2 &\to x_3 \ast (x_4 x_3)^{(j-3)/2}x_4 x_2^{-1}
      &\to x_3 \ast (x_2 x_3)^{(j-1)/2}x_2^{-1} \\  
  x_3 &\to x_4 \ast (x_3 x_4)^{(j-1)/2}x_2^{-1}
      &\to x_2 \ast (x_3 x_2)^{(j+1)/2}x_2^{-2} \\  
  x_4 &\to x_2
      &\to x_2.
 \end{array} \right. \end{eqnarray*}
  Hence the Boltzmann weight is
$$B= \theta (
c(x_3 \ast (x_2 x_3)^{(j-1)/2}x_2^{-1}),
c(x_2 \ast (x_3 x_2)^{(j+1)/2}x_2^{-2}),
c(x_2))^{-1}.$$

  Put $j= 2k-1$, then
$$B= \theta (
c(x_3 \ast (x_2 x_3)^{k-1}x_2^{-1}),
c(x_2 \ast (x_3 x_2)^{k}x_2^{-2}),
c(x_2))^{-1},$$
where $k=1,2,\dots, n$ if $m=2n+1$,
or $k=1,2,\dots, n-1$ if $m=2n$.        

\item
For even $j$ ($1\leq j \leq m-2$),  the composition of
$Q(\sigma_2^{-(j-1)}\sigma_3\sigma_2)$ and $\iota$ maps
  \begin{eqnarray*}  \left\{  \begin{array}{lll}
  x_1 &\to x_1
      &\to x_2 \ast x_3 x_2^{-1} \\ 
  x_2 &\to x_4 \ast (x_3 x_4)^{(j-2)/2}x_2^{-1}
      &\to x_2 \ast (x_3 x_2)^{j/2}x_2^{-1} \\ 
  x_3 &\to x_3 \ast (x_4 x_3)^{(j-2)/2}x_4 x_2^{-1}
      &\to x_3 \ast (x_2 x_3)^{j/2}x_2^{-1} \\ 
  x_4 &\to x_2
      &\to x_2.
 \end{array} \right. \end{eqnarray*}
  Hence the Boltzmann weight is
$$B= \theta (
c(x_2 \ast (x_3 x_2)^{j/2}x_2^{-2}),
c(x_3 \ast (x_2 x_3)^{j/2}x_2^{-1}),
c(x_2))^{-1}.$$

  Put $j= 2k$, then
$$B= \theta (
c(x_2 \ast (x_3 x_2)^{k}x_2^{-2}),
c(x_3 \ast (x_2 x_3)^{k}x_2^{-1}),
c(x_2))^{-1},$$
where $k=1,2,\dots, n-1$.
\end{enumerate}

\item
$w_{\Gamma}(\beta_{2, j}) = \sigma_1^{-(j-1)}\sigma_3\sigma_2,
\quad  {\rm sgn} =-1,
{\rm labels} = \{1,2\},
(j= 1, \dots, m-1).$          \bigskip

\begin{enumerate}
\item
For odd $j$ ($1\leq j \leq m-1$),  the composition of
$Q(\sigma_1^{-(j-1)}\sigma_3\sigma_2)$ and $\iota$ maps
  \begin{eqnarray*}  \left\{  \begin{array}{lll}
  x_1 &\to x_1 \ast (x_2 x_3 x_2^{-1} x_1)^{(j-1)/2}x_1^{-1}
      &\to x_2 \ast (x_3 x_2)^{(j-5)/2}x_3 \\
  x_2 &\to x_3 \ast (x_2^{-1} x_1 x_2 x_3)^{(j-1)/2}x_2^{-1}
      &\to x_3 \ast (x_2 x_3)^{(j-3)/2} \\
  x_3 &\to x_4 \ast x_2^{-1}
      &\to x_2 \ast x_3 x_2^{-1} \\
  x_4 &\to x_2
      &\to x_2.
 \end{array} \right. \end{eqnarray*}
  Hence the Boltzmann weight is
$$B= \theta (
c(x_2 \ast (x_3 x_2)^{(j-5)/2}x_3),
c(x_3 \ast (x_2 x_3)^{(j-3)/2}),
c(x_2 \ast x_3 x_2^{-1}))^{-1}.$$  
  Put $j= 2k-1$, then
$$B= \theta (
c(x_2 \ast (x_3 x_2)^{k-3}x_3),
c(x_3 \ast (x_2 x_3)^{k-2}),
c(x_2 \ast x_3 x_2^{-1}))^{-1},$$  
where $k=1,2,\dots, n$.

\item
For even $j$ ($1\leq j \leq m-1$),  the composition of
$Q(\sigma_1^{-(j-1)}\sigma_3\sigma_2)$ and $\iota$ maps
  \begin{eqnarray*}  \left\{  \begin{array}{lll}
  x_1 &\to x_3 \ast (x_2^{-1} x_1 x_2 x_3)^{(j-2)/2}x_2^{-1}
      &\to x_3 \ast (x_2 x_3)^{(j-4)/2} \\
  x_2 &\to x_1 \ast (x_2 x_3 x_2^{-1} x_1)^{j/2}x_1^{-1}
      &\to x_2 \ast (x_3 x_2)^{(j-4)/2}x_3 \\
  x_3 &\to x_4 \ast x_2^{-1}
      &\to x_2 \ast x_3 x_2^{-1} \\ 
  x_4 &\to x_2
      &\to x_2.
 \end{array} \right. \end{eqnarray*}
  Hence the Boltzmann weight is
$$B= \theta (
c(x_3 \ast (x_2 x_3)^{(j-4)/2}),
c(x_2 \ast (x_3 x_2)^{(j-4)/2}x_3),
c(x_2 \ast x_3 x_2^{-1}))^{-1}.$$ 
  Put $j= 2k$, then
$$B= \theta (
c(x_3 \ast (x_2 x_3)^{k-2}),
c(x_2 \ast (x_3 x_2)^{k-2}x_3),
c(x_2 \ast x_3 x_2^{-1}))^{-1}, $$ 
where $k=1,2,\dots, n$ if $m=2n+1$,
$k=1,2,\dots, n-1$ if $m=2n$.
\end{enumerate}

\item
$w_{\Gamma}(\beta_{3, j}) = \sigma_2^{-(j-2)},
\quad  {\rm sgn} =+1,
{\rm labels} = \{1,2\},
(j= 1, \dots, m-1).$               

\begin{enumerate}\item
For odd $j$ ($1\leq j \leq m-1$),  the composition of
$Q(\sigma_2^{-(j-2)})$ and $\iota$ maps
  \begin{eqnarray*}  \left\{  \begin{array}{lll}
  x_1 &\to x_1
      &\to x_2 \ast x_3 x_2^{-1} \\ 
  x_2 &\to x_3 \ast (x_2 x_3)^{(j-3)/2}
      &\to x_3 \ast (x_2 x_3)^{(j-3)/2} \\
  x_3 &\to x_2 \ast (x_3 x_2)^{(j-3)/2}x_3
      &\to x_2 \ast (x_3 x_2)^{(j-3)/2}x_3 \\
  x_4 &\to x_4
      &\to x_2 \ast x_3.
 \end{array} \right. \end{eqnarray*}
  Hence the Boltzmann weight is
$$B= \theta (
c(x_2 \ast x_3 x_2^{-1}), 
c(x_3 \ast (x_2 x_3)^{(j-3)/2}),
c(x_2 \ast (x_3 x_2)^{(j-3)/2}x_3))^{+1}.$$
  Put $j= 2k-1$, then
$$B= \theta (
c(x_2 \ast x_3 x_2^{-1}), 
c(x_3 \ast (x_2 x_3)^{k-2}),
c(x_2 \ast (x_3 x_2)^{k-2}x_3))^{+1},$$
where $k=1,2,\dots, n$.          \bigskip

\item
For even $j$ ($1\leq j \leq m-1$),  the composition of
$Q(\sigma_2^{-(j-2)})$ and $\iota$ maps
  \begin{eqnarray*}  \left\{  \begin{array}{lll}
  x_1 &\to x_1
      &\to x_2 \ast x_3 x_2^{-1} \\  
  x_2 &\to x_2 \ast (x_3 x_2)^{(j-4)/2}x_3
      &\to x_2 \ast (x_3 x_2)^{(j-4)/2}x_3 \\
  x_3 &\to x_3 \ast (x_2 x_3)^{(j-2)/2}
      &\to x_3 \ast (x_2 x_3)^{(j-2)/2} \\
  x_4 &\to x_4
      &\to x_2 \ast x_3.
 \end{array} \right. \end{eqnarray*}
  Hence the Boltzmann weight is
$$B= \theta (
c(x_2 \ast x_3 x_2^{-1}), 
c(x_2 \ast (x_3 x_2)^{(j-4)/2}x_3),
c(x_3 \ast (x_2 x_3)^{(j-2)/2}))^{+1}.$$
  Put $j= 2k$, then
$$B= \theta (
c(x_2 \ast x_3 x_2^{-1}),  
c(x_2 \ast (x_3 x_2)^{k-2}x_3),
c(x_3 \ast (x_2 x_3)^{k-1}))^{+1}, $$
where $k=1,2,\dots, n$ if $m=2n+1$,
$k=1,2,\dots, n-1$ if $m=2n$.
\end{enumerate}

\item
$w_{\Gamma}(\beta_{4, j}) = \sigma_3^{-(j-1)},
\quad  {\rm sgn} =+1,
{\rm labels} = \{2,3\},
(j= 1, \dots, m-2).$                \bigskip

\begin{enumerate}\item
For odd $j$ ($1\leq j \leq m-2$),  the composition of
$Q(\sigma_3^{-(j-1)})$ and $\iota$ maps
  \begin{eqnarray*}  \left\{  \begin{array}{lll}
  x_1 &\to x_1
      &\to x_2 \ast x_3 x_2^{-1} \\ 
  x_2 &\to x_2
      &\to x_2 \\
  x_3 &\to x_3 \ast (x_4 x_3)^{(j-3)/2}x_4
      &\to x_3 \ast (x_2 x_3)^{(j-1)/2} \\
  x_4 &\to x_4 \ast (x_3 x_4)^{(j-1)/2}
      &\to x_2 \ast (x_3 x_2)^{(j-1)/2}x_3.
 \end{array} \right. \end{eqnarray*}
  Hence the Boltzmann weight is
$$B= \theta (
c(x_2),
c(x_3 \ast (x_2 x_3)^{(j-1)/2}),
c(x_2 \ast (x_3 x_2)^{(j-1)/2}x_3))^{+1}.$$
  Put $j= 2k-1$, then
$$B= \theta (
c(x_2), 
c(x_3 \ast (x_2 x_3)^{k-1}),  
c(x_2 \ast (x_3 x_2)^{k-1}x_3))^{+1},$$  
where $k=1,2,\dots, n$ if $m=2n+1$,
$k=1,2,\dots, n-1$ if $m=2n$.         

\item
For even $j$ ($1\leq j \leq m-2$),  the composition of
$Q(\sigma_3^{-(j-1)})$ and $\iota$ maps
  \begin{eqnarray*}  \left\{  \begin{array}{lll}
  x_1 &\to x_1
      &\to x_2 \ast x_3 x_2^{-1} \\
  x_2 &\to x_2
      &\to x_2  \\
  x_3 &\to x_4 \ast (x_3 x_4)^{(j-2)/2}
      &\to x_2 \ast (x_3 x_2)^{(j-2)/2}x_3 \\
  x_4 &\to x_3 \ast (x_4 x_3)^{(j-2)/2}x_4
      &\to x_3 \ast (x_2 x_3)^{j/2}.
 \end{array} \right. \end{eqnarray*}
  Hence the Boltzmann weight is
$$B= \theta (
c(x_2),
c(x_2 \ast (x_3 x_2)^{(j-2)/2}x_3),
c(x_3 \ast (x_2 x_3)^{j/2}))^{+1}.$$
  Put $j= 2k$, then
$$B= \theta (
c(x_2),
c(x_2 \ast (x_3 x_2)^{k-1}x_3),
c(x_3 \ast (x_2 x_3)^{k}))^{+1},$$
where $k=1,2,\dots, n-1$.
\end{enumerate}\end{enumerate}

By replacing $c(x_2)$ and $c(x_3)$ by $y_1$ and
$y_2$, we have the theorem.
$\Box$

\begin{sect} {\bf Examples.\/} {\rm 
{\bf The case of $m=3$ ($n=1$)}
\begin{eqnarray*}
  \sum_{y_1,y_2}
  &&
  \theta (
y_2 \ast y_1,
y_2 ,
y_1
  )^{-1} \\
  &\times&
  \theta (
y_2 ,
y_2 \ast y_1,
y_2
  )^{-1} \\
  &\times&
  \theta (
y_2 \ast y_1,
y_1 ,
y_2
  )^{-1} \\
  &\times&
  \theta (
y_2,
y_2 \ast y_1 ,
y_1
  )^{+1} \\
  &\times&
  \theta (
y_2,
y_1 ,
y_2
  )^{+1} \\
  &\times&
  \theta (
y_1,
y_2 ,
y_1 \ast y_2
  )^{+1} ,
\end{eqnarray*}
where  $y_1, y_2$ run over all elements of  $X$ satisfying
$y_2 \ast y_1 y_2 = y_1$
and
$y_1 \ast y_2^2 = y_1$.

\bigskip

Put $z_1 = y_1 \ast y_2 = y_2 \ast y_1$ and
$z_2 = y_1$ (and $z_1 \ast z_2 = y_2$), then

\begin{eqnarray*}
  \sum_{y_1,y_2}
  &&
  \theta (
z_1,
z_1 \ast z_2 ,
z_2
  )^{-1} \\
  &\times&
  \theta (
z_1 \ast z_2 ,
z_1,
z_1 \ast z_2
  )^{-1} \\
  &\times&
  \theta (
z_1,
z_2 ,
z_1 \ast z_2
  )^{-1} \\
  &\times&
  \theta (
z_1 \ast z_2,
z_1 ,
z_2
  )^{+1} \\
  &\times&
  \theta (
z_1 \ast z_2,
z_2 ,
z_1 \ast z_2
  )^{+1} \\
  &\times&
  \theta (
z_2,
z_1 \ast z_2 ,
z_1
  )^{+1} ,
\end{eqnarray*}
where  $z_1, z_2$ run over all elements of  $X$ satisfying
$z_2 \ast z_1 z_2 = z_1$
and
$z_1 \ast z_2^2 = z_1$.
This formula is the same as that in \cite{CJKLS}.

{\bf The case of $m=4$ ($n=2$)}
\begin{eqnarray*}
  \sum_{y_1,y_2}
  &&
  \theta (
y_2 \ast y_1^{-1} ,
y_1 \ast y_2 ,
y_1
  )^{-1} \\
  &\times&
  \theta (
y_1 \ast y_2 ,
y_2 \ast y_1 y_2 y_1^{-1} ,  
y_1
  )^{-1} \\
  &\times&
  \theta (
y_1 \ast y_2 ,
y_2 \ast y_1^{-1} ,
y_1 \ast y_2
  )^{-1} \\
  &\times&
  \theta (
y_1 ,
y_2 ,
y_1 \ast y_2
  )^{-1} \\
  &\times&
  \theta (
y_2 \ast y_1^{-1} ,
y_1 ,
y_1 \ast y_2
  )^{-1} \\
  &\times&
  \theta (
y_1 \ast y_2 ,
y_2 \ast y_1^{-1} ,
y_1
  )^{+1} \\
  &\times&
  \theta (
y_1 \ast y_2 ,
y_2 ,
y_1 \ast y_2
  )^{+1} \\
  &\times&
  \theta (
y_1 \ast y_2 ,
y_1 ,
y_2
  )^{+1} \\
  &\times&
  \theta (
y_1,
y_2 ,
y_1 \ast y_2
  )^{+1} \\
  &\times&
  \theta (
y_1,
y_1 \ast y_2,
y_2 \ast y_1 y_2
  )^{+1},
\end{eqnarray*}
where  $y_1, y_2$ run over all elements of  $X$ satisfying
$y_1 \ast (y_2 y_1)^{2} = y_1$
and
$y_1 \ast y_2^2 = y_1$.

Notice that $y_1 \ast y_2 y_1^{-1} = y_1 \ast y_2^{-1} =
y_1 \ast y_2$ and
$y_1 \ast (y_2 y_1)^{-1} = y_1 \ast y_2 y_1$.

{\bf The case of $m=5$ ($n=2$)}
\begin{eqnarray*}
  \sum_{y_1,y_2}
  &&
  \theta (
y_2 \ast y_1,
y_1 \ast y_2 y_1,
y_1
  )^{-1} \\
  &\times&
  \theta (
y_1 \ast y_2 ,
y_2 ,
y_1
  )^{-1} \\
  &\times&
  \theta (
y_1 \ast y_2 y_1 ,
y_1 \ast y_2 ,
y_1
  )^{-1} \\
  &\times&
  \theta (
y_1 \ast y_2 y_1 ,
y_2 \ast y_1 ,
y_1 \ast y_2 y_1
)^{-1} \\
  &\times&
  \theta (
y_1 ,
y_2 ,
y_1 \ast y_2 y_1
)^{-1} \\
  &\times&
  \theta (
y_2 \ast y_1 ,
y_1 ,
y_1 \ast y_2 y_1
  )^{-1} \\
  &\times&
  \theta (
y_2 ,
y_1 \ast y_2 ,
y_1 \ast y_2 y_1
  )^{-1} \\
  &\times&
  \theta (
y_1 \ast y_2 y_1 ,
y_2 \ast y_1 ,
y_1
  )^{+1} \\
  &\times&
  \theta (
y_1 \ast y_2 y_1 ,
y_2 ,
y_1 \ast y_2
  )^{+1} \\
  &\times&
  \theta (
y_1 \ast y_2 y_1,
y_1 ,
y_2
  )^{+1} \\
  &\times&
  \theta (
y_1 \ast y_2 y_1 ,
y_1 \ast  y_2 ,
y_1 \ast y_2 y_1
  )^{+1} \\
  &\times&
  \theta (
y_1 ,
y_2 ,
y_1 \ast y_2
  )^{+1} \\
  &\times&
  \theta (
y_1 ,
y_1 \ast y_2 y_1 ,
y_2 \ast y_1
  )^{+1} \\
  &\times&
  \theta (
y_1 ,
y_1 \ast y_2 ,
y_1 \ast y_2 y_1
  )^{+1},
\end{eqnarray*}
where  $y_1, y_2$ run over all elements of  $X$ satisfying
$y_2 \ast (y_1 y_2)^{2} = y_1$
and
$y_1 \ast y_2^2 = y_1$.

} \end{sect}

\begin{sect}{\bf Remark.\/}
{\rm The values of the partition function
that we computed in this section 
differ from those in the preceding  section. 
For example, using the dihedral quandle $R_3$
 for the 2-twist-spun trefoil would result in $3+6t^2$ here for $(3-2-A)$.
This is because the orientations chosen are opposite.

Even if the same orientation is chosen, the state-sum expressions
(before evaluation with specific cocycles) would be different.
The expressions depend on the choice of diagram, and we used different
diagrams based on different methods. 
The different expressions are related by coboundaries,
that correspond to Roseman moves relating the two diagrams.

}\end{sect}

\section{Symmetry and Cocycle Invariants}
\label{symsect}

In this section, we  discuss how the partition function behaves
 under mirror images and 
change of orientation.

For an element $\sum a_i g_i$ of a group ring
${\bf Z}[A]$ 
 (where $a_i \in {\bf Z}$ and $g_i \in A$), we
denote by
$\overline{\sum a_i g_i}$ the element
$\sum a_i g_i^{-1}$ in the group ring.

For a link $L$, we denote by
$-L$  the same link with the opposite
orientation, by  $L^{\ast}$ the mirror image
of $L$.

\begin{sect}{\bf Theorem.\ }
For any link $L$ and any quandle 2-cocycle
$\phi \in Z^2(Q;A)$,
$$\Phi_\phi ( -L^{\ast}) =              
\overline{\Phi_\phi ( L)}.$$ \end{sect}
 {\it Proof.}  Let $D$ be a link diagram of $L$. We may
assume that the arcs of $D$ around each crossing point
are oriented downward as in Fig.~\ref{twocrossings}. Let $D^{\ast}$ be
the link diagram which is obtained from $D$ by reversing
the vertical direction, and $-D^{\ast}$ the link diagram
obtained from $D^{\ast}$ by reversing the orientation of
the arcs of $D^{\ast}$.  Obviously, the diagram
$-D^{\ast}$ presents the link $-L^{\ast}$.
Each positive (negative, resp.) crossing of $D$,
which looks the left (right) side of Fig.~\ref{twocrossings},
changes to a negative (positive) crossing of
$-D^{\ast}$, which looks
the right (left) side of Fig.~\ref{twocrossings}.
We notice that the labels $x$, $y$, $x \ast y$
around the crossing point of $D$ are inherited
to the corresponding crossing of $-D^{\ast}$.
So the colorings of $D$ by a quandle $Q$ are
naturally in one-to-one correspondence to
the colorings of $-D^{\ast}$, and if the
Boltzmann weight is $\phi(x,y)^{\epsilon}$, then
the corresponding crossing point has
Boltzmann weight $\phi(x,y)^{-\epsilon}$.
Therefore we have $\Phi_\phi ( -L^{\ast}) =
\overline{\Phi_\phi ( L)}$.  
$\Box$

For a surface link $F$, we denote by
$-F$  the same surface link with the opposite
orientation, by  $F^{\ast}$ the mirror image
of $F$.

\begin{sect}{\bf Theorem. \ } 
For any surface link $F$ and any quandle 3-cocycle
$\theta \in Z^3(Q;A)$,
$$\Phi_\theta ( -F^{\ast}) =
\overline{\Phi_\theta ( F)}.$$ 
\end{sect}
{\it Proof.}
The proof is similar to the classical case.
Let $D$ be a broken surface diagram of $F$. We may
assume that every triple point of $D$
looks like one of Fig.~\ref{alltypeIII} in a movie.
Let $D^{\ast}$ be
the diagram which is obtained from $D$ by reversing
the vertical direction in each cross-section of the movie, and
$-D^{\ast}$ the link diagram  obtained from $D^{\ast}$ by
reversing the orientation of  the arcs of each
cross-sectional link diagram of $D^{\ast}$.
Obviously, the diagram
$-D^{\ast}$ presents the surface link $-F^{\ast}$.
Each positive (negative, resp) triple point of $D$,
which looks like one of Fig.~\ref{alltypeIII},
changes to a negative (positive) crossing of
$-D^{\ast}$, which looks
like another of the figure.
For example, a triple point looking like the (1,1)-entry of
Fig.~\ref{alltypeIII} changes to one like the  (1,6)-entry.
A triple point looking like (2,2)-entry of
Fig.~\ref{alltypeIII} changes to one like (5,3)-entry, {\it etc.}
We notice that the labels $p$, $q$, $r$
around the crossing point of $D$ indicated in Fig.~\ref{alltypeIII}
are inherited
to the corresponding crossing of $-D^{\ast}$.
So the colorings of $D$ by a quandle $Q$ are
naturally in one-to-one correspondence to
the colorings of $-D^{\ast}$, and if the
Boltzmann weight is $\theta(p,q,r)^{\epsilon}$, then
the corresponding crossing point has
Boltzmann weight $\theta(p,q,r)^{-\epsilon}$.
Therefore we have $\Phi_\theta ( -F^{\ast}) =
\overline{\Phi_\theta ( F)}$.  $\Box$

\begin{sect}{\bf Example.\ }{\rm The surface link $\tau^2(T(2,m))$ is
isotopic to its mirror image $\tau^2(T(2,m))^{\ast}$.
(This is well-known for the case that $m$ is an odd
integer, and studied in more general cases, cf.
\cite{Litherland}.)  Hence
$$\Phi_\theta (\tau^2(T(2,m))) =
\Phi_\theta (\tau^2(T(2,m))^{\ast}) =
\overline{\Phi_\theta ( -\tau^2(T(2,m)))}
= \overline{\Phi_\theta ( -\tau^2(T(2,m))^{\ast})}.$$

Thus, if we know $\Phi_\theta (\tau^2(T(2,m)))$, we
do not need to calculate the invariants of
$-\tau^2(T(2,m))$, $\tau^2(T(2,m))^{\ast}$
and $-\tau^2(T(2,m))^{\ast}$.
}\end{sect}

\end{document}